\documentclass{article}

\usepackage{amssymb,amscd,amsthm,verbatim,fancyhdr}
\usepackage{mathrsfs}
\usepackage{amsmath}

\begin{document}
\numberwithin{equation}{section}
\newtheorem{thm}{Theorem}
\newtheorem{lemma}{Lemma}
\newtheorem{clm}{Claim}
\newtheorem{remark}{Remark}
\newtheorem{definition}{Definition}
\newtheorem{cor}{Corollary}
\newtheorem{prop}{Proposition}
\newtheorem{statement}{Statement}


\newcommand{\hdt}{{\dot{\mathrm{H}}^{1/2}}}
\newcommand{\hdtr}{{\dot{\mathrm{H}}^{1/2}(\mathbb{R}^3)}}
\newcommand{\R}{\mathbb{R}}
\newcommand{\ei}{\mathrm{e}^{it\Delta}}
\newcommand{\ltrt}{{L^3(\mathbb{R}^3)}}
\newcommand{\ldrd}{{L^d(\mathbb{R}^d)}}
\newcommand{\lprd}{{L^p(\mathbb{R}^d)}}
\newcommand{\lt}{{L^3}}
\newcommand{\ld}{{L^d}}
\newcommand{\lp}{{L^p}}
\newcommand{\rt}{\mathbb{R}^3}
\newcommand{\rd}{\mathbb{R}^d}
\newcommand{\X}{\mathfrak{X}}
\newcommand{\F}{\mathfrak{F}}
\newcommand{\hdhalf}{{\dot H^\frac{1}{2}}}
\newcommand{\hdthalf}{{\dot H^\frac{3}{2}}}
\newcommand{\hdo}{\dot H^1}
\newcommand{\rthmiz}{\R^3\times(-\infty,0)}
\newcommand{\q}[2]{{#1}_{#2}}
\renewcommand{\t}{\theta}
\newcommand{\lxt}[2]{L_{x,\,t}^{#1}}
\newcommand{\rr}{\sqrt{x_1^2+x_2^2}}
\newcommand{\ve}{\varepsilon}
\newcommand{\hdhrt}{\dot H^\frac{1}{2}(\mathbb{R}^3)}
\renewcommand{\P}{\mathbb{P} }
\newcommand{\RR}{\mathcal{R} }
\newcommand{\TT}{\overline{T} }
\newcommand{\e}{\epsilon }
\newcommand{\D}{\Delta }
\renewcommand{\d}{\delta }
\renewcommand{\l}{\lambda }
\newcommand{\To}{\TT_1 }
\newcommand{\ukt}{u^{(KT)} }
\newcommand{\ttil}{\tilde{T_1} }
\newcommand{\etl}{e^{t\Delta} }
\newcommand{\et}{\mathscr{E}_T }
\newcommand{\se}{\mathscr{E}}
\newcommand{\ft}{\mathscr{F}_T }
\newcommand{\eti}{\mathscr{E}^{\infty}_T }
\newcommand{\fti}{\mathscr{F}^{\infty}_T }
\newcommand{\xjn}{x_{j,n}}
\newcommand{\xjpn}{x_{j',n}}
\newcommand{\ljn}{\l_{j,n} }
\newcommand{\ljpn}{\l_{j',n} }
\newcommand{\lkn}{\l_{k,n} }
\newcommand{\voj}{U_{0,j} }
\newcommand{\uoj}{U_{0,j} }
\newcommand{\voo}{U_{0,1} }
\newcommand{\uon}{u_{0,n} }
\newcommand{\N}{\mathbb{N} }
\newcommand{\Z}{\mathbb{Z} }
\newcommand{\E}{\mathscr{E} }
\newcommand{\tu}{\tilde{u} }
\newcommand{\tU}{\tilde{U} }
\newcommand{\etj}{{\E_{T^*_j}}}
\newcommand{\B}{\mathcal{B}}
\newcommand{\tujn}{\tilde{U}_{j,n}}
\newcommand{\tujpn}{\tilde{U}_{j',n}}
\newcommand{\soj}{\sum_{j=1}^J}
\newcommand{\soi}{\sum_{j=1}^\infty}
\newcommand{\hnj}{H_{n,J}}
\newcommand{\enj}{e_{n,J}}
\newcommand{\pnj}{p_{n,J}}
\newcommand{\lfoi}{L^5_{(0,\infty)}}
\newcommand{\lfhoi}{L^{{5/2}}_{(0,\infty)}}
\newcommand{\wnj}{w_n^J}
\newcommand{\rnj}{r_n^J}
\newcommand{\lfij}{{L^5_{I_j}}}
\newcommand{\lfiz}{{L^5_{I_0}}}
\newcommand{\lfio}{{L^5_{I_1}}}
\newcommand{\lfit}{{L^5_{I_2}}}
\newcommand{\lfhij}{{L^{{5/2}}_{I_j}}}
\newcommand{\lfhiz}{{L^{{5/2}}_{I_0}}}
\newcommand{\lfhio}{{L^{{5/2}}_{I_1}}}
\newcommand{\lfhit}{{L^{{5/2}}_{I_2}}}
\newcommand{\lfi}{{L^5_{I}}}
\newcommand{\lfhi}{{L^{{5/2}}_{I}}}
\newcommand{\doh}{D^\frac{1}{2}}
\newcommand{\tjn}{t_{j,n}}
\newcommand{\tjpn}{t_{j',n}}
\newcommand{\wnlj}{w_n^{l,J}}
\renewcommand{\O}{\mathcal{O}}
\newcommand{\bes}{{\dot B^{s_p}_{p,q}}}
\newcommand{\bespp}{{\dot B^{s_p}_{p,p}}}
\newcommand{\besa}{{\dot B^{s_{a}}_{a,b}}}
\newcommand{\besb}{{\dot B^{s_{p}}_{p,q}}}
\newcommand{\besr}{{\dot B^{s_{b}+\frac{2}{\rho}}_{b,r}}}
\newcommand{\tn}{{\tilde \|}}
\newcommand{\tp}{{\tilde{\phi}}}
\newcommand{\xto}{\xrightarrow[n\to\infty]{}}
\newcommand{\LL}{\Lambda}
\newcommand{\binf}{{\dot B^{-{d/p}}_{\infty,\infty}}}
\newcommand{\brq}{{\dot B^{s_{p,r}}_{r,q}}}
\newcommand{\bfin}{{\dot B^{s_{p,r}}_{r,q}}}
\newcommand{\trnj}{{ \tilde{r}_n^J}}
\def\longformule#1#2{
\displaylines{ \qquad{#1} \hfill\cr \hfill {#2} \qquad\cr } }
\def\sumetage#1#2{ \sum_{\binom{\scriptstyle {#1}}{\scriptstyle {#2}}}}
\def\longformule#1#2{ \displaylines{\qquad{#1} \hfill\cr \hfill {#2} \qquad\cr } }
\newcommand{\avint}[1]{{\int_{#1} \!\!\!\!\!\!\!\!\!\! -} \ \ }
\newcommand{\avintttt}[1]{{\int \!\!\!\!\!\!\! - \!\! \int_{#1} \!\!\!\!\!\!\!\!\!\!\!\!\!\!\!\!\!\!\!\!\!\!\!\! - } \ \ }
\newcommand{\avinttt}[1]{{\int \!\!\!\!\!\!\! - \!\! \int_{#1} \!\!\!\!\!\!\!\!\!\!\!\!\!\!\!\! - } \ \ }
\newcommand{\avintt}[1]{{\int \!\!\!\!\!\!\! - \!\! \int_{#1} \!\!\!\!\!\!\!\!\!\! - } \ \ }
\newcommand{\intt}[1]{{\int \!\!\! \int_{#1}}}
\newcommand{\g}{{\gamma}}
\newcommand{\n}{{\nabla}}
\newcommand{\p}{{\partial}}
\newcommand{\esssup}{\mathop{\mathrm{ess\,sup}}}

\title{Partial regularity for Navier-Stokes and liquid crystals inequalities without maximum principle}
\author{{\sc Gabriel S. Koch}\\ \small University of Sussex\\ \small Brighton BN1 9QH, UK \\ \small {\em g.koch@sussex.ac.uk}}

\date{}
\maketitle

\begin{abstract}
In 1985, V. Scheffer discussed partial regularity results for what he called solutions to the ``Navier-Stokes inequality''.  These maps essentially satisfy the incompressibility condition as well as the local and global energy inequalities and the pressure equation which may be derived formally from the Navier-Stokes system of equations, but they are not required to satisfy the Navier-Stokes system itself.

We extend this notion to a system considered by Fang-Hua Lin and Chun Liu in the mid 1990s related to models of the flow of nematic liquid crystals, which include the Navier-Stokes system when the ``director field'' $d$ is taken to be zero.  In addition to an extended Navier-Stokes system, the Lin-Liu model includes a further parabolic system which implies an a priori maximum principle for $d$ which they use to establish partial regularity (specifically, $\mathcal{P}^{1}(\mathcal{S})=0$) of solutions.

For the analogous ``inequality'' one loses this maximum principle, but here we nonetheless establish certain partial regularity results (namely $\mathcal{P}^{\frac 92 + \d}(\mathcal{S})=0$, so that in particular the putative singular set $\mathcal{S}$ has space-time Lebesgue measure zero).  Under an additional assumption on $d$ for any fixed value of a certain parameter $\sigma \in (5,6)$ (which for $\sigma =6$  reduces precisely to the boundedness of $d$ used by Lin and Liu), we obtain the same partial regularity ($\mathcal{P}^{1}(\mathcal{S})=0$) as do Lin and Liu.   In particular, we recover  the partial regularity result ($\mathcal{P}^{1}(\mathcal{S})=0$) of Caffarelli-Kohn-Nirenberg (1982) for ``suitable weak solutions'' of the Navier-Stokes system, and we verify Scheffer's assertion that the same holds for solutions of the weaker ``inequality'' as well.

We remark that the proofs of partial regularity both here and in the work of Lin and Liu largely follow the proof in Caffarelli-Kohn-Nirenberg, which in turn used many ideas from an earlier work of Scheffer (1975).
\end{abstract}

\section{Introduction}
\noindent
In \cite{linliu95} and \cite{linliu}, Fang-Hua Lin and Chun Liu consider the following system, which reduces to the classical Navier-Stokes system in the case $d\equiv 0$ (here we have set various parameters equal to one for simplicity):
\begin{equation}\label{maineq}
\!\!\boxed{\begin{array}{rcl}
u_t  - \Delta u + \nabla^T \cdot [u\otimes u+\n d \odot \n d ]+ \nabla p & = & 0\\\\
\nabla \cdot u & = & 0\\\\
d_t -\Delta d + (u\cdot \nabla)d  +f(d) &=&0
\end{array}}
\end{equation}
with $f=\n F$ for a scalar field $F$ given by $$F(x) := (|x|^2-1)^2\, ,$$
so that
$$f(x)= 4(|x|^2-1)x$$
(and in particular $f(0) = 0$).  We take the spatial dimension to be three, so that for some $\Omega \subseteq \R^3$ and $T>0$, we are considering maps of the form
$$u,d : \Omega \times (0,T) \to \R^3\, , \quad p: \Omega \times (0,T) \to \R\, ,$$
and here
$$F:\R^3 \to \R \, , \quad f:\R^3 \to \R^3$$
are fixed as above.  As usual, $u$ represents the velocity vector field of a fluid, $p$ is the scalar pressure in the fluid, and, as in nematic liquid crystals models, $d$ corresponds roughly\footnote{In principle, for $d$ to only represent a ``direction'' one should have $|d|\equiv 1$. As proposed in  \cite{linliu95},  F(d) is used to model a  Ginzburg-Landau type of relaxation of the pointwise constraint $|d|\equiv 1$.   For further discussions on the modeling assumptions leading to systems such as the one above, see e.g. \cite{linwang14} or the appendix of \cite{linliu95} and the references mentioned therein.} to the ``director field'' representing the local orientation of rod-like molecules, with $u$ also giving the velocities of the centers of mass of those anisotropic molecules.
\\\\
In (\ref{maineq}), for vector fields $v$ and $w$, the matrix fields $v\otimes w$ and $\n v \odot \n w$ are defined to be the ones with entries
$$(v\otimes w)_{ij} = v_i w_j \quad \textrm{and} \quad (\n v \odot \n w)_{ij}=v_{,i}\cdot w_{,j}:= \frac{\partial v_k}{\partial x_i}\frac{\partial w_k}{\partial x_j}$$
(summing over the repeated index $k$ as per the Einstein convention), and for a matrix field $J= (J_{ij})$, we define\footnote{Many authors simply write $\n \cdot J$, which is perhaps more standard.} the vector field $\n^T \cdot J$ by
$$(\n^T \cdot J)_i:=J_{ij,j}:= \frac {\partial J_{ij}}{\partial x_j}$$
(summing again over $j$).  We think formally of $\n$ (as well as any vector field) as a column vector and $\n^T$ as a row vector, so that each entry of (the column vector) $\n^T \cdot J$ is the divergence of the corresponding {\em row} of $J$.  In what follows, for a vector field $v$ we similarly denote by $\n^T v$ the matrix field with $i$-th row given by $\n^T v_i:=(\n v_i)^T$, i.e.,
$$(\n^T v)_{ij}=v_{i,j}:=\frac{\partial v_i}{\partial x_j}\, ,$$
so that for smooth vector fields $v$ and $w$ we always have
\begin{equation}\label{vecprodrule}
\n^T \cdot (v\otimes w) = (\n^T v)w + v(\n \cdot w)= (w\cdot \n)v + v(\n \cdot w)\, .
\end{equation}
For a scalar field $\phi$ we set $\n^2 \phi:= \n^T (\n \phi)$,
and for matrix fields $J=(J_{ij})$ and $K=(K_{ij})$, we let $J:K:=J_{ij}K_{ij}$ (summing over repeated indices) denote the (real) Frobenius inner product of the matrices ($J:K=\mathrm{tr} (J^T K)$).   We set $|J|:=\sqrt{J:J}$ and $|v|:=\sqrt{v\cdot v}$, and to minimize cumbersome notation will often abbreviate by writing $\n v:=\n^T v$ for a vector field $v$ where the precise structure of the {\em matrix} field $\n^T v$ is not crucial; for example, $|\n v|:=|\n^T v|$.
\\\\
We note that by formally taking the divergence $\n \cdot$ of the first line in (\ref{maineq}) we obtain the usual ``pressure equation''
\begin{equation}\label{preseqa}
-\D p = \n \cdot (\n^T \cdot [u\otimes u+\n d \odot \n d ])\, .
\end{equation}
As in the Navier-Stokes ($d\equiv 0$) setting, one may formally deduce (see Section \ref{motivation} for more details) from (\ref{maineq}) the following global and local energy inequalities which one may expect ``sufficiently nice'' solutions of (\ref{maineq}) to satisfy:\footnote{For sufficiently regular solutions one can show that equality holds. }
\begin{equation}\label{globenineq}
\displaystyle{\frac{d}{d t}\int_{\Omega} \left[ \frac{|u|^2}2 +\frac{|\n d|^2}2 + F(d)  \right]\, dx + \int_{\Omega}\left[|\n u|^2 + |\D d-f(d)|^2 \right]}\, dx\leq 0
\end{equation}
for each $t\in (0,T)$, as well as a localized version\footnote{Note that in \cite{linliu}, the term ``$-\mathcal{R}_f(d,\phi)$'' in (\ref{locenineq}) actually appears incorrectly as ``$+\mathcal{R}_f(d,\phi)$''.  See Section \ref{motivation} for more details.\label{minussign}}
\begin{equation}\label{locenineq}
\begin{array}{l}
\displaystyle{\frac{d}{d t}\int_{\Omega}\left[\left(\frac{|u|^2}{2} + \frac{|\n d|^2}2 \right)\phi \right] \, dx + \int_{\Omega}\left(|\n u|^2+ |\n^2 d|^2\right)\phi \, dx}
\\\\
\qquad \qquad \qquad
\begin{array}{l}
\displaystyle{\leq  \int_{\Omega} \bigg[\left(\frac{|u|^2}{2}+ \frac{|\n d|^2}{2}\right)(\phi_t + \D \phi) +\left(\frac{|u|^2}2 +\frac{|\n d|^2}2 + p\right)u\cdot \n \phi} \\\\
\displaystyle{ \qquad  \qquad  \qquad \qquad \qquad +\ \  u\otimes \n \phi : \n d \odot \n d  \ \  -
\ \ \underbrace{\phi\n^T [f(d)]:  \n^T d}_{=:\mathcal{R}_f(d,\phi)}\bigg]\, dx}
\end{array}
\end{array}
\end{equation}
for $t\in (0,T)$ and each smooth, compactly supported in $\Omega$ and non-negative scalar field $\phi \geq 0$.  (For Navier-Stokes, i.e. when $d\equiv 0$, one may omit all terms involving $d$, even though $0\neq F(0)\notin L^1(\R^3)$.)
\\\\
In \cite{linliu95}, for smooth and bounded $\Omega$, the global energy inequality (\ref{globenineq}) is used to construct global weak solutions to (\ref{maineq}) for initial velocity in $L^2(\Omega)$, along with a similarly appropriate condition on the initial value of $d$ which allows (\ref{globenineq}) to be integrated over $0<t<T$.  This is consistent with the pioneering result of J. Leray \cite{leray} for Navier-Stokes (treated later by many other authors using various methods, but always relying on the natural energy as in \cite{leray}).
\\\\
In \cite{linliu}, the authors establish a partial regularity result for weak solutions to (\ref{maineq}) belonging to the natural energy spaces which moreover satisfy the local energy inequality (\ref{locenineq}).  The result is of the same type as known partial regularity results for ``suitable weak solutions'' to the Navier-Stokes equations.  The program for such partial regularity results for Navier-Stokes was initiated in a series of papers by V. Scheffer in the 1970s and 1980s (see, e.g., \cite{scheffer77,scheffer80} and other works mentioned in \cite{caf}), and subsequently improved by various authors (e.g. \cite{caf,lin,ladyser,vasseur}), perhaps most notably by L. Caffarelli, R. Kohn and L. Nirenberg in \cite{caf}. They show (as do \cite{linliu}) that the one-dimensional parabolic Hausdorff measure of the (potentially empty) singular set $S$ is zero ($\mathcal{P}^1(S)=0$, see Definition \ref{phausdorff} below), implying that singularities (if they exist) cannot for example form any smooth one-parameter curve in space-time.  The method of proof in \cite{linliu} largely follows the method of \cite{caf}.
\\\\
Of course the general system (\ref{maineq}) is (when $d \neq 0$) substantially more complex than the Navier-Stokes system, and one therefore could not expect a stronger result than the type in \cite{caf}.  In fact, it is surprising that one even obtains the same type of result ($\mathcal{P}^1(S)=0$) as in \cite{caf}.  The explanation for this seems to be that although (\ref{maineq}) is more complex than Navier-Stokes in view of the additional $d$ components, one can derive an a priori maximum principle for $d$ because of the third equation in (\ref{maineq}) which substantially offsets this complexity from the viewpoint of regularity.  Therefore, under suitable boundary and initial conditions on $d$, one may assume that $d$ is in fact bounded, a fact which is significantly exploited in \cite{linliu}.  More recently, the authors of the preprint \cite{duhuwang} establish the same type of result for a related but more complex ``Q-tensor'' system; however there, as well, one may obtain a maximum principle which is of crucial importance for proving partial regularity.  One is therefore led to the following natural question, which we will address below:
\\\\
{\bf Can one deduce any partial regularity for systems similar in structure to (\ref{maineq}) but which lack any maximum principle?}
\\\\
In the Navier-Stokes setting, it was asserted by Scheffer in \cite{scheffer3} that in fact the proof of the partial regularity result in \cite{caf} does not require the full set of equations in (\ref{maineq}).  He mentions that the key ingredients are membership of the global energy spaces, the local energy inequality (\ref{locenineq}), the divergence-free condition $\n \cdot u =0$ and the {\em pressure} equation (\ref{preseqa}) (with $d\equiv 0$ throughout).  Scheffer called vector fields satisfying these four requirements solutions to the ``Navier-Stokes inequality'', equivalent to solutions to the Navier-Stokes equations with a forcing $f$ which satisfies $f\cdot u \leq 0$ everywhere.  In contrast, the results in \cite{linliu} do very strongly use the third equation in (\ref{maineq}) in that it implies a maximum principle for $d$.
\\\\
In this paper, we explore what happens if one considers the analog of Scheffer's ``Navier-Stokes inequality''  for the system (\ref{maineq}) when $d\neq 0$.  That is, we consider triples $(u,d,p)$ with global regularities implied (at least when $\Omega$ is bounded and under suitable assumptions on the initial data) by (\ref{globenineq}) which satisfy (\ref{preseqa}) and $\n \cdot u=0$ weakly as well as (a formal consequence of) (\ref{locenineq}), but are {\em not} necessarily weak solutions of the first and third equations (i.e., the two vector equations) in (\ref{maineq}).  In particular, we will {\em not} assume that $d\in L^\infty(\Omega \times (0,T))$, which would have been reasonable in view of the third equation in (\ref{maineq}).  We see that without further assumptions, the result is substantially weaker than the $\mathcal{P}^1(S)=0$ result for Navier-Stokes: following the methods of \cite{linliu,caf} we  obtain (see Theorem \ref{mainthm} below) $\mathcal{P}^{\frac 92 + \delta}(S)=0$ for any $\delta >0$.  This reinforces our intuition that the situation here is substantially more complex than that of Navier-Stokes.  On the other hand, we show that under a suitable uniform local decay condition on $|d|^\sigma (|u|^3 +|\n d|^3)^{(1-\frac \sigma 6)}$ with $\sigma \in (5,6)$
(see (\ref{dmorreysmall}) below, which in particular holds when $d\equiv 0$ as in \cite{caf}), one in fact obtains $\mathcal{P}^1(S)=0$ as in \cite{linliu} and \cite{caf}.  In particular, we verify the above-mentioned assertion made by Scheffer in \cite{scheffer3} regarding partial regularity for Navier-Stokes inequalities.
\\\\
Our key observation  which allows us to work without any maximum principle is that,  in view of the global energy (\ref{globenineq}) and the particular forms of $F$ and $f$, it is reasonable (see Section \ref{motivation}) to assume (\ref{enspaces}); this implies\footnote{In fact, one can also show that $d\in L^s_{\mathrm{loc}}(0,T;L^\infty(\Omega))$ for any $s\in [2,4)$.} that ${d \in L^\infty(0,T;L^6(\Omega))}$ which is sufficient for our purposes.
\\\\
As alluded to above, for our purposes we actually do not require all of the information which appears in (\ref{locenineq}) above.  In view of the fact that
\begin{equation}\label{expandrfdphipre}
\left| \mathcal{R}_f(d,\phi) \right| =|\phi\n^T [f(d)]:  \n^T d| \leq 12|d|^2|\n d|^2\phi + 8\left(\frac{|\n d|^2}2 \phi \right)
\end{equation}
(see (\ref{expandrfdphi}) below), a consequence of (\ref{locenineq}) is that
\begin{equation}\label{locenineqrougher}
\mathcal{A}'(t) + \mathcal{B}(t) \leq 8\mathcal{A}(t) + \mathcal{C}(t)\ \quad \textrm{for} \ \ 0<t<T\, ,
\end{equation}
with $\mathcal{A}, \mathcal{B}, \mathcal{C} \geq 0$ defined (denoting $\int_{\Omega \times \{t\}}g:=\int_{\Omega} g(\cdot, t)\, dx$) as
$$\mathcal{A}(t):=\int_{\Omega \times \{t\}}\left(\frac{|u|^2}{2} + \frac{|\n d|^2}2 \right)\phi  \, , \quad
\mathcal{B}(t):=\int_{\Omega \times \{t\}}\left(|\n u|^2+ |\n^2 d|^2\right)\phi
$$
and
$$\mathcal{C}(t):=\int_{\Omega \times \{t\}} \bigg[\left(\frac{|u|^2}{2}+ \frac{|\n d|^2}{2}\right)|\phi_t + \D \phi| + 12 |d|^2 |\n d|^2\phi\bigg]
 \qquad \qquad \qquad$$
 $$\qquad \qquad \qquad +\ \ \left|\int_{\Omega \times \{t\}} \left[\left(\frac{|u|^2}2 +\frac{|\n d|^2}2 + p\right)u\cdot \n \phi + u\otimes \n \phi : \n d \odot \n d\right]\right|\, .
$$
(\ref{locenineqrougher}) is nearly sufficient, with the appearance of $\mathcal{A}(t)$   on the right-hand side  (in fact, even with $u$ omitted,  which cannot be avoided as ``$\mathcal{R}_f(d,\phi)$'' appears on the right-hand side of (\ref{locenineq}) with a minus\footnote{See Footnote \ref{minussign}.} sign) actually being,  for technical reasons,  the only\footnote{In fact, the appearance of $|d|^2$ on the right-hand side of (\ref{expandrfdphipre}), and hence of (\ref{locenineqrougher}) as well, is handled precisely by the assumption that ${d \in L^\infty(0,T;L^6(\Omega))}$, and is the reason for the slightly weaker results compared to the Navier-Stokes setting (i.e., when $d\equiv 0$).} troublesome term.\footnote{Note that if  $\mathcal{R}_f(d,\phi)$ had appeared with a plus sign in (\ref{locenineq}), one could have simply dropped this troublesome term as a non-positive quantity.}  We therefore use a Gr\"onwall-type argument to hide this term to the left-hand side of (\ref{locenineqrougher}) so that (if $\phi|_{t=0} \equiv 0$)
\begin{equation}\label{locenineqsuf}
\mathcal{A}'(t) + \mathcal{B}(t) \leq  \mathcal{C}(t) + 8e^{8T}\int_0^t\mathcal{C}(\tau)\, d\tau\ \quad \textrm{for} \ \ 0<t<T\, .
\end{equation}
The (formally derived) local energy inequality (\ref{locenineqsuf})  implies (\ref{locenta}) below (for an appropriate constant ${\bar C} \sim 8Te^{8T} + 1$), which is sufficient for our purposes. (In fact, for all elements of the proof other than Proposition \ref{lc}, a weaker form as in (\ref{locent}) is sufficient.)
\ \\\\
Our main result is the following:

\begin{thm}\label{mainthm}
Fix any open set $\Omega \subset \R^3$ and any $T, {\bar C}\in (0,\infty)$, set $\Omega_T:=\Omega \times (0,T)$ and suppose
$u,d:\Omega_T  \to \R^3$ and $p: \Omega_T \to \R$ satisfy the following four assumptions:
\begin{enumerate}
\item $u$, $d$ and $p$  belong to the following spaces:\footnote{For a vector field $f$ or matrix field $J$ and scalar function space $X$, by $f\in X$ or $J\in X$ we mean that all components or entries of $f$ or $J$ belong to $X$; by $\n^2 f \in X$ we mean all second partial derivatives of all components of $f$ belong to $X$; etc.}
\begin{equation}\label{enspaces}
u, d, \n d \in L^\infty(0,T;L^2(\Omega))\, , \quad
\nabla u, \n d, \n^2 d  \in L^2(\Omega_T)
\end{equation}
and
\begin{equation}\label{pspace}
p\in L^{\frac 32}(\Omega_T)\, ;
\end{equation}
\item $u$ is weakly divergence-free:\footnote{Locally integrable functions will always be associated to the standard distribution whose action is integration against a suitable test function so that, e.g., $[\n \cdot u](\psi)= -[u](\n \psi):= -\int u\cdot \n \psi$ for $\psi \in \mathcal{D}(\Omega_T)$.}
\begin{equation}\label{divfree}
\nabla \cdot u = 0 \quad \textrm{in} \quad \mathcal{D}'(\Omega_T)\, ;
\end{equation}
\item the following pressure equation holds weakly:\footnote{Note that $u\otimes u + \n d \odot \n d \in L^{\frac 53}(\Omega_T) \subset L^1_{\mathrm{loc}}(\Omega_T)$, see (\ref{prespacetimeinterpaaa}) - (\ref{graddmixedleb}).}
\begin{equation}\label{preseq}
-\Delta p = \nabla \cdot [\nabla^T \cdot (u\otimes u+\n d \odot \n d)]\quad \textrm{in} \quad \mathcal{D}'(\Omega_T)\, ;
\end{equation}
\item the following  local energy inequality holds:\footnote{For brevity, for $\omega \subset \R^3$, we set  $$\int_{\omega \times \{t\}} g\, dx:=\int_{\omega}g(x, t)\, dx\, .$$ \label{fnintdef}}
\begin{equation}\label{locenta}
 \boxed{\begin{array}{l}
 \int_{\Omega \times \{t\}} \left(|u|^2 + |\n d|^2\right) \phi \, dx + \int_{0}^t \int_{\Omega} \left(|\n u|^2+ |\n^2 d|^2\right) \phi \, dx\, d\tau  \\\\
\qquad \leq {\bar C}\int_{0}^t \big\{ \int_{\Omega \times \{\tau\}} \left[\left(|u|^2 + |\n d|^2\right)|\phi_t + \D \phi|  + |d|^2 |\n d|^2\phi \right]\, dx\\\\
 \qquad \qquad \quad +\  \  \big|\int_{\Omega \times \{\tau\}} \big[\big(\frac{|u|^2}2 +\frac{|\n d|^2}2 + p\big)u\cdot \n \phi + u\otimes \n \phi : \n d \odot \n d\big]\, dx\big|\ \big\}\, d\tau\ \\\\
\qquad \textrm{for}\ \textrm{a.e.}\ t\in (0,T)\  \quad \textrm{and}\quad \forall \ \phi \in  \mathcal{C}_0^\infty(\Omega \times (0,\infty))\ \textrm{s.t.}\ \phi \geq 0\, .
\end{array}}
\end{equation}
\end{enumerate}
Let $\mathcal{S} \subset \Omega_T$ be the (potentially empty) set of singular points where $|u|$ and $|\n d|$ are not essentially bounded in any neighborhood of each $z\in \mathcal{S}$, and let $\mathcal{P}^k$ be the $k$-dimensional parabolic Hausdorff outer measure (see Definition \ref{phausdorff} below).  The following are then true:
\begin{enumerate}
\item $\mathcal{P}^{\frac 92 + \d}(\mathcal{S})=0$, for any $\d >0$ arbitrarily small.
\item If \footnote{In general we set $z=(x,t)\in \Omega_T$, $dz:=dx\, dt$, and recall from Definition \ref{phausdorff} that $Q_r(x_0,t_0):=B_r(x_0)\times (t_0-r^2,t_0)$.}

\begin{equation}\label{dmorreysmall}
g_\sigma:=\sup_{z_0\in \Omega_T}\left( \limsup_{r\searrow 0} \frac 1{r^{2+\frac \sigma2}}
 \intt{Q_r(z_0)}|d|^\sigma(|u|^3 +|\n d|^3)^{(1-\frac \sigma 6)}\, dz\right)< \infty
\end{equation}
for some $\sigma \in (5,6)$, then $\mathcal{P}^1(\mathcal{S})=0$.
\end{enumerate}
\end{thm}
\ \\
Note that in the case $d\equiv 0$, we regain the classical result of $\mathcal{P}^1(\mathcal{S})=0$ for Navier-Stokes as obtained in, for example, \cite{caf}, and more specifically for the (weaker) Navier-Stokes inequalities mentioned in \cite{scheffer3}.
\\\\
We recall that the definition of the outer parabolic Hausdorff measure $\mathcal{P}^k$ is given as follows (see \cite[pp.783-784]{caf}):

\begin{definition}[Parabolic Hausdorff measure]\label{phausdorff}
For any $\mathcal{S} \subset \R^3 \times \R$ and $k\geq 0$, define
$$\mathcal{P}^k(\mathcal{S}):=\lim_{\delta \searrow 0}\mathcal{P}^k_\delta(\mathcal{S})\, ,$$
where
$$\mathcal{P}^k_\delta(\mathcal{S}):=\inf\left\{\, \sum_{j=1}^\infty r_j^k\ \bigg| \ \mathcal{S}\subset \bigcup_{j=1}^\infty Q_{r_j}\, , r_j < \delta\  \forall j\in \N\, \right\}$$
and $Q_r$ is any parabolic cylinder of radius $r>0$, i.e. $$Q_r=Q_r(x,t):= B_r(x) \times (t-r^2,t) \subset \R^3 \times \R$$ for some $x\in \R^3$ and $t\in \R$.  $\mathcal{P}^k$ is an outer measure, and all Borel sets are $\mathcal{P}^k$-measurable.
\end{definition}
\noindent

\begin{remark}\label{pspacesremark}
In the case $\Omega = \rt$, the condition (\ref{pspace}) on the pressure follows (locally, at least) from (\ref{enspaces}) and (\ref{preseq}) if $p$ is taken to be the potential-theoretic solution to (\ref{preseq}), since (\ref{enspaces}) implies that  ${u,\n d\in L^{\frac{10}3}(\Omega_T)}$ by interpolation (see  (\ref{prespacetimeinterpaaa})) and Sobolev embeddings, and then (\ref{preseq}) gives \linebreak $p\in L^{\frac{5}3}(\Omega_T)\subset L^{\frac{3}2}_{\mathrm{loc}}(\Omega_T)$ by Calderon-Zygmund estimates.  For a more general $\Omega$, the existence of such a $p$ can be derived from the motivating equation (\ref{maineq})  (e.g. by estimates for the Stokes operator), see \cite{linliu} and the references therein.  Here, however, we will not refer to (\ref{maineq}) at all and simply {\em assume} $p$ satisfies (\ref{pspace}) and address the partial regularity of such a hypothetical set of functions satisfying \linebreak (\ref{enspaces}) - (\ref{locenta}).
\end{remark}
\ \\
We note that  Theorem \ref{mainthm} does not immediately recover the result of \cite{linliu} (which would correspond to $\sigma = 6$ in (\ref{dmorreysmall}), which holds when $d\in L^\infty$ as assumed in \cite{linliu}).  Heuristically, however, one can argue\footnote{We assume this is roughly the argument in \cite{linliu}, although the details are not explicitly given; see, in particular,  \cite[(2.45)]{linliu} which appears without the ``remainder'' term denoted in \cite{linliu} by $\mathbf{R}(f,\phi)$, and here by $\mathcal{R}_f(d,\phi)$. } as follows:
\\\\
If $d$ were bounded, then taking for example $D:=24 \|d\|_{L^\infty(\Omega_T)}^2 +8 <\infty$ one would deduce from (\ref{expandrfdphipre}) that
$$|\mathcal{R}_f(d,\phi)| \leq  D\left( \frac{|\n d|^2}2 \right)\phi\, .$$
Adjusting the Gr\"onwall-type argument leading to (\ref{locenineqsuf}),  one could then deduce from (\ref{locenineq}) that (if $\mathcal{A}(0)=0$)
$$\mathcal{A}'(t) + \mathcal{B}(t) \leq  \widetilde{\mathcal{C}}(t) + De^{DT}\int_0^T\widetilde{\mathcal{C}}(\tau)\, d\tau\ \quad \textrm{for} \ \ 0<t<T\, ,
$$
where
$$\widetilde{\mathcal{C}}(t):=\int_{\Omega \times \{t\}} \left(\frac{|u|^2}{2}+ \frac{|\n d|^2}{2}\right)|\phi_t + \D \phi|
 +\ \ \left|\int_{\Omega \times \{t\}} \left[\left(\frac{|u|^2}2 +\frac{|\n d|^2}2 + p\right)u\cdot \n \phi + u\otimes \n \phi : \n d \odot \n d\right]\right|\, .
$$
Using such an energy inequality, one would  not need to include the $|d|^6$ term in $E_{3,6}$ (see (\ref{etdefn})) as one would not need to consider the term coming from $\mathcal{R}_f(d,\phi)$  at all in Proposition \ref{lb}, and (noting that the $L^\infty$ norm is invariant under the re-scaling on $d$ in (\ref{recenteredscaling})) one could then adjust Lemmas \ref{thma} and \ref{thmb} appropriately to recover the result in \cite{linliu} using the proof of Theorem \ref{mainthm} below.
\\\\
Finally, we remark that the majority of the arguments in the proofs given below are not new, with many essentially appearing in \cite{linliu} or \cite{caf}. However we feel that our presentation is particularly transparent and may be a helpful addition to the literature, and we include all details so that our results are easily verifiable.
\\\\
{\bf Acknowledgment:} \quad The author would like to offer his sincere thanks to Prof. Arghir Zarnescu for many insightful discussions, for introducing him to the field of liquid crystals models, and for suggesting a problem which led to this publication.  The author would also like to thank Prof. Camillo De Lellis for introducing him to Scheffer's notion of Navier-Stokes inequalities. Finally, the author would like to thank the anonymous referee for insightful comments about a previous draft of this article.

\section{Motivation}\label{motivation}
\noindent
We will show in this section that the assumptions in Theorem \ref{mainthm} are at least formally satisfied by smooth solutions to the system (\ref{maineq}).

\subsection{Energy identities}
As in \cite{linliu}, let us assume that we have smooth solutions to (\ref{maineq}) which vanish or decay sufficiently at $\partial \Omega$ (assumed smooth, if non-empty) and at spatial infinity as appropriate so that all boundary terms vanish in the following integrations by parts, and proceed to establish smooth versions of (\ref{globenineq}) and (\ref{locenineq}).  First, noting the simple identities
\begin{equation}\label{llloca}
\n^T \cdot (\n d \odot \n d)=
  \n \left(\frac{|\n d|^2}2\right)+
(\n^T d)^T\D d
\end{equation}
and
\begin{equation}\label{symmetryetc}
[(\n^T d)^T\D d]\cdot u
=[(\n^T d)u]\cdot \D d
=[(u\cdot \n)d]\cdot \D d\, ,
\end{equation}
at a fixed $t$ one may perform various integrations by parts (keeping in mind that $\n \cdot u =0$) to see that
\begin{equation}\label{preglobllenueq}
\begin{array}{rcl}
0& =&\displaystyle{ \int_{\Omega}[u_t - \D u + \n^T \cdot (u\otimes u) + \n p + \n^T \cdot (\n d \odot \n d)]\cdot u \, dx}\\\\
& =& \displaystyle{ \int_{\Omega}\left[\frac{\p}{\p t} \left(\frac{|u|^2}2 \right) + |\n u|^2 + \underbrace{[(u\cdot \n)d]\cdot \D d}\right]\, dx }\\\\
\end{array}
\end{equation}
and, recalling that $f = \nabla F$ so that $[d_t + (u\cdot \n)d]\cdot f(d) = \left(\tfrac{\partial}{\partial t} + u\cdot \n\right)[F(d)]$, that
\begin{equation}\label{preglobllendeq}
\!\!\!\!\!\!\!\!\begin{array}{rcl}
0& =& -\displaystyle{ \int_{\Omega}[d_t + (u\cdot \n)d - (\D d - f(d))]\cdot (\D d-f(d))\, dx}\\\\
&=& \displaystyle{-\int_{\Omega}\left[-\frac{\p}{\p t}\left(\frac{|\n d|^2}2 + F(d)\right) + \underbrace{[(u\cdot \n)d]\cdot \D d} - |\D d-f(d)|^2\right] \, dx\, . }
\end{array}
\end{equation}
Adding the two gives the
\\\\
\underline{Global energy identity for (\ref{maineq}):}
\begin{equation}\label{globllen}
\displaystyle{\frac{d}{d t}\int_{\Omega} \left[ \frac{|u|^2}2 +\frac{|\n d|^2}2 + F(d)  \right]\, dx + \int_{\Omega}\left[|\n u|^2 + |\D d-f(d)|^2 \right]}\, dx=0
\end{equation}
in view of the cancelation of the indicated terms in (\ref{preglobllenueq}) and (\ref{preglobllendeq}).
\\\\
It is not quite straightforward to localize the calculations in (\ref{preglobllenueq}) and (\ref{preglobllendeq}), for example replacing the (global) multiplicative factor $(\D d - f(d))$ by $(\D d - f(d))\phi$ for a smooth and compactly supported $\phi$.  Arguing as in \cite{linliu}, one can deduce a local energy identity by instead replacing $(\D d - f(d))$ by only a part of its localized version in divergence-form, namely by $\n^T \cdot(\phi \n^T d)$, at the expense of the appearance of $|\D d-f(d)|^2$ anywhere in the local energy.
\\\\
Recalling (\ref{llloca}) and (\ref{symmetryetc}) and noting further that
$$
\!\!\!\!\!\!\!\!\!\!\begin{array}{rcl}
[(u\cdot \n)d] \cdot [\n^T \cdot (\phi \n^T d)] & = &
\displaystyle{
[(u\cdot \n)d] \cdot [\phi \D d] + [(u\cdot \n)d] \cdot [(\n\phi \cdot \n)d]
 }\\\\
& = &
\displaystyle{
[(u\cdot \n)d] \cdot [\phi \D d] + u\otimes \n \phi : \n d \odot \n d
 }
\end{array}
$$
and that
$$[\Delta(\n^Td)]:\n^Td = \Delta\left(\frac{|\n d|^2}2\right) - |\n^2 d|^2\, ,$$
one may  perform various integrations by parts to deduce (as $\n \cdot u=0$) that
$$
\begin{array}{rcl}
0& =& \displaystyle{ \int_{\Omega}[u_t - \D u + \n^T \cdot (u\otimes u) + \n p + \n^T \cdot (\n d \odot \n d)]\cdot u\phi\, dx }\\\\
& =& \displaystyle{ \int_{\Omega}\left[\frac{\p}{\p t}\left(\frac{|u|^2}{2}\phi\right) + |\n u|^2\phi -  \frac{|u|^2}{2}(\phi_t + \D \phi) \right .}\\\\
& & \displaystyle{ \qquad \left . -\left(\frac{|u|^2}2 +\frac{ |\n d|^2}2+ p\right)u\cdot \n \phi + \underbrace{[(u\cdot \n)d]\cdot (\D d)\phi} \right]\, dx}
\end{array}
$$
and
$$
\begin{array}{rcl}
0& =&  -\displaystyle{ \int_{\Omega}[d_t + (u\cdot \n)d - (\D d - f(d))]\cdot [\n^T \cdot (\phi \n^T d)]\, dx}\\\\
& =& \displaystyle{ -\int_{\Omega}\left[   - \frac{\p}{\p t}\left(\frac{|\n d|^2}2 \phi\right) - |\n^2 d|^2\phi + \frac{|\n d|^2}2 (\phi_t +  \D \phi)  \right .}\\\\
& & \displaystyle{ \qquad \qquad -
\n^T [f(d)]: \phi \n^T d + \underbrace{[(u\cdot \n)d] \cdot (\D d)\phi} + u\otimes \n \phi : \n d \odot \n d  \bigg]\, dx}
\end{array}
$$
for smooth and compactly-supported $\phi$, upon adding which and noting again the cancelation of the indicated terms we obtain the
\\\\
\underline{Local energy identity for (\ref{maineq}):}
\begin{equation}\label{locllen}\frac{d}{d t}\int_{\Omega}\left[\left(\frac{|u|^2}{2} + \frac{|\n d|^2}2 \right)\phi \right] \, dx + \int_{\Omega}\left(|\n u|^2+ |\n^2 d|^2\right)\phi \, dx=\qquad \qquad
\end{equation}
$$ \qquad =  \int_{\Omega} \bigg[\left(\frac{|u|^2}{2}+ \frac{|\n d|^2}{2}\right)(\phi_t + \D \phi) +\left(\frac{|u|^2}2 +\frac{|\n d|^2}2 + p\right)u\cdot \n \phi
$$
\ \medskip
$$\qquad  \qquad +\ \  u\otimes \n \phi : \n d \odot \n d  \ \  -
\ \ \underbrace{\phi\n^T [f(d)]:  \n^T d}_{=:\mathcal{R}_f(d,\phi)}\bigg]\, dx \, .$$
Note that we have corrected the omission in \cite{linliu} of the ``$-$'' preceding $\mathcal{R}_f(d,\phi)$, and  the term ``${((u\cdot \n)d\odot \n d)\cdot \n \phi}$'' which appears in \cite{linliu} has been more accurately written here as \linebreak ${u\otimes \n \phi : \n d \odot \n d}$, and that
$u\otimes \n \phi : \n d \odot \n d = [(\n d \odot \n d)\n \phi]\cdot u = [(u\cdot \n)d]\cdot[(\n \phi \cdot \n)d]$.
\subsection{Global energy regularity heuristics}
Let us first see where the  {\em  global} energy identity (\ref{globllen}) leads us to expect  weak solutions to (\ref{maineq}) to live (and hence why we assume (\ref{enspaces}) in Theorem \ref{mainthm}).
\\\\
To ease notation, in what follows let's fix $\Omega \subset \R^3$, and for $T\in (0,\infty]$ let us set $\Omega_T:=\Omega \times (0,T)$ and
$$L^r_tL^q_x(T):=L^r(0,T;L^q(\Omega)\, .$$
According to (\ref{globllen}), we expect, so long as
$$M_0:=\tfrac 12 \|u(\cdot, 0)\|_{L^2(\Omega)}^2 + \tfrac 12 \|\n d(\cdot, 0)\|_{L^2(\Omega)}^2 + \|F(d(\cdot, 0))\|_{L^1(\Omega)} <\infty\, ,$$
(which we would assume  as a requirement on the initial data), to construct solutions with $u$ in the usual Navier-Stokes spaces:
\begin{equation}\label{spacesans}
u \in L^\infty_tL^2_x(\infty) \quad \textrm{and} \quad \n u\in L^2_tL^2_x(\infty)\, .
\end{equation}
As for  $d$  we expect as well in view of (\ref{globllen}) that
\begin{equation}\label{spacesa}
\n d \in L^\infty_tL^2_x(\infty)\, ,\ \  F(d)\in L^\infty_tL^1_x(\infty) \quad \textrm{and} \quad [\D d -f(d)] \in L^2_tL^2_x(\infty)\, .
\end{equation}
The norms of all quantities in the spaces given in (\ref{spacesans}) and (\ref{spacesa}) are controlled by either $M_0$ (the $F(d)$ term) or $(M_0)^{\frac 12}$ (all other terms), by integrating  (\ref{globllen}) over $t\in (0,\infty)$.  Recalling that
\begin{equation}\label{ffdefn}
F(d) := (|d|^2-1)^2 \quad \textrm{and} \quad f(d):=4(|d|^2-1)d\, ,
\end{equation}
one sees that $|f(d)|^2 = 16F(d)|d|^2$, and
one can easily confirm the following simple estimates:

\begin{equation}\label{ffesta}
\|d\|^2_{L^\infty_tL^4_x(\infty)} \leq \|F(d)\|^{1/2}_{L^\infty_tL^1_x(\infty)} + \|1\|_{L^\infty_tL^2_x(\infty)}\, ,
\end{equation}
\begin{equation}\label{ffestb}
\|F(d)\|^{1/2}_{L^\infty_tL^{3/2}_x(\infty)} \leq \|d\|^{2}_{L^\infty_tL^6_x(\infty)} + \|1\|_{L^\infty_tL^3_x(\infty)}\, ,
\end{equation}
\begin{equation}\label{ffestc}
\|f(d)\|^2_{L^\infty_tL^2_x(\infty)} \leq 16\|F(d)\|_{L^\infty_tL^{3/2}_x(\infty)}\|d\|^2_{L^\infty_tL^6_x(\infty)}
\end{equation}
and
\begin{equation}\label{ffestd}
\|\D d\|_{L^2(\Omega_T)} \leq \|\D d - f(d)\|_{L^2(\Omega_T)} + T^{1/2} \|f(d)\|_{L^\infty_tL^{2}_x(\infty)}\, .
\end{equation}
Therefore, if we assume that
\begin{equation}\label{bdddomain}
|\Omega|<\infty\, ,
\end{equation}
and hence
$$1\in L^\infty(0,\infty;L^2(\Omega)) \cap L^\infty(0,\infty;L^3(\Omega))\, ,$$
(\ref{spacesa}) along with (\ref{ffesta}) implies that
\begin{equation}\label{spacesaa}
d \in L^\infty(0,\infty;L^4(\Omega)) \stackrel{(\ref{bdddomain})}{\subset} L^\infty(0,\infty;L^2(\Omega))\, .
\end{equation}
so that (\ref{spacesa}) and (\ref{spacesaa}) imply
\begin{equation}\label{spacesf}
d\in L^\infty(0,\infty;H^1(\Omega)) \hookrightarrow L^\infty(0,\infty;L^6(\Omega))
\end{equation}
by the Sobolev embedding, from which (\ref{ffestb}) implies that
$$F(d)\in L^\infty_tL^{3/2}_x(\infty)$$
which, along with (\ref{ffestc}) and (\ref{spacesf}),  implies that
$$f(d)\in L^\infty_tL^{2}_x(\infty)\, ,$$
from which, finally, (\ref{ffestd}) and the last inclusion in (\ref{spacesa}) implies that
\begin{equation}\label{spacesg}
\D d \in L^2(\Omega_T)   \quad \textrm{for any} \quad T<\infty\, ,
\end{equation}
with the explicit estimate (\ref{ffestd})
which can then further be controlled by $M_0$ via (\ref{spacesa}), (\ref{ffesta}), (\ref{ffestb})  and (\ref{ffestc}).
\\\\
We therefore see that it is reasonable (in view of the usual elliptic regularity theory) to expect that weak solutions to (\ref{maineq}) should have the regularities in (\ref{enspaces}) of Theorem \ref{mainthm}.
\\\\
Note further that various interpolations of Lebesgue spaces imply, for example, that for any interval $I\subset \R$ one has
\begin{equation}\label{prespacetimeinterpaaa}
L^\infty(I;L^2(\Omega))\cap L^2(I;L^6(\Omega)) \subset L^{\frac 2\alpha}(I;L^{\frac 6{3-2\alpha}}(\Omega)) \quad \textrm{for any} \ \ \alpha \in [0,1]
\end{equation}
(for example, one may take $\alpha=\frac 35$ so that $\frac 2\alpha=\frac 6{3-2\alpha} =\frac {10}3$).  Using this along with the Sobolev embedding we expect (as mentioned in Remark \ref{pspacesremark}) that
\begin{equation}\label{graddmixedleb}
(\, u \ \ \textrm{and}\, ) \ \ \n d \in L^{\frac 2\alpha}(0,T;L^{\frac 6{3-2\alpha}}(\Omega)) \quad \textrm{for any} \ \ \alpha \in [0,1]\, ,\ \  T<\infty
\end{equation}
with the explicit estimate\footnote{$A\lesssim B$ means that $A\leq CB$ for some suitably universal constant $C>0$. }
$$\|\n d\|_{L^{\frac 2\alpha}_tL^{\frac 6{3-2\alpha}}_x(T)}^{\frac 2 \alpha } \lesssim
T \|\n d\|_{L^\infty_tL^2_x(\infty)}^{\frac 2\alpha }
 + \|\n d\|_{L^\infty_tL^2_x(\infty)}^{\frac 2\alpha -2}
\|\n^2 d\|_{L^{2}(\Omega_T)}^2\, .
$$

\subsection{Local energy regularity heuristics}
Here, we will justify the well-posedness of the terms appearing in the {\em local} energy equality (\ref{locllen}), based on the expected {\em global} regularity discussed in the previous section. In fact, all but the final term in (\ref{locllen}) (where one can furthermore take the essential supremum over $t\in (0,T)$) can be seen to be well-defined by (\ref{graddmixedleb}) under the assumptions in (\ref{enspaces}) and  (\ref{pspace}).
\\\\
The  $\mathcal{R}_f(d,\phi)$ term of (\ref{locllen}) requires some further consideration: in view of (\ref{ffdefn}) we see that
\begin{equation}\label{expandgradfd}
\tfrac 14 \n^T [f(d)] =
\n^T [(|d|^2-1)d] =
2d \otimes [d\cdot (\n^T d)] +(|d|^2-1)\n^T d\, ,
\end{equation}
Recalling that
$$\mathcal{R}_f(d,\phi) :=
\phi\n^T [f(d)]:  \n^T d\, ,$$
we therefore have
\begin{equation}\label{expandrfdphi}
\tfrac 14 \mathcal{R}_f(d,\phi) =
\phi\bigg(2d \otimes [d\cdot (\n^T d)]:  \n^T d +|d|^2 |\n d|^2\bigg)-\phi |\n d|^2
\end{equation}
where we have to be careful how we handle the appearance of, essentially, $|d|^2$ in the first term (the second term is integrable in view of (\ref{spacesa})).  We have, for example, that
$$\|\phi |d|^2 |\n d|^2\|_{L^1(\Omega_T)} \leq \|\phi \|_{L^\infty(\Omega_T)}\|d \|_{L^6(\Omega_T)}^2\| \n d \|_{L^{3}(\Omega_T)}^2$$
and that
\begin{equation}\label{dlsix}
\|d \|_{L^6(\Omega_T)} <\infty \quad \textrm{for any} \ T\in (0,\infty)
\end{equation}
by (\ref{spacesf}), and either
$$\|\phi  |\n d|^2\|_{L^1(\Omega_T)} \leq \|\phi \|_{L^\infty(\Omega_T)}\|\n d \|_{L^2(\Omega_T)}^2$$
or
$$\|\phi  |\n d|^2\|_{L^1(\Omega_T)} \leq \|\phi \|_{L^3(\Omega_T)}\|\n d \|_{L^3(\Omega_T)}^2\, ,$$
(recall that $\phi$ is assumed to have compact support) and, for example, that
\begin{equation}\label{ndlsix}
\|\n d \|_{L^{10/3}(\Omega_T)} <\infty \quad \textrm{for any} \ T\in (0,\infty)
\end{equation}
by (\ref{graddmixedleb}).

\section{Proof of Theorem \ref{mainthm}}
\noindent
The first part of Theorem \ref{mainthm} is a consequence of the following ``$L^3$ $\e$-regularity" Lemma \ref{thma}, while the second part is a consequence of  the ``$\dot H^1$ $\e$-regularity" Lemma \ref{thmb} below which is itself a consequence of Lemma \ref{thma}.  In the following, for a given $z_0 = (x_0, t_0) \in \rt \times \R$ and $r>0$, as in  \cite{caf} we will adopt the following the notation for the standard parabolic cylinder $Q_r(z_0)$ as well as the following time intervals and their ``centered"
versions\footnote{These are defined in such a way that $Q^*_r(x_0,t_0)=Q_r(x_0,t_0+\tfrac{r^2}8)$, and subsequently
$$Q_{\frac r2}(x_0,t_0+\tfrac{r^2}8) = B_{\frac r2}(x_0)\times (t_0-\tfrac{r^2}8,t_0+\tfrac{r^2}8)$$
is a ``centered'' cylinder with center $(x_0,t_0)$.} (indicated with a star):
\begin{equation}\label{defnqr}
\begin{array}{c}
I_r(t_0):= (t_0 - r^2,t_0)\ , \quad I^*_r(t_0):= (t_0 - \frac 78 r^2,t_0 +  \frac 18 r^2)\ , \\\\
 Q_r(z_0):= B_r(x_0) \times I_r(t_0)\  \quad
\textrm{and} \quad Q^*_r(z_0):= B_r(x_0) \times I^*_r(t_0)\, .
\end{array}
\end{equation}

\begin{lemma}[$L^3$ $\e$-regularity, cf. Theorem 2.6 of \cite{linliu} and Proposition 1 of \cite{caf}]\label{thma}
Fix any $\bar C \in (0,\infty)$.  For each $q \in (5,6]$, there exists $\bar \e_q = \bar \e_q(\bar C) \in (0,1)$ sufficiently small\footnote{Roughly speaking, $\bar \e_q \lesssim (\bar C)^{-9}(2^{\alpha_q}-1)^9$ with $\alpha_q:=\frac{2(q-5)}{q-2}$; in particular, $\bar \e_q  \to 0$ as $q \searrow 5$.} so that for any   ${\bar z = (\bar x,\bar t) \in \rt \times \R}$ and $\bar \rho \in (0,1]$, the following holds:
\\\\
Suppose (see (\ref{defnqr}))
$u,d:Q_1(\bar z)  \to \R^3$ and $p: Q_1(\bar z) \to \R$ with
\begin{equation}\label{f}
\begin{array}{c}
u, d, \n d \in L^\infty(I_1(\bar t);L^2(B_1(\bar x)))\, , \quad
\nabla u, \n d, \n^2 d  \in L^2(Q_1(\bar z))\\\\
\textrm{and} \quad p\in L^{\frac 32}(Q_1(\bar z))
\end{array}
\end{equation}
satisfy
\begin{equation}\label{udivfree}
\nabla \cdot u =0  \quad \ \textrm{in}\ \  \mathcal{D}'(Q_1(\bar z))\, ,
\end{equation}
\begin{equation}\label{preseqJ}
 -\Delta p = \nabla \cdot (\nabla^T \cdot [u\otimes u+\n d \odot \n d ]) \quad \ \textrm{in}\ \  \mathcal{D}'(Q_1(\bar z))
\end{equation}
and the following local energy inequality holds:\footnote{See Footnote \ref{fnintdef}, and note that (\ref{locenta}) implies (\ref{locent}) with  $\bar \rho =1$ if $Q_1(\bar z)\subseteq \Omega_T$, since
$$
 \left|\left(\tfrac{|u|^2}2 +\tfrac{|\n d|^2}2\right)u\cdot \n \phi + u\otimes \n \phi : \n d \odot \n d\right|
 \leq \left(\tfrac{1}2|u|^3 +\tfrac 32 |u||\n d|^2\right)|\n \phi|
 \leq \left(|u|^3 +|\n d|^3\right)|\n \phi|\, .
$$\label{locenconstsame}}
\begin{equation}\label{locent}
 \!\boxed{\begin{array}{l}
 \int_{B_1(\bar x) \times \{t\}} \left(|u|^2 + |\n d|^2\right) \phi\, dx  + \int_{\bar t -1}^t \int_{B_1(\bar x)} \left(|\n u|^2+ |\n^2 d|^2\right) \phi \, dx\, d\tau  \\\\
\qquad \leq  \bar C\int_{\bar t -1}^t \big\{ \int_{B_1(\bar x)\times \{\tau\}} \left[ \left(|u|^2 + |\n d|^2\right)|\phi_t + \D \phi|  + (|u|^3 + |\n d|^3)|\nabla \phi| + \bar \rho |d|^2|\n d|^2 \phi \right] \, dx \\\\
\qquad \qquad \qquad \qquad \quad +\  \big|\int_{B_1(\bar x) \times \{\tau\}}  pu \cdot \nabla \phi \, dx\big|\ \big\}\, d\tau \\\\
 \textrm{for}\ \textrm{a.e.}\ t\in I_1(\bar t)\  \quad \textrm{and}\quad \forall \ \phi \in \mathcal{C}^\infty_{0}(B_1(\bar x)\times (\bar t -1,\infty)) \ \ \textrm{s.t.}\ \ \phi \geq 0\, .
\end{array}}
\end{equation}
Set\footnote{Note that $E_{3,q}<\infty$ by (\ref{f}) and standard embeddings, see Section \ref{motivation} along with (\ref{gsiginterpest}) with $\sigma =6$.}
\begin{equation}\label{etdefn}
E_{3,q}:=\intt{Q_1(\bar z)} (|u|^3+ |\n d|^3+ |p|^{\frac 3 2} +|d|^q|\n d|^{3(1-\frac q6)})\ dz\, .
\end{equation}
If $E_{3,q} \leq \bar \e_q$, then $u, \n d \in L^\infty(Q_{\frac 12}(\bar z))$ with
$$\|u\|_{L^\infty(Q_{1/2}(\bar z))}, \|\n d\|_{L^\infty(Q_{1/2}(\bar z))} \leq {\bar \e_q}^{2/9}\, .$$
\end{lemma}
\noindent
In order to prove Lemma \ref{thma}, we will require the following two technical propositions.  In order to state them, let us fix (recalling (\ref{defnqr})), for a given $z_0=(x_0,t_0)$ (to be clear by the context), the abbreviated notations

\begin{equation}\label{balls}
r_k:= 2^{-k}\ , \quad B^k:=B_{r_k}(x_0)\ , \quad I^k:=I_{r_k}(t_0)\quad \textrm{and} \quad Q^k:=B^k \times I^k
\end{equation}
(so that $Q^k= Q_{2^{-k}}(z_0)$) and, for each $k\in \N$, we define the quantities
$${L_k = L_k(z_0)} \quad \textrm{and}\quad  R_k = R_k(z_0)$$
(again, the dependence on $z_0 = (x_0,t_0)$ will be clear by context) by\footnote{We use the standard notation for averages, e.g.
$$\int_{B} \!\!\!\!\!\!\!\!\!- \ \ f(x)\ dx:=\frac 1{|B|}\int_B f(x)\, dx\, .$$}
\begin{equation}\label{Lk}
\!\!\!\!\!L_k:=  \esssup_{t\in I^k} \avint{B^k}\left(|u(t)|^2+ |\n d(t)|^2 \right)\ dx + \int_{I^k} \!  \avint{B^k} \left( |\nabla u|^2 + |\n^2 d|^2 \right)\ dx\, dt
\end{equation}
 and
\begin{equation}\label{Rk}
R_k:=\avintt{Q^k} \left( |u|^3 + |\n d|^3\right) \ dz + r_{k}^{1/3} \avintt{Q^{k}} |u||p-\bar p_{k}|\ dz
\end{equation}
$$\textrm{where} \quad \bar p_k(t):= \avint{B^k}\ p(x,t)\ dx\ .$$
$L_k$ and $R_k$ correspond roughly to the left- and right-hand sides of the local energy inequality (\ref{locent}).
We now state the technical propositions, whose proofs we will give in Section \ref{technical}:
\begin{prop}[Cf. Lemma 2.7 of \cite{linliu}]\label{la}
There exists a large universal constant $C_A >0$ such that the following holds:
\\\\
Fix any $\bar z = (\bar x, \bar t) \in \R^3 \times \R$, suppose $u$, $d$ and $p$  satisfy (\ref{f}) and (\ref{preseqJ}).
\\\\
Then  for any $z_0 \in Q_{\frac 12}(\bar z)$ we have (see (\ref{balls}), (\ref{Lk}), (\ref{Rk}))
\begin{equation}\label{c}
R_{n+1}(z_0) \leq C_A \bigg(\max_{1\leq k \leq n} L_k^{3/2}(z_0) + \underbrace{\|p\|^{3/2}_{L^{3/2}(Q_{1/2}(z_0))}}_{\leq E_{3,q} \ \forall q\geq 0,\ \mathrm{cf.}\ (\ref{etdefn})} \bigg) \qquad \forall \ \ n\geq 2\ .
\end{equation}
\end{prop}
\noindent
The proof of Proposition \ref{la} uses only the H\"older and Poincar\'e inequalities, Sobolev embedding and Calderon-Zygmund estimates along with a local decomposition  of the pressure (see (\ref{pinversion})) using the pressure equation  (\ref{preseqJ}).
\begin{prop}[Cf. Lemma 2.8 of \cite{linliu}]\label{lb} There exists a large universal constant $C_B >0$ such that the following holds:
\\\\
Fix any $\bar z = (\bar x, \bar t) \in \R^3 \times \R$, suppose $u$, $d$ and $p$  satisfy (\ref{f}), (\ref{udivfree}) and (\ref{locent}), and set $E_{3,q}$ as in (\ref{etdefn}).
\\\\
Then for any $z_0 \in Q_{\frac 12}(\bar z)$ and any $q\in (5,6]$, we have (see (\ref{balls}), (\ref{Lk}), (\ref{Rk}))
\begin{equation}\label{d}
L_n(z_0) \leq \bar C \cdot C_B \bigg(\frac 1{2^{\alpha_q}-1}\cdot\max_{k_0\leq k \leq n}R_k(z_0)
+  E_{3,q}^{2/3} + (1+k_02^{5k_0})E_{3,q} \bigg) \quad \forall \,  n \geq 2
\end{equation}
for any $k_0\in \{1, \dots, n-1\}$, where $\bar C$ is the constant from (\ref{locent}) and
$$\alpha_q:=\frac{2(q-5)}{q-2} >0\, .$$
\end{prop}
\noindent
The proof of Proposition \ref{lb} uses only the local energy inequality (\ref{locent}), the divergence-free condition (\ref{udivfree}) on $u$ and elementary estimates.  The quantities on either side of (\ref{d}) do not scale (in the sense of (\ref{recenteredscaling}))  the same way (as do those in (\ref{c})), which is why the energy inequality is necessary.
\\\\
Let us now prove Lemma \ref{thma} using Propositions \ref{la} and \ref{lb}.\\

\noindent
{\bf Proof of Lemma \ref{thma}:} \quad
Let us fix some $q\in (5,6]$ and $\bar C \in (0,\infty)$.  We first note that for any
 $\phi \geq 0$ as in (\ref{locent}) we have\footnote{The inequality in fact holds for any $q\in (2,6]$.} (recalling that $\bar \rho \leq 1$)
$$
\bar \rho \intt{Q^1}  |d|^2 |\n d|^2\phi
\leq \tfrac 2q \intt{Q^1}  |d|^q|\n d|^{3(1-\frac q6)} + (1-\tfrac 2q)\intt{Q^1}  |\n d|^3 \phi^{\frac 13(5-\alpha_q)}\, ,
$$
with $\alpha_q:=\frac{2(q-5)}{q-2} \in (0,\tfrac 12]$.
Taking $\phi$ in particular such that $\phi \equiv 1$ on $Q^1 = Q_{1/2}(z_0)$, we see easily from this  that
\begin{equation}\label{lnolessE}
\frac{L_1}{\bar C} \stackrel{(\ref{locent})\ \ }{\lesssim} E_{3,q} +  E_{3,q}^{2/3} \qquad \forall \ z_0\in Q_{1/2}(\bar z)\, .
\end{equation}
It is also easy to see that
\begin{equation}\label{lnolessln}
L_{n+1} \leq 8 L_{n} \quad  \textrm{for any}\ \  n\in \N\, .
\end{equation}
Hence  we may pick $C_0 = C_0(q,\bar C) >>1$ such that for any $z_0 \in Q_{\frac 12}(\bar z)$ (and suppressing the dependence on $z_0$ in what follows) we have
\begin{equation}\label{lonetothree}
L_1, L_2, L_3 \stackrel{(\ref{lnolessE}),(\ref{lnolessln})}{\leq} \tfrac 12(C_0)^{2/3}\left(E_{3,q} +  E_{3,q}^{2/3}\right)\, ,
\end{equation}
 $$C_A \leq \frac {C_0}2  \quad \textrm{and} \quad ((2^{\alpha_q}-1)^{-1} + 2 + 3\cdot 2^{15})\bar C \cdot C_B \leq (C_0)^{2/3}$$
for $C_A$ and $C_B$ as in Propositions \ref{la} and \ref{lb}.   Having fixed $C_0$ (uniformly over $z_0\in Q_{1/2}(\bar z)$), we then choose $\bar \e_q \in (0,1)$ so small that
$$\bar \e_q < \frac 1{(C_0)^6} \qquad \iff \qquad C_0^2 \bar \e_q < \bar \e_q^{2/3}\, .$$
Noting first that $\bar \e_q \leq (\bar \e_q)^{2/3}$, under the assumption $E_{3,q} \leq \bar \e_q$ we in particular see from (\ref{lonetothree}) that
$$L_1, L_2, L_3 \leq (C_0 \bar \e_q)^{2/3}\, .$$
Then, by Proposition \ref{la} with $n\in \{2,3\}$ we have
$$R_3,R_4 \stackrel{(\ref{c})}{\leq} \frac{C_0}{2} (\max \{L_1^{3/2},L_2^{3/2},L_3^{3/2} \} + \bar \e_q)
\leq \frac{C_0(C_0 +1)}{2}  \bar \e_q
\leq C_0^2 \bar \e_q < \bar \e_q^{2/3}$$
which implies due to Proposition \ref{lb} with $n=4$ and $k_0 =3$ that
$$L_4 \stackrel{(\ref{d})}{\leq}
C_B((2^{\alpha_q}-1)^{-1}\max \{R_3,R_4\} + E_{3,q}^{2/3} + (1+3\cdot 2^{15})E_{3,q})  \leq (C_0 \bar \e_q)^{2/3}\ .$$
Then in turn, Proposition \ref{la} with $n=4$ gives
$$L_1,L_2,L_3,L_4 \leq (C_0 \bar \e_q)^{2/3} \quad \stackrel{(\ref{c})}{\Longrightarrow} \quad R_5 < \bar \e_q^{2/3}\, ,$$
from which Proposition \ref{lb} with $n=5$  and, again, $k_0 =3$ gives
$$R_3,R_4,R_5 < \bar \e_q^{2/3} \quad \stackrel{(\ref{d})}{\Longrightarrow} \quad L_5 \leq (C_0 \bar \e_q)^{2/3}\, ,$$
and continuing we see by induction that Proposition \ref{la} and Proposition \ref{lb} (with $k_0=3$ fixed throughout) imply that
$$R_n(z_0) < \bar \e_q^{2/3} \ , \quad L_n(z_0) \leq (C_0 \bar \e_q)^{2/3}\  \qquad \forall \ n \geq 3\, .$$
This, in turn, implies (for example) that (see, e.g., \cite[Theorem 7.16]{wheeden})
$$|u(x_0,t_0)|^3 + |\n d(x_0,t_0)|^3 \leq {\bar \e_q}^{2/3}$$
$$ \textrm{for all Lebesgue points}\ z_0 \in Q_{\frac 12}(\bar z)\ \textrm{of}\ |u|^3 + |\n d|^3$$
which implies the $L^\infty$ statement, and   Lemma \ref{thma} is proved. \hfill $\Box$
\\\\
Lemma \ref{thma} will be used to prove the first assertion in Theorem \ref{mainthm} as well as the next lemma, which in turn will be used to prove the second assertion in Theorem \ref{mainthm}.
\begin{lemma}[$\dot H^1$ $\e$-regularity, cf. Theorem 3.1 of \cite{linliu} and Proposition 2 of \cite{caf}]\label{thmb}
Fix any $\bar C \in (0,\infty)$ and $\bar g \in [1,\infty)$.  For each  $\sigma \in (5,6)$,  there exists a small  constant $\e_\sigma = \e_\sigma (\bar C, \bar g)>0$ such that the following holds.  Fix  $\Omega_T:=\Omega \times (0,T)$ as in Theorem \ref{mainthm}, and suppose $u$, $d$ and $p$ satisfy assumptions (\ref{enspaces}) - (\ref{locenta}).  If (recall (\ref{defnqr}))
\begin{equation}\label{dsmallatzo}
 \limsup_{r\searrow 0} \frac 1{r^{2+ \frac \sigma 2}}  \intt{Q^*_r(z_0)}|d|^\sigma \left(|u|^3+|\n d|^3\right)^{(1-\frac \sigma 6)} \, dz \leq \bar g
\end{equation}
and
\begin{equation}\label{k}
\limsup_{r\searrow 0} \frac 1 r \intt{Q^*_r(z_0)} \left(|\nabla u|^2+|\nabla^2 d|^2\right) \, dz \leq \e_\sigma\ ,
\end{equation}
for some $z_0\in \Omega_T$, then $z_0$ is a regular point, i.e. $|u|$ and $|\n d|$ are essentially bounded in some neighborhood of $z_0$.
\end{lemma}
\noindent
For the proof of Lemma \ref{thmb}, for $z_0 =(x_0,t_0) \in  \Omega_T$ and for $r>0$ sufficiently small, we define  $A_{z_0}$, $B_{z_0}$, $C_{z_0}$, $D_{z_0}$, $E_{z_0}$, $F_{z_0}$ (cf.  \cite[(3.3)]{linliu}) and $G_{z_0}$ using the cylinders $Q^*_r(z_0)$  (whose ``centers'' $z_0$ are in the interior, see (\ref{defnqr})) by

\begin{equation}\label{athroughgdef}
\!\!\!\!\!\!\!\!\!\!\!\!\!\!\!\!\!\!\begin{array}{c}
\displaystyle{A_{z_0}(r):= \frac 1r\esssup_{t\in I_r^*(t_0)}  \int_{B_r(x_0)} \left(|u(t)|^2 + |\n d(t)|^2\right)\, dx \, , \quad }
\\\\
\displaystyle{B_{z_0}(r):= \frac 1r \intt{Q^*_r(z_0)}\left(|\nabla u|^2+|\n^2 d|^2\right)\, dz \, ,}
\\\\
\displaystyle{C_{z_0}(r):= \frac 1{r^2}\intt{Q^*_r(z_0)} \left(|u|^3 + |\n d|^3\right)\, dz \, ,\qquad  D_{z_0}(r):= \frac 1{r^2}\intt{Q_r^*(z_0)}|p|^{3/2}\, dz \, ,}
\\\\
\!\!\!\!\!\!\!\!\!\!\displaystyle{E_{z_0}(r):=\frac 1{r^2} \intt{Q^*_r(z_0)}|u| \left\{\left| |u|^2 - \overline{|u|^2}^r \right| + \left| |\n d|^2 - \overline{|\n d|^2}^r \right|\right\}\, dz}
\\\\
\left( \textrm{where} \qquad \displaystyle{\overline{g}^r(t):= \int_{B_r(x_0)}\!\!\!\!\!\!\!\!\!\!\!\!\!\!\!\!\!\!\!- \ \ \quad g(y,t) \, dy }\right)\, , \qquad\displaystyle{ F_{z_0}(r):= \frac 1{r^2} \intt{Q^*_r(z_0)} |u||p| \, dz }
\\\\
\displaystyle{  \qquad \textrm{and} \qquad  G_{q,z_0}(r):= \frac 1{r^{2+ \frac q2}}  \intt{Q^*_r(z_0)}|d|^q\left(|u|^3+|\n d|^3\right)^{(1-\frac q6)} \, dz}
\end{array}
\end{equation}
(note that $G_{0,z_0}\equiv C_{z_0}$) and, for $q\in [0,6)$, define
\begin{equation}\label{athroughgdefm}
M_{q,z_0}(r):= \tfrac 12 \left[C_{z_0}(r) + G^{\frac 6{6-q}}_{q,z_0}(r)\right] + D^2_{z_0}(r) + E_{z_0}^{\frac 3 2}(r) + F_{z_0}^{\frac 3 2}(r) \, .
\end{equation}
The statement in Lemma \ref{thmb} will follow from Lemma \ref{thma} along with the following technical ``decay estimate" which will be proved in Section \ref{technical}.
\begin{prop}[Decay estimate, cf. Lemma 3.1 of \cite{linliu} and Proposition 3 of \cite{caf}]\label{lc}
Fix any ${\bar C} \in (0,\infty)$.  There exists some constant $\bar c =\bar c({\bar C}) >0$ such that the following holds: fix any $q,\sigma \in \R$ with $2 \leq q < \sigma < 6$, and define
\begin{equation}\label{esqasqdef}
\alpha_{\sigma,q}:=\frac 6\sigma \cdot \frac {\sigma-q}{6-q} \in (0,1) \, .
\end{equation}
If $u$, $d$ and $p$ satisfy   (\ref{enspaces}) - (\ref{locenta}) for $\Omega_T$ as in Theorem \ref{mainthm}, and  $z_0\in \Omega_T$ and $\rho_0 \in (0,1]$ are such that $Q_{\rho_0}^*(z_0) \subseteq \Omega_T$ and furthermore
\begin{equation}\label{grhofinite}
\sup_{\rho \in (0,\rho_0]} B_{z_0}(\rho) \leq 1 \qquad \textrm{and} \qquad \sup_{\rho \in (0,\rho_0]} G_{\sigma, z_0}(\rho) \leq {\bar g}
\end{equation}
for some finite ${\bar g} \in [1,\infty)$, then for any $\rho \in (0,\rho_0]$ and $\gamma \in (0, \frac 1 4]$ we
have
\begin{equation}\label{l}
M_{q,z_0}(\gamma \rho) \leq \bar{c}\cdot {\bar g}^{\frac {6}{6-\sigma}} \left[\gamma^{\frac18\cdot \alpha_{\sigma,q}}
(M_{q,z_0} + M^{\alpha_{\sigma,q}}_{q,z_0}) + \gamma^{-15} B_{z_0}^{\frac 3 4 \cdot \alpha_{\sigma,q}} \sum_{k=0}^2 (M_{q,z_0}^{\frac 1 {2^k}} +
M_{q,z_0}^{\frac 1 {2^k}\cdot \alpha_{\sigma,q}})\right](\rho)\, .
\end{equation}
\end{prop}
\noindent
(In fact,  in the sum over $k$ in (\ref{l}), one can omit the term with $\alpha_{\sigma,q}$ when $k=0$.)
\\\\
The key new element in our statement and proof of Proposition \ref{lc} (and hence in achieving Lemma \ref{thmb}) is the fact that, for certain $q>0$ (so that $G_{q,z_0} \neq C_{z_0}$ and hence $M_{q,z_0}$ is notably different from the quantity found in the standard literature, namely $M_{0,z_0}$), we can still derive an estimate for $M_{q,z_0}$ of the form (\ref{l}), with a constant depending only on $\bar C$, $\sigma$ and ${\bar g}$ (and not on $q$).  This is made possible (see Claim \ref{clmalocen} and its applications in Section \ref{claimsproofs}) by the following interpolation-type estimate for the range of the quantities $G_{q,z_0}$ (including $G_{0,z_0}=C_{z_0}$), a simple consequence of H\"older's inequality:
\begin{equation}\label{gsiginterpest}
0\leq q \leq \sigma \leq 6 \quad \Longrightarrow \quad G_{q,z_0}(r) \leq G_{\sigma,z_0}^{\frac q\sigma}(r)C_{z_0}^{1-\frac q\sigma}(r) \quad \forall \ r>0\, .
\end{equation}
The estimate (\ref{gsiginterpest}) follows by writing
$$
|d|^q\left(|u|^3+|\n d|^3\right)^{(1-\frac q6)} =
\left[|d|^\sigma\left(|u|^3+|\n d|^3\right)^{(1-\frac \sigma 6)}\right]^{\frac q\sigma}\cdot \left(|u|^3+|\n d|^3\right)^{\frac {\sigma -q}\sigma}
$$
and applying H\"older's inequality with
$$1 = \frac q\sigma + \frac{\sigma -q}\sigma$$
to $G_{q,z_0}$, and noting that $r^{2+\frac q2}= [r^{2+\frac \sigma 2}]^{\frac q\sigma} \cdot [r^{2}]^{1-\frac q\sigma}$.  In particular, if $0\leq q \leq \sigma < 6$, setting
$$\alpha_{\sigma,q}:=\left(1-\frac q\sigma\right)\cdot \frac{6}{6-q}
\quad \textrm{and} \quad \beta_{\sigma,q}:=\frac q\sigma \cdot \frac {6}{6-q}$$
and noting that
$$\beta_{\sigma,q}=\frac {6}{6-\sigma}\cdot \left(1- \alpha_{\sigma,q} \right) \leq \frac {6}{6-\sigma}\, ,$$
we see that
\begin{equation}\label{gsiginterpesta}
G^{\frac{6}{6-q}}_{q,z_0}(r) \stackrel{(\ref{gsiginterpest})}{\leq}
G_{\sigma,z_0}^{\beta_{\sigma,q}}(r)C_{z_0}^{\alpha_{\sigma,q}}(r)
\stackrel{(\ref{dsmallatzo})}{\leq}
{\bar g}^{\frac{6}{6-\sigma}}\cdot \left[2 M_{q,z_0}^{\alpha_{\sigma,q}}(r)\right] \quad \forall \ r>0
\end{equation}
as long as ${\bar g}\geq 1$; this leads  to the constants appearing in (\ref{l}).
\ \\\\
Let's now use Proposition \ref{lc} and Lemma \ref{thma} to prove Lemma \ref{thmb}.
\\\\
{\bf Proof of Lemma \ref{thmb}:} \quad  Fix any $\bar C \in (0,\infty)$, $\sigma \in (5,6)$ and ${\bar g} \in [1,\infty)$, and choose\footnote{In the requirement that $q  \in (5,\min\{\sigma,\bar q\})$,  the choice of $\bar q:=\frac {11}2$ is somewhat arbitrary and taken only for concreteness; one could similarly choose any $\bar q \in (5,6)$ and adjust the subsequent constants accordingly.}    any ${q = q(\sigma) \in (5,\min\{\sigma,\tfrac {11}2\})}$ which we now also fix, noting that $\frac 6{6-q} < 12$ and $2(6-q) > 1$; for the chosen $q$, let $\bar \e_q = \bar \e_q({\bar C}) \in (0,1)$ be the corresponding small constant from Lemma \ref{thma}.
\\\\
Let us first note the following important consequence of Lemma \ref{thma}. Fix $\Omega_T$ as in Lemma \ref{thma} and $z_0:=(x_0,t_0)\in \Omega_T$, and suppose that
\begin{equation}\label{mgsmall}
M_{q,z_0}(r) \leq \frac 12\left( \frac {\bar \e_q}{3}\right)^{12}
\end{equation}
for some $r\in (0,1]$ such that $Q^*_r(z_0) \subseteq \Omega_T$.  Setting
\begin{equation}\label{recenteredscaling}
\begin{array}{c}
\!\!\!\!\!\!\!\!u_{z_0,r}(x,t):=ru(x_0+rx, t_0+ r^2t)\, , \quad p_{z_0,r}(x,t):=r^2p(x_0+rx, t_0+ r^2t)\\\\
\textrm{and}\quad d_{z_0,r}(x,t):=d(x_0+rx, t_0+ r^2t)\, ,
\end{array}\end{equation}
a change of variables from $z=(x,t)$ to
\begin{equation}\label{varchange}
(y,s):=(x_0+rx,t_0+ r^2t)
\end{equation}
 implies that

$$\int_{Q_1^*(0,0)}\left(|u_{z_0,r}|^3  +|\n d_{z_0,r}|^3+ |p_{z_0,r}|^{\frac 32} +|d_{z_0,r}|^q\left(|u_{z_0,r}|^3+|\n d_{z_0,r}|^3\right)^{(1-\frac q6)}\right)\, dz
\qquad \qquad \qquad \qquad \qquad $$
$$\qquad \qquad \qquad \qquad \qquad \qquad  =
C_{z_0}(r) + D_{z_0}(r) + G_{q,z_0}(r) \  \leq \ \left(\frac {\bar \e_q}3\right)^{12} + \left(\frac {\bar \e_q}3\right)^{6} + \left(\frac {\bar \e_q}3\right)^{2(6-q)} < \ \bar \e_q\, .
$$
Since $Q^*_1(0,0) = Q_1(0,\tfrac 18)$, it follows from assumptions (\ref{enspaces}) - (\ref{locenta}) that  $u_{z_0,r}$, $d_{z_0,r}$ and $p_{z_0,r}$ satisfy the assumptions\footnote{For example, if one fixes an arbitrary $\phi \in \mathcal{C}_0^\infty(Q_1^*(0,0))$ and sets $$\phi^{z_0,r}(x,\tau):= \phi\left(\frac{x-x_0}r, \frac{\tau-t_0}{r^2} \right)\, ,$$ then
$\phi^{z_0,r}\in \mathcal{C}_0^\infty(Q_r^*(z_0)) \subset \mathcal{C}_0^\infty(\Omega_T)$.  One can therefore use the test function $\phi^{z_0,r}$ in (\ref{locenta}), make the change of variables $(\xi,s):=\left(\frac{x-x_0}r, \frac{\tau-t_0}{r^2} \right)$ (so $(x,\tau)=(x_0+r\xi, t_0+ r^2s)$) and divide both sides of the result by $r$ to obtain the local energy inequality (\ref{locent}) for the re-scaled functions with $\bar \rho = r^2$ (as all terms scale the same way except for $|d|^2|\n d|^2\phi^{z_0,r}$) and $\bar z =(0,\tfrac 18)$. The other assumptions are straightforward.} of Lemma \ref{thma} with $\bar z = (\bar x, \bar t):=(0,\tfrac 18)$ and $\bar \rho:=r^2 \in (0,1]$, with the same constant $\bar C$ (see Footnote \ref{locenconstsame}).    Since we have just seen that $$E_{3,q}=E_{3,q}(u_{z_0,r},d_{z_0,r},p_{z_0,r},\bar z) < \bar \e_q\, ,$$ we therefore conclude by Lemma \ref{thma} that
$$|u_{z_0,r}(z)|,|\n d_{z_0,r}(z)| \leq {\bar \e_q}^{\frac 29} \qquad \textrm{for a.e.} \  z\in Q_{\frac 12}(0,\tfrac 18)
=B_{\frac 12}(0)\times(- \tfrac {1}8,\tfrac {1}8)$$
and hence
$$|u(y,s)|,|\n d(y,s)| \leq \frac {{\bar \e_q}^{\frac 29}}r \qquad \textrm{for a.e.}\   (y,s) \in
B_{\frac r2}(x_0)\times(t_0 - \tfrac {r^2}8,t_0 + \tfrac {r^2}8)\, .$$
In particular, by definition, $z_0 = (x_0,t_0)$ is a {\em regular} point, i.e. $|u|$ and $|\n d|$ are essentially bounded in a neighborhood of $z_0$, so long as (\ref{mgsmall}) holds for some sufficiently small $r>0$.
\ \\\\
In view of this fact, setting
$$\delta_\sigma:=\frac 12\left(\frac {\bar \e_{q(\sigma)}}3\right)^{12} \quad \textrm{and} \quad \bar c_\sigma := \bar{c}\cdot {\bar g}^{\frac {6}{6-\sigma}}\, ,$$
we choose $\gamma_\sigma \in (0,\tfrac 14]$ so small that furthermore
\begin{equation}\label{gsmall}
\bar{c}_\sigma  \gamma_\sigma^{\frac18\cdot \alpha_{\sigma,q}} \leq \frac 1{4} \left( \frac{\delta_\sigma^{[1-\alpha_{\sigma,q}]}}2 \right)\, ,
\end{equation}
where $\bar c = \bar c (\bar C)$ is the constant from Proposition \ref{lc} and  $\alpha_{\sigma,q}$ is defined as in (\ref{esqasqdef}); finally, we choose $\e_\sigma \in (0,1]$ so small that
\begin{equation}\label{eonesmall}
\bar c_\sigma \gamma_\sigma^{-15} \e_\sigma^{\frac 3 4 \cdot \alpha_{\sigma,q}} \leq \frac 14 \left( \frac{ \delta_\sigma^{\left[1-\frac 14 \cdot \alpha_{\sigma,q}\right]}}6\right)\, .
\end{equation}
If $z_0 \in \Omega_T$ is such that (\ref{dsmallatzo}) and (\ref{k}) hold, it implies in particular that there exists some $\rho_0\in (0,1]$ such that $Q_{\rho_0}^*(z_0) \subseteq \Omega_T$ and, furthermore,
\begin{equation}\label{rhonaughtgg}
\sup_{\rho \in (0,\rho_0]} G_{\sigma,z_0}(\rho) \leq  {\bar g}
\end{equation}
and
\begin{equation}\label{rhonaughtgb}
\sup_{\rho \in (0,\rho_0]} B_{z_0}(\rho) < \e_\sigma\, .
\end{equation}
It then follows from (\ref{gsmall}), (\ref{eonesmall}) and (\ref{rhonaughtgb})  (and the facts that $\alpha_{\sigma,q}, \delta_\sigma \leq 1$) that
$$
\bar{c}_\sigma  \gamma_\sigma^{\frac18\cdot \alpha_{\sigma,q}} \stackrel{(\ref{gsmall})}{\leq} \frac 1{4} \left( \frac{\delta_\sigma^{[1-\alpha_{\sigma,q}]}}2 \right)
= \frac 1{4} \left( \frac{\min \left\{1,\delta_\sigma^{[1-\alpha_{\sigma,q}]}\right\}}2 \right)\, ,
$$
and that

$$\!\!\!\!\!\!\!\!\!\! \bar c_\sigma\gamma_\sigma^{-15}B_{z_0}^{\frac 3 4 \cdot \alpha_{\sigma,q}}(\rho)
\stackrel{(\ref{rhonaughtgb})}{\leq}\bar c_\sigma \gamma_\sigma^{-15} \e_\sigma^{\frac 3 4 \cdot \alpha_{\sigma,q}}
\stackrel{(\ref{eonesmall})}{\leq}
\frac 14 \left( \frac{ \delta_\sigma^{\left[1-\frac 14 \cdot \alpha_{\sigma,q}\right]}}6\right)
\qquad \qquad \qquad \qquad \qquad \qquad $$$$\qquad \qquad \qquad \qquad \qquad \qquad \qquad \qquad =
\frac 14 \left( \frac{ \min_{k\in \{0,2\}} \left\{\min \left\{ \delta_\sigma^{\left[1-\frac 1{2^k} \right]},   \delta_\sigma^{\left[1-\frac 1{2^k} \cdot \alpha_{\sigma,q}\right]}  \right\}\right\}}6\right)
$$
for all $\rho \leq \rho_0$.   Suppose now that $z_0$ is not a regular point.  Then we must have
\begin{equation}\label{mbig}
\delta_\sigma < M_{q,z_0}(\rho) \qquad \textrm{for all}\ \ \rho \in (0,\rho_0]\, ,
\end{equation}
or else (\ref{mgsmall}) would hold for some $r\in (0,\rho_0]$ which would imply that $z_0$ is a regular point as we established above using Lemma \ref{thma}.
\\\\
In view of (\ref{rhonaughtgg}) and (\ref{rhonaughtgb}) (so that in particular (\ref{grhofinite}) holds, as we chose $\e_\sigma \leq 1$),  we conclude by  the estimate (\ref{l}) of  Proposition \ref{lc} (along with   (\ref{gsmall}), (\ref{eonesmall}), (\ref{rhonaughtgb}), (\ref{mbig}) and our calculations above) that
$$M_{q,z_0}(\gamma_\sigma \rho) \leq  \frac 1 2 M_{q,z_0}(\rho) \qquad \textrm{for all} \quad
\rho \in (0, \rho_0]$$
for any $z_0$ which is not a regular point.  However,  since $\g_\sigma^{k}\rho_0 \in (0,\rho_0]$ for any $k\in \N$, by iterating the estimate above we would  conclude for such $z_0$ that
$$M_{q,z_0}(\gamma_\sigma^n \rho_0) \leq   \frac 1 2 M_{q,z_0}(\g_\sigma^{n-1}\rho_0)  \leq   \frac 1 {2^2} M_{q,z_0}(\g_\sigma^{n-2}\rho_0) \leq \cdots \leq \frac 1 {2^n} M_{q,z_0}(\rho_0) < \delta_\sigma$$
for a sufficiently large $n\in \N$ which contradicts (\ref{mbig}) (with $\rho = \gamma_\sigma^n \rho_0$), and hence contradicts our assumption that $z_0$ is not a regular point.  Therefore $z_0$ must indeed be regular whenever (\ref{rhonaughtgg}) and (\ref{rhonaughtgb}) hold for our choice of $\e_\sigma$, which proves  Lemma \ref{thmb}. \hfill $\Box$
\\\\
In order to prove Theorem \ref{mainthm}, we now prove the following  general lemma, from which Lemma \ref{thma} and Lemma \ref{thmb} will have various consequences (including Theorem \ref{mainthm} as well as various other historical results, which we point out for the reader's interest).  As a motivation, note first that, for $r>0$ and $z_1:= (x_1,t_1)  \in \R^3 \times \R$, according to the notation in (\ref{recenteredscaling})
a change of variables gives
$$\int_{Q_1^*(0,0)}|u_{z_1,r}|^q + |p_{z_1,r}|^{\frac q2} = \frac 1{r^{5-q}} \int_{Q_r^*(x_1,t_1)}|u|^q + |p|^{\frac q2}\, ,
\quad
\int_{Q_1^*(0,0)}|\n u_{z_1,r}|^q = \frac 1{r^{5-2q}} \int_{Q_r^*(x_1,t_1)}|\n u|^q$$
and
\begin{equation}\label{scalinggq}
\int_{Q_1^*(0,0)}|d_{z_1,r}|^q|\n d_{z_1,r}|^{3(1-\frac q6)} = \frac 1{r^{2+\frac q2}} \int_{Q_r^*(x_1,t_1)}|d|^q|\n d|^{3(1-\frac q6)}
\end{equation}
for any $q\in [1,\infty)$.

\begin{lemma}\label{partialreggen}
Fix any open and bounded $\Omega  \subset \subset \R^3$,  $T\in (0,\infty)$, $k\geq 0$ and $C_k>0$, and suppose $\mathcal{S} \subseteq \Omega_T:=\Omega \times (0,T)$  and that $U:\Omega_T \to [0,\infty]$ is a non-negative Lebesgue-measurable function such that the following property holds in general:
\begin{equation}\label{partialreggene}
(x_0,t_0)\in \mathcal{S} \quad \Longrightarrow \quad \limsup_{r\searrow 0} \frac 1{r^{k}} \int_{Q_r^*(x_0,t_0)}\!\!\!U \, dz \ \geq\  C_k\, .
\end{equation}
If, furthermore,
\begin{equation}\label{partialreggenf}
U\in L^1(\Omega_T)\, ,
\end{equation}
then (recall Definition \ref{phausdorff}) $\mathcal{P}^{k}(\mathcal{S})<\infty$ (and hence the parabolic Hausdorff dimension of $\mathcal{S}$ is at most $k$) with the explicit estimate
\begin{equation}\label{partialreggeng}
\mathcal{P}^{k}(\mathcal{S}) \leq \frac {5^5}{C_{k}} \int_{\Omega_T}  \!U\, dz \, ;
\end{equation}
moreover, if $k=5$, then
\begin{equation}\label{partialreggengkfive}
\mu(\mathcal{S}) \leq \frac{4 \pi}3\mathcal{P}^{5}(\mathcal{S}) \leq \frac {5^5\cdot 4 \pi}{3C_5}\int_{\Omega_T}U\, dz
\end{equation}
where $\mu$ is the Lebesgue outer measure, and if $k<5$, then  in fact $\mathcal{P}^{k}(\mathcal{S})=\mu(\mathcal{S})=0$.
\end{lemma}
\ \\
Before proving Lemma \ref{partialreggen}, let's first use it along with Lemma \ref{thma} and Lemma \ref{thmb} to give the
\\
\\
{\bf Proof of Theorem \ref{mainthm}:} \quad
\\\\
First note that for any $r>0$ and $z_1:= (x_1,t_1)  \in \R^3 \times \R$ such that $Q_r(z_1) \subseteq \Omega_T$, it follows (as in the proof of Lemma \ref{thmb}) that the re-scaled triple $(u_{z_1,r},d_{z_1,r},p_{z_1,r})$ (see (\ref{recenteredscaling})) satisfies the conditions of Lemma \ref{thma} with $\bar z:=(0,0)$ and $\bar \rho :=r^2$.  Therefore if $q\in (5,6]$ and
$$\frac 1{r^{2}} \int_{Q_r(x_1,t_1)}|u|^3 +|\n d|^3 + |p|^{\frac 32}\  +\  \frac 1{r^{2+\frac q2}} \int_{Q_r(x_1,t_1)}|d|^q|\n d|^{3(1-\frac q6)}
=\qquad \qquad \qquad \qquad $$
\begin{equation}\label{rescaledsmallthreefive}
\qquad \qquad = \int_{Q_1(0,0)}|u_{z_1,r}|^3  +|\n d_{z_1,r}|^3+ |p_{z_1,r}|^{\frac 32} +|d_{z_1,r}|^q|\n d_{z_1,r}|^{3(1-\frac q6)} < \bar \e_q
\end{equation}
(with $\bar \e_q=\e_q(\bar C)$ as in Lemma \ref{thma}), it follows that $|u_{z_1,r}|,|\n d_{z_1,r}|\leq C$ on $Q_{\frac 12}(0,0)$ for some $C>0$, and hence   $|u|,|\n d|\leq \frac Cr$ on $Q_{\frac r2}(x_1,t_1)$; in particular, every interior point of $Q_{\frac r2}(x_1,t_1)$ is a regular point, assuming (\ref{rescaledsmallthreefive}) holds.  Therefore, taking $z_0:=(x_0,t_0)$ such that
$$Q_{\frac r2}(x_1,t_1) = Q^{*}_{\frac r2}(x_0,t_0)\, ,$$
(so $x_0 = x_1$ and $t_0$ is slightly lower than $t_1$ so that $(x_0,t_0)$ is in the interior of the cylinder $Q_{\frac r2}(x_1,t_1)$) and letting $\mathcal{S} \subset \Omega_T$ be the singular set of the solution $(u,d,p)$, we see (in particular) that, since $r^{2+\frac q2} < r^2$ for $r<1$,
\begin{equation}\label{upthreespacessingfive}
\left .\begin{array}{c} (x_0,t_0)\in \mathcal{S} \\ q\in (5,6] \end{array} \right\} \quad \Longrightarrow \quad \limsup_{r\searrow 0}\frac 1{r^{2+\frac q2}} \int_{Q_r^*(x_0,t_0)}|u|^3 +|\n d|^3 + |p|^{\frac 32} + |d|^q|\n d|^{3(1-\frac q6)} \geq \bar \e_q
\end{equation}
(in fact, (\ref{upthreespacessingfive}) must hold with $\liminf$ instead of $\limsup$).  Therefore, since (\ref{enspaces}) - (\ref{pspace}) imply that
\begin{equation}\label{upthreespacesfive}
|u|^3 +|\n d|^3 + |p|^{\frac 32} + |d|^q|\n d|^{3(1-\frac q6)}\in L^1(\Omega_T)
\end{equation}
(for $T<\infty$), we may apply Lemma \ref{partialreggen} (it is not hard to see, by using a suitable covering argument, that without loss of generality we can assume $\Omega$ is bounded)
with $U:= |u|^3 +|\n d|^3 + |p|^{\frac 32} + |d|^q|\n d|^{3(1-\frac q6)} $, $k=2 + \frac q2$  and $C_k:=\bar \e_q$ to see (setting $\delta := \frac{q-5}2 \in (0,\frac 12) \iff 5<q < 6$ with $2+\frac q2 = \frac 92 + \d$) that
$$
\boxed{\ \mathcal{P}^{\frac 92 + \delta}(\mathcal{S})   = 0 \quad \textrm{for any} \quad \delta \in (0,\tfrac 12)\  \, .}
$$
Before continuing with the proof of Theorem \ref{mainthm}, we describe some intermediate results (using only Lemma \ref{thma}), with historical relevance, for the interest of the reader:
\\\\
Suppose that (\ref{dmorreysmall}) holds for some $\sigma \in (5,6)$ which we now fix. We further fix any $q\in (5,\sigma)$, and choose
$\g_{\sigma,q}>0$ small enough that
$$
\g_{\sigma,q}^{1-\frac q\sigma}(\g_{\sigma,q}^{\frac q\sigma} +
({g_\sigma})^{\frac q\sigma} )< \bar \e_q\, .$$
As in the proof of (\ref{gsiginterpest}), H\"older's inequality (along with (\ref{scalinggq})) implies that
$$
\int_{Q_1(0,0)}|d_{z_1,r}|^q|\n d_{z_1,r}|^{3(1-\frac q6)} \leq ({g_\sigma})^{\frac q\sigma} \left(\int_{Q_1(0,0)}|\n d_{z_1,r}|^3\right)^{1-\frac q\sigma}\, ,
$$
so that if
\begin{equation}\label{rescaledsmallthree}
\!\!\!\!\!\frac 1{r^{2}} \int_{Q_r(x_1,t_1)}|u|^3 +|\n d|^3 + |p|^{\frac 32} = \int_{Q_1(0,0)}|u_{z_1,r}|^3  +|\n d_{z_1,r}|^3+ |p_{z_1,r}|^{\frac 32}  < \g_{\sigma,q}\, ,
\end{equation}
it follows that
$$\int_{Q_1(0,0)}|u_{z_1,r}|^3  +|\n d_{z_1,r}|^3+ |p_{z_1,r}|^{\frac 32} +|d_{z_1,r}|^q|\n d_{z_1,r}|^{3(1-\frac q6)} <
\bar \e_q$$
and hence $(x_0,t_0)\notin \mathcal{S}$ for $(x_0,t_0)$ as above.
\\\\
Therefore under the general assumption (\ref{dmorreysmall}) with $\sigma \in (5,6)$, there exists $\g_{\sigma}>0$ (e.g., $\g_\sigma:= \g_{\sigma,\frac{5+\sigma}2}$) such that
\begin{equation}\label{upthreespacessing}
(x_0,t_0)\in \mathcal{S} \quad \Longrightarrow \quad \limsup_{r\searrow 0}\frac 1{r^{2}} \int_{Q_r^*(x_0,t_0)}|u|^3 + |\n d|^3 + |p|^{\frac 32} \geq \g_\sigma\, .
\end{equation}
Therefore, as long as
\begin{equation}\label{upthreespaces}
(u,\n d,p)\in L^3(\Omega_T) \times L^3(\Omega_T) \times L^{\frac 32}(\Omega_T)  \, ,
\end{equation}
we may apply Lemma \ref{partialreggen} with $U:= |u|^3 +|\n d|^3 + |p|^{\frac 32} $, $k=2$  and $C_k:=\g_\sigma$ to see (similar to Scheffer's result in \cite{scheffer77}) that
$$\boxed{\ \mathcal{P}^2(\mathcal{S})=0\, .\ }$$
On the other hand, we know slightly more than (\ref{upthreespaces}).  The assumptions on $u$ and $d$ in (\ref{enspaces}) imply (for example, by (\ref{prespacetimeinterpaaa}) with $\alpha = \frac 35$, along with Sobolev embedding) that
$u, \n d\in L^{\frac {10}3}(\Omega_T)$.
{\em Suppose} we also knew (as in the case when $\Omega = \R^3$) that
$p\in L^{\frac 53}(\Omega_T)$
(which essentially follows from (\ref{enspaces}) and (\ref{preseq}), see \cite[Theorem 2.5]{linliu}).  Then (\ref{partialreggenf}) holds with $U:=|u|^{\frac{10}3} + |\n d|^{\frac{10}3} + |p|^{\frac 53}$, and moreover H\"older's inequality implies that
$$\left\{\frac 1{r^{2}} \int_{Q_r^*(z_0)}|u|^3 + |\n d|^{3}+ |p|^{\frac 32}\right\}^{\frac {10}9} \leq 2^{\frac{10}9}|Q_1|^{\frac 19} \left[\frac 1{r^{\frac 53}} \int_{Q_r^*(z_0)}|u|^{\frac{10}3} + |\n d|^{\frac{10}3}+ |p|^{\frac 53} \right]$$
($|Q_1|$ is the Lebesgue measure of the unit parabolic cylinder). In view of (\ref{upthreespacessing}), one could therefore apply Lemma \ref{partialreggen}  with
$$U:=|u|^{\frac{10}3} + |\n d|^{\frac{10}3}+ |p|^{\frac 53}\, , \quad k=\frac{5}3 \quad \textrm{and} \quad  C_k = \frac{ {\g_\sigma}^{\frac {10}9}}{2^{\frac{10}9}|Q_1|^{\frac 19}}\, .$$
to deduce (similar to Scheffer's result in \cite{scheffer80}) that
$$\boxed{\ \mathcal{P}^{\frac 53}(\mathcal{S})=0\, .\ }$$
All of the above follows from Lemma \ref{thma} alone.  We will now show that Lemma \ref{thmb} allows one (under assumption (\ref{dmorreysmall}) for some $\sigma\in (5,6)$, and even if $p\notin L^{\frac 53}(\Omega_T)$) to further decrease the dimension of the parabolic Hausdorff measure, with respect to which the singular set has measure zero, from $\frac 53$ to $1$.  This was essentially the most significant contribution of \cite{caf} in the Navier-Stokes setting $d\equiv 0$.
\ \\\\
Let us now proceed with the proof of the second assertion in Theorem \ref{mainthm}.
Suppose $d$ satisfies (\ref{dmorreysmall}) for some $\sigma \in (5,6)$.  Taking $\e_\sigma=\e_\sigma(\bar C, g_\sigma)>0$ as in Lemma \ref{thmb} with $\bar g:=g_\sigma$, we see from (\ref{k}) 
that
$$(x_0,t_0)\in \mathcal{S} \quad \Longrightarrow \quad \limsup_{r\searrow 0} \frac 1{r} \int_{Q_r^*(x_0,t_0)}\left(|\n u|^2+ |\n^2 d|^2\right) \geq \e_\sigma\, ,$$
so that (\ref{partialreggene}) holds with $U:=|\n u|^2+ |\n^2 d|^2$ and $k=1$.  The second assumption in (\ref{enspaces}) implies that (\ref{partialreggenf}) holds as well with $U:=|\n u|^2+ |\n^2 d|^2$.  Therefore Lemma \ref{partialreggen} with $U:=|\n u|^2+ |\n^2 d|^2$, $k=1$ and $C_k=\e_\sigma$ implies that
$$\boxed{\ \mathcal{P}^{1}(\mathcal{S})=0\, .\ }$$
This completes the proof of Theorem \ref{mainthm} (assuming Lemma \ref{partialreggen}).  \hfill $\Box$
\\\\
Let us now give the
\\\\
{\bf Proof of Lemma \ref{partialreggen}.} \quad  Fix any $\delta >0$, and any open set $V$ such that
\begin{equation}\label{sinv}
\mathcal{S}\subseteq V \subseteq \Omega \times (0,T)\ .
\end{equation}
For each $z:=(x,t)\in \mathcal{S}$, according to (\ref{partialreggene}) we can choose $r_{z} \in (0,\delta)$ sufficiently small so that $Q^*_{r_{z}}(z) \subset V$ and
\begin{equation}\label{partialreggenb}
\frac 1{r_{z}^{k}} \int_{Q_{r_{z}}^*(z)}U \geq C_k\, .
\end{equation}
By a Vitalli covering argument (see \cite[Lemma 6.1]{caf}), there exists a sequence $(z_j)_{j=1}^\infty \subseteq \mathcal{S}$ such that
\begin{equation}\label{coverofs}
\mathcal{S} \subseteq \bigcup_{j=1}^\infty Q^*_{5r_{z_j}}(z_j)
\end{equation}
and such that the set of cylinders $\{Q^*_{r_{z_j}}(z_j)\}_j$ are pair-wise disjoint. We therefore see from (\ref{partialreggenb}) that
\begin{equation}\label{hausdorffest}
\sum_{j=1}^\infty r_{z_j}^{k}\  \leq\  \frac 1{C_k}\sum_{j=1}^\infty\int_{Q_{r_{z_j}}^*(z_j)}U
\ \leq \ \frac 1{C_k}\int_{V}U\  \leq \ \frac 1{C_k}\int_{\Omega_T}U
\end{equation}
which is finite (and uniformly bounded in $\delta$) by (\ref{partialreggenf}).  Note that according to Definition \ref{phausdorff} of the parabolic Hausdorff measure $\mathcal{P}^k$, (\ref{hausdorffest}) implies
\begin{equation}\label{hausdorffestaaa}
\mathcal{P}^{k}(\mathcal{S}) \leq \frac {5^k}{C_k}\int_{V}U \leq \frac {5^k}{C_k}\int_{\Omega_T}U
\end{equation}
due to (\ref{hausdorffest}), which  establishes (\ref{partialreggeng}).
\\\\
Let us now assume that $k\leq 5$.  Letting $\mu$ be the Lebesgue (outer) measure, note that
$$\mu(Q^*_{5r_{z_j}}) \leq |B_1|(5r_{z_j})^5$$
so that
\begin{equation}\label{hausdorffesta}
\mu(\mathcal{S}) \stackrel{(\ref{coverofs})}{\leq} |B_1|\sum_{j=1}^\infty (5r_{z_j})^5 \leq 5^5|B_1|\delta^{5-k}\sum_{j=1}^\infty r_{z_j}^{k}
\stackrel{(\ref{hausdorffest})}{\leq} \delta^{5-k}\frac {5^5|B_1|}{C_k}\int_{\Omega_T}U\, ,
\end{equation}
since we have chosen $r_{z} < \delta$ for all $z\in \mathcal{S}$.  If $k=5$, (\ref{hausdorffesta}) along with Definition \ref{phausdorff} gives the explicit estimate (\ref{partialreggengkfive}) on $\mu(S)$. If $k<5$, since $\delta >0$ was arbitrary, sending $\delta \to 0$ we conclude (by (\ref{partialreggenf})) that $\mu(\mathcal{S}) =0$
and hence $\mathcal{S}$ is Lebesgue measurable with Lebesgue measure zero.  We may therefore take $V$ to be an open set such that $\mu(V)$ is arbitrarily small but so that (\ref{sinv}) still holds, and deduce that $\mathcal{P}^{k}(\mathcal{S})=0$ by (\ref{partialreggenf}) and (\ref{hausdorffestaaa}). \hfill $\Box$
\ \\

\section{Proofs of technical propositions}\label{technical}
\noindent
In order to prove Proposition \ref{la} as well as Proposition \ref{lc}, we will require certain local decompositions of the pressure (cf. \cite[(2.15)]{caf}) as follows:

\subsection{Localization of the pressure}
\noindent
\begin{clm}\label{smoothlocalpressure} Fix open sets $\Omega_1 \subset \subset \Omega_2 \subset \subset \Omega \subset \R^3$ and $\psi \in \mathcal{C}^\infty_0(\Omega_2;\R)$ with $\psi \equiv 1$ on $\Omega_1$.  Let
\begin{equation}\label{fundsoln}
G^x(y):=\frac 1{4\pi}\frac{1}{|x-y|}
\end{equation}
be the fundamental solution of $-\D$ in $\R^3$ so that, in particular,
$$\n G^x \in L^q(\Omega_2) \quad  \textrm{for any} \ \  q\in [1,\tfrac 32) $$
for any fixed $x\in \R^3$, and set
$$G^x_{\psi,1}:=-G^x \n \psi$$
$$G^x_{\psi,2}:=2  \n G^x \cdot \n \psi + G^x  \D \psi
$$
$$G^x_{\psi,3}:=
 \n G^x \otimes \n \psi+  \n \psi \otimes \n G^x  +G^x \n^2  \psi\, ,
$$
so that
$$G^x_{\psi,1}, G^x_{\psi,2}, G^x_{\psi,3} \in \mathcal{C}^\infty_0(\Omega_2) \quad  \textrm{for any fixed} \ x\in \Omega_1\, .$$
\ \\
Suppose $\Pi \in \mathcal{C}^2(\Omega;\R)$,
$v \in \mathcal{C}^1(\Omega;\R^{3})$ and
$K \in \mathcal{C}^2(\Omega;\R^{3\times 3})$.
\\\\
If
\begin{equation}\label{clasrepformav}
-\D \Pi = \n \cdot v \quad \textrm{in} \ \Omega\, ,
\end{equation}
then for any $x\in \Omega_1$,
\begin{equation}\label{clasrepformcv}
\Pi(x) = -\int \n G^x \cdot v\psi + \int G^x_{\psi,1} \cdot v + \int G^x_{\psi,2} \Pi \, .
\end{equation}
Similarly, if
\begin{equation}\label{clasrepforma}
-\D \Pi = \n \cdot (\n^T \cdot K) \quad \textrm{in} \ \Omega\, ,
\end{equation}
then for any $x\in \Omega_1$,
\begin{equation}\label{clasrepformc}
\Pi(x) = S[\psi K](x)  + \int G^x_{\psi,3} : K + \int G^x_{\psi,2} \Pi
\end{equation}
where
$$S[\widetilde K ](x):= \n_x \cdot \left(\n_x^T \cdot \int G^{x} \widetilde K\right)=  \int G^{x} \n \cdot \left(\n^T \cdot \widetilde K\right) \quad \forall \ \widetilde K \in \mathcal{C}_0^2(\Omega_2;\R^{3\times 3})\, ;$$
in particular (noting $\n^2G^x \notin L^1_{\mathrm{loc}}$), $S:\left[L^q(\Omega_2)\right]^{3\times 3} \to L^q(\Omega_2)$ for any $q\in (1,\infty)$ is a bounded, linear Calderon-Zygmund operator.
\end{clm}
\begin{remark}\label{smoothlocalpressurerk}
We note, therefore, that under the assumptions (\ref{enspaces}), (\ref{pspace}) and (\ref{preseq}), by suitable regularizations one can see that for almost every fixed   $t\in (0,T)$, (\ref{clasrepformcv}) and (\ref{clasrepformc}) hold for a.e. $x\in \Omega_1$ with $\Pi:=p(\cdot, t)$, $K:=J(\cdot, t)$ and $v:=\n^T \cdot J(\cdot, t)$ where $$J:=u\otimes u+\n d \odot \n d\, .$$
Indeed, under the assumptions (\ref{enspaces}), we have $u,\n d \in L^{\frac{10}3}(\Omega_T)$ so that (omitting the $x$-dependence)
\begin{equation}\label{inclta}
J(t)\in L^{\frac 53}(\Omega) \quad \textrm{for a.e.} \ \  t\in (0,T)\, .
\end{equation}
Moreover, since $u,\n d \in L^\infty(0,T;L^2(\Omega))\cap L^{\frac{10}3}(\Omega_T)$ and $\n u, \n^2 d\in L^2(\Omega_T)$, we have
$$\n^T \cdot J \in L^2(0,T;L^1(\Omega)) \cap L^{\frac 54}(\Omega_T)$$
so that
\begin{equation}\label{incltb}
\n^T \cdot J(t) \in L^1(\Omega)\cap L^{\frac 54}(\Omega) \quad \textrm{for a.e.} \ \  t\in (0,T)\, .
\end{equation}
Finally, (\ref{pspace}) implies that
\begin{equation}\label{incltc}
p(t)\in L^{\frac 32}(\Omega) \quad \textrm{for a.e.} \ \  t\in (0,T)\, .
\end{equation}
Fix now any $t\in (0,T)$ such that the inclusions in (\ref{inclta}), (\ref{incltb}) and (\ref{incltc}) hold.
Since $G^x_{\psi,j} \in \mathcal{C}^\infty_0$ for $x\in \Omega_1$, the terms in (\ref{clasrepformcv}) and (\ref{clasrepformc}) containing $G^x_{\psi,j}$ are all well-defined for every $x\in \Omega_1$ since $J(t), \n^T \cdot J(t), p(t) \in L^1_{loc}(\Omega)$.
The term in (\ref{clasrepformcv}) containing $\n G^x$ is in $L^r_x(\Omega_2)$ for any $r\in [1,\frac {15}7)$ by Young's convolution inequality (since $\Omega_2$ is bounded), so that term is well-defined for a.e.   $x\in \Omega_2$.  Indeed, for $R>0$ such that $\Omega_2 \subseteq B_{\frac R2}(x_0)$ for some $x_0\in \R^3$, we have $x-y \in B_R:=B_R(0)$ for all $x,y\in \Omega_2$.  Letting $G(y):=G^0(y)$ and $\chi_{B_R}$ the indicator function of $B_R$, since $\psi$ is supported in $\Omega_2$ we therefore have
$$-\int \n G^x \cdot v\psi = [([\n G] \chi_{B_R})*(v\psi)](x)$$
for all $x\in \Omega_2$.  Therefore
$$\begin{array}{rcl}
\displaystyle{\left\|\int \n G^x \cdot v\psi\right\|_{L^r_x(\Omega_2)}}
&\leq &\|([\n G] \chi_{B_R})*v\psi\|_{L^r(\R^3)}\\\\
&\leq &
\|[\n G] \chi_{B_R}\|_{L^q(\R^3)}\|v\psi\|_{L^s(\R^3)}
\\\\
&=&\|\n G \|_{L^q(B_R)}\|v\psi\|_{L^s(\Omega_2)} <\infty
\end{array}$$
by Young's inequality for any $q\in [1,\frac 32)$, $s\in [1,\frac 54)$ and $r$ such that   $1 + \frac 1r = \frac 1q + \frac 1s$ (note that $\frac 23 + \frac 45 -1 = \frac 7{15}$).
Finally, $S[\psi J(t)]\in L^{\frac 53}(\Omega_2)$ by the Calderon-Zygmund estimates (as  $1< \frac 53 < \infty$), so again that term is defined for a.e. $x\in \Omega_2$.
\\\\
Regularizing the linear equation (\ref{preseq}) using a standard spatial mollifier at any $t\in (0,T)$ where (\ref{preseq}) holds in $\mathcal{D}'(\Omega)$ and where the inclusions in (\ref{inclta}), (\ref{incltb}) and (\ref{incltc}) hold, applying Claim \ref{smoothlocalpressure} and passing to limits gives the almost-everywhere convergence (after passing to a suitable subsequence) due, in particular, to the boundedness of the linear operator $S$ on $L^{\frac 53}(\Omega_2)$.
\end{remark}
\ \\\\
\noindent
{\bf Proof of Claim \ref{smoothlocalpressure}.} \quad
Since (extending $\Pi$ by zero outside of $\Omega$) $\psi \Pi \in \mathcal{C}^\infty_0(\R^3)$,  by the classical representation formula (see, e.g.,  \cite[(2.17)]{gilbarg}), for any $x\in \R^3$ we have
\begin{equation}\label{clasrepformb}
\psi(x)\Pi(x)
= -\int G^x \D (\psi \Pi)
= -\int G^x  (\psi \D \Pi + 2\n \psi \cdot \n \Pi + \Pi \D \psi)\, .
\end{equation}
In particular, for a fixed $x\in \Omega_1$ where $\psi \equiv 1$, we have $G^x\n \psi \in \mathcal{C}_0^\infty(\R^3)$ so that
 integrating by parts in (\ref{clasrepformb}) we see that
\begin{equation}\label{clasrepform}
\Pi(x)= \int G^x \psi (- \D \Pi) +\int  G^x_{\psi,2}\Pi\, .
\end{equation}
If (\ref{clasrepformav}) holds, then by (\ref{clasrepform}) we have
\begin{equation}\label{clasrepformv}
\Pi(x)= \int G^x \psi \n \cdot v +\int  G^x_{\psi,2}\Pi
\end{equation}
for any $x\in \Omega_1$.  One can then carefully integrate by parts once in the first term of (\ref{clasrepformv}) as follows: for a small $\e>0$,
$$
\!\!\!\!\!\!\!\!\!\!\!\!\!\!\!\int_{|y-x|>\e} G^x \psi \n \cdot v\, dy =
-\int_{|y-x|>\e} [\n (G^x \psi)]\cdot v\, dy
+\frac 1{4\pi\e}\underbrace{\int_{|y-x|=\e}  \psi v\cdot \nu_y\, dS_y}_{=\mathcal{O}(\e^2)}
$$
and since the second term vanishes as $\e \to 0$ due to the fact that $|\partial B_\e(x)| \lesssim \e^2$, we conclude (since  $\n G^x \in L^1_{loc}$) that
$$\int G^x \psi \n \cdot v = -\int [\n (G^x \psi)]\cdot v= -\int \n G^x \cdot v\psi + \int G^x_{\psi,1} \cdot v$$
which, along with (\ref{clasrepformv}), implies (\ref{clasrepformcv}) for any $x\in \Omega_1$.
\\\\
On the other hand, if (\ref{clasrepforma}) holds, then by (\ref{clasrepform}) we have
\begin{equation}\label{clasrepformk}
\Pi(x)= \int G^x \psi \n \cdot (\n^T \cdot K) +\int  G^x_{\psi,2}\Pi
\end{equation}
and one can write
$$
\n \cdot (\n^T \cdot (\psi K))= [\n^2 \psi]^T:K + \n^T \psi \cdot [\n \cdot K]+ \n \psi \cdot [\n^T \cdot K]+ \psi \n \cdot (\n^T \cdot K)
$$
so that (as $\n^2 \psi = \n^T(\n \psi) = \n (\n^T \psi) = [\n^2 \psi]^T$ since $\psi \in \mathcal{C}^2$)
$$
\int G^x [\psi \n \cdot (\n^T \cdot K)] =
\int G^x [\n \cdot (\n^T \cdot (\psi K))] -
\int G^x [\n^2 \psi:K ]
\qquad \qquad
$$
$$
\qquad \qquad \qquad \qquad \qquad
-\int  \bigg([G^x\n^T \psi] \cdot [\n \cdot K] + [G^x\n \psi] \cdot [\n^T \cdot K]\bigg)\, .
$$
Since $G^x\n \psi \in \mathcal{C}^\infty_0$ for $x\in \Omega_1$, one can again integrate by parts in the final term to obtain
$$
\Pi(x)= \int G^x [\n \cdot (\n^T \cdot (\psi K))] +\int G^x_{\psi,3} : K + \int G^x_{\psi,2} \Pi
$$
for $x\in \Omega_1$ in view of (\ref{clasrepformk}).  Moreover, since $\psi K \in \mathcal{C}^2_0$ and $G^x \in L^1_{loc}$, as usual for convolutions one can change variables to obtain
$$
\int G^x \n \cdot \left(\n^T \cdot (\psi K)\right)= \left[\n_x \cdot \left(\n_x^T \cdot \int G^x \psi K\right)\right](x) =:S[\psi K](x)$$
which gives us (\ref{clasrepformc}) for any $x\in \Omega_1$, where (see, e.g., \cite[Theorem 9.9]{gilbarg}) $S$ is a singular integral operator as claimed.  (Note that $\n^2 G^x \notin L^1_{loc}$ so that one cannot simply integrate by parts twice in this term putting all derivatives on $G^x$, but $\int G^x \psi K$ is the Newtonian potential of $\psi K$ which can be twice differentiated in various senses depending on the regularity of $K$.) \hfill $\Box$
\\

\subsection{Proof of Proposition \ref{la}}
In what follows, for  $\mathcal{O}\subseteq \rt$ and $I\subseteq \R$, we will use the notation
$$\|\cdot \|_{q;\mathcal{O}}:=\|\cdot \|_{L^q(\mathcal{O})}\, , \quad \|\cdot \|_{s;I}:=\|\cdot \|_{L^s(I)}\, ,$$
$$\|\cdot \|_{q,s;\mathcal{O}\times I}:=\|\cdot \|_{L^s(I;L^q(\mathcal{O}))} =\left\|\|\cdot \|_{L^q(\mathcal{O})}\right\|_{L^s(I)}
$$
and we will abbreviate by writing
$$\|\cdot \|_{q;\mathcal{O} \times I}:=\|\cdot \|_{q,q;\mathcal{O}\times I}=\|\cdot \|_{L^q(\mathcal{O}\times I)}\, .$$
We first note some simple inequalities.  Letting $B_r \subset \rt$ be a ball of radius $r>0$, from the embedding  $W^{1,2}(B_1) \hookrightarrow L^6(B_1)$  applied to functions of the form $g_r(x) = g(rx)$ (or suitably shifted, if the ball is not centered as zero), we obtain
$$\|g_r\|_{6;B_1} \lesssim \|g_r\|_{2;B_1} + \|\n g_r\|_{2;B_1} = \|g_r\|_{2;B_1} + r\|(\n g)_r\|_{2;B_1} $$
whereupon, noting by a simple change of variables that
$$\|g_r\|_{q;B_1} = r^{-\frac 3q}\|g\|_{q;B_r}$$
for any $q\in [1,\infty)$, we obtain for any ball $B_r$ of radius $r>0$  and any $g$ that
\begin{equation}\label{scaledemb}
\|g\|_{6;B_r} \lesssim \tfrac 1r\|g\|_{2;B_r} + \|\n g\|_{2;B_r}
\end{equation}
where the constant is independent of $r$ as well as the center of $B_r$.  Next, for any $v(x,t)$, using H\"older to interpolate between $L^2$ and $L^6$ we have
\begin{equation}\label{interpsobball}
\|v(t)\|_{3;B_r} \leq \|v(t)\|_{2;B_r}^{\frac 12}\|v(t)\|_{6;B_r}^{\frac 12}
\stackrel{(\ref{scaledemb})}{\lesssim} r^{-\frac 12} \|v(t)\|_{2;B_r} + \|v(t)\|_{2;B_r}^{\frac 12}\|\n v(t)\|_{2;B_r}^{\frac 12}\, .
\end{equation}
Then for $I_r \subset \R$ with $|I_r|=r^2$ and $Q_r:=B_r \times I_r$, H\"older in the $t$ variable gives
$$\|v\|_{3;Q_r} \lesssim
r^{-\frac 12}|I_r|^{\frac 13} \|v\|_{2,\infty;Q_r} + \|v\|_{2,\infty;Q_r}^{\frac 12}
\left[|I_r|^{\frac 16}\|\n v\|_{2;Q_r}\right]^{\frac 12}
$$
so that
$$r^{-\frac 16}\|v\|_{3;Q_r} \lesssim  \|v\|_{2,\infty;Q_r} + \|v\|_{2,\infty;Q_r}^{\frac 12}
\|\n v\|_{2;Q_r}^{\frac 12} \lesssim \|v\|_{2,\infty;Q_r} +
\|\n v\|_{2;Q_r}
$$
(the first of which is sometimes called the ``multiplicative inequality")
with a constant independent of $r$.  From these, noting that $|B_r|\sim r^3$, $|Q_r|\sim r^5$,  it follows easily that, for example,
\begin{equation}\label{multineq}
\avintt{Q^n} |v|^3\ dz \lesssim  \left(\esssup_{t\in I^n} \avint{B^n}|v(t)|^2\ dx\right)^{\frac 3 2}+\left( \int_{I^k} \!  \avint{B^k} |\nabla v|^2\, dx\, dt \right)^{\frac 3 2}\, .
\end{equation}
Note also that a similar scaling argument applied to Poincar\'e's inequality gives the estimate
\begin{equation}\label{poincareball}
\|g-\overline{g_{B_r}}\|_{q;B_r} \lesssim r\|\nabla g\|_{q;B_r} \sim |B_r|^{\frac 13}\|\nabla g\|_{q;B_r}
\end{equation}
for any  $r>0$ and $q\in [1,\infty]$, where $\overline{g_{\mathcal{O}}}$ is the average of $g$ in $\mathcal{O}$ for any $\mathcal{O} \subset \R^3$ with $|\mathcal{O}|<\infty$.  Note finally that  a simple application of H\"older's inequality  gives
\begin{equation}\label{lqavg}
\|\overline{g_{\mathcal{O}}}\|_{q;\mathcal{O}} \leq \|g\|_{q;\mathcal{O}}\, .
\end{equation}
Proceeding now with the proof, fix some $\tilde \phi \in \mathcal{C}^\infty_0(\rt)$ such that
$$\tilde \phi \equiv 1 \quad  \textrm{in} \quad B_{r_2}(0)=B_{\frac 14}(0)$$ and
$$\mathrm{supp}(\tilde \phi) \subseteq B_{r_1}(0)=B_{\frac 12}(0)\, .$$
Now fix
$\bar z = (\bar x,\bar t) \in \R^3 \times \R$ and $z_0 = (x_0,t_0) \in Q_{\frac 12}(\bar z)$,
define $B^k$, $I^k$ and $Q^k$ by (\ref{balls}) for this $z_0$ and define $\phi$ by $\phi(x):=\tilde \phi (x- x_0)$.  So
$$\phi \equiv 1 \quad  \textrm{in} \quad  B^2 = B_{\frac 14}(x_0)$$
and
$$\mathrm{supp}(\phi) \subseteq B^1 =B_{\frac 12}(x_0) \subset B_{1}(\bar x) \, ,$$
since $x_0 \in B_{\frac 12}(\bar x)$.  The following estimates will clearly depend only on $\tilde \phi$,  i.e. constants will be uniform for all  $z_0 \in Q_{\frac 12}(\bar z)$).
\\\\
First, applying (\ref{multineq}) to $v\in \{u,\n d\}$ and recalling (\ref{Lk})  we see that
\begin{equation}\label{a}
\frac 1{r_n^5}\left(\|u\|_{3;Q^n}^3+\|\n d\|_{3;Q^n}^3\right) \lesssim \avintt{Q^n} (|u|^3 + |\n d|^3) \ dz  \stackrel{(\ref{multineq})}{\lesssim} L_n^{3/2}
\end{equation}
for any $n$, with a constant independent of $n$.  In particular,
\begin{equation}\label{uthreelnest}
\|u\|_{3;Q^n} +\|\n d\|_{3;Q^n} \lesssim r_n^{\frac 53} L_n^{1/2}
\end{equation}
for any $n$.
\\\\
Next, by Claim \ref{smoothlocalpressure} and Remark \ref{smoothlocalpressurerk} with $\psi:=\phi$, $\Omega_2:= B^1$ and $\Omega_1:= B^2$, at almost every \linebreak $(x,t) \in Q^2 = Q_{\frac 14}(z_0) = B_{\frac 14}(x_0)\times (t_0 - (\tfrac 14)^2, t_0)$ (where $p = \phi p$),  as in  (\ref{clasrepformc})  we have
\begin{equation}\label{pinversion}
\begin{array}{rcl}
p(x,t)
& = & \displaystyle{S[\phi J(t)](x) +\int_{B^1 \backslash B^2} (2\nabla G^x \otimes_\sigma \nabla \phi + G^x \nabla^2 \phi): J(t)\ dy}\\\\
&   & \qquad \qquad \qquad +
\displaystyle{\int_{B^1 \backslash B^2} (2\nabla G^x \cdot \nabla \phi + G^x \D \phi)p(t) \ dy}\, ,
\end{array}
\end{equation}
where
\begin{equation}\label{Juu}
J:=u\otimes u + \n d \odot \n d \, ,
\end{equation}
$2a\otimes_\sigma b:= a\otimes b + b\otimes a$ and the operator $S$  consisting of second derivatives of the Newtonian potential given by
$$S[\tilde K](x):= \n_x \cdot \left(\n_x^T \cdot \int_{B^1} G^x \tilde K\right)$$
for $\tilde K\in L^q(B^1)$ is a bounded linear Calderon-Zygmund operator on $L^q(B^1)$ for $1<q<\infty$.
Hence for any  $n\in \N$, denoting by $\chi_n$ the indicator function for the set  $B^n=B_{2^{-n}}(x_0)$ and splitting  $\phi = \chi_n \phi +(1-\chi_n)\phi$ in the first term of (\ref{pinversion}), we can write
$$p=p^{1,n} +p^{2,n} +p^{3,n}\equiv p^{1,n} +p^{2,n} +p^{3}\ ,$$
where, for almost every $(x,t) \in Q^2$,
$$p(x,t)
=\underbrace{S[\chi_{n}\phi J(t)](x)}_{=:p^{1,n}(x,t)} +
\underbrace{S[(1-\chi_{n}) \phi J(t)](x)}_{=:p^{2,n}(x,t)} +
\qquad \qquad \qquad \qquad \qquad \qquad \qquad$$
$$\qquad \qquad
+\underbrace{\int_{B^1 \backslash B^2} (2\nabla G^x \otimes_\sigma \nabla \phi + G^x \nabla^2 \phi): J(t)\ dy +
\int_{B^1 \backslash B^2} (2\nabla G^x \cdot \nabla \phi + G^x \D \phi)p(t) \ dy}_{=:p^{3,n}(x,t) \equiv p^{3}(x,t)}
$$
(where the last term is clearly independent of $n$, but we keep the notation $p^{3,n}$ for convenience).
\\\\
Note first that, by the classical Calderon-Zygmund estimates, there is a universal constant $C_{cz}>0$ such that, for all $n\in \N$, we have
\begin{equation}\label{firstpest}
\|p^{1,n}(t)\|_{\frac 32;B^{n+1}} \leq C_{cz}
\|\chi_{n}\phi J(t)\|_{\frac 32;\R^3}
\leq C_{cz}\|\tilde \phi\|_{\infty;\rt} \|J(t)\|_{\frac 32;B^n}\, .
\end{equation}
Next, since the appearance of $\n \phi$ in $p^3$ exactly cuts off a neighborhood of the singularity of $G^x$ (see (\ref{fundsoln})) uniformly for all $x\in B_{\frac 18}(x_0)$ (as we integrate over $|x_0-y| \geq \frac 14$, hence $|x-y| \geq \frac 18$), we see that $p^{3,n}(\cdot , t) \in \mathcal{C}^\infty(B_{\frac 18}(x_0))$ for $t\in I_{\frac 18}(t_0)$ with, in particular,
\begin{equation}\label{secondpest}
\|\n_x p^{3,n}(t)\|_{\infty;B^{n+1}} \stackrel{(n\geq 2)}{\leq}\|\n_x p^{3,n}(t)\|_{\infty;B_{\frac 18}(x_0)} \leq c(\tilde \phi)\left( \|J(t)\|_{1;B^1}+\|p(t)\|_{1;B^1} \right)\, .
\end{equation}
In the term $p^{2,n}$,  the singularity coming from $G^x$ is also isolated due to the appearance of $\chi_n$, but it is no longer uniform in $n$ so we must be more careful.  As we are integrating over a region which avoids a neighborhood of the singularity at $y=x$ of $G^x$, we can pass the derivatives in $S$ under the integral sign to write
$$\nabla_x p^{2,n}(x,t) = \int_{B^1 \setminus B^n} \nabla_x[(\nabla^2_x G^x)^T:\phi J(t)]\, dy = \sum_{k=1}^{n-1} \int_{B^k \setminus B^{k+1}}\nabla_x[(\nabla^2_x G^x)^T:\phi J(t)]\, dy$$
and note, in view of (\ref{fundsoln}) that
$$\left|\nabla^3_xG^x(y)\right| \lesssim \frac 1{|x-y|^4} \leq (2^{k+2})^4 \lesssim \frac{2^k}{|B^k|} \quad \forall\ x\in B^{k+2},\ y\in \left(B^{k+1}\right)^c\, .$$
Therefore, since
$$B^{n+1} = B^{(n-1)+2} \subseteq B^{k+2} \quad \textrm{for}\quad 1\leq k \leq n-1\, ,$$
we see that
\begin{equation}\label{thirdpest}
\|\n_x p^{2,n}(\cdot, t)\|_{\infty, B^{n+1}} \lesssim c(\tilde \phi)\sum_{k=1}^{n-1} 2^k\avint{B^k}|J(y, t)|\, dy
\end{equation}
for all $t\in I_{\frac 18}(t_0)$.
\\\\
Now, recalling the notation
$$\bar f_k(t):= \avint{B^k} f(x,t)\ dx$$
for a function $f(x,t)$ and $k\in \N$,
for any $t\in I^2=(t_0 - (\tfrac 14)^2,t_0)$ and $n\geq 2$, we estimate
\begin{equation}\label{calcuppbara}
\int_{B^{n+1}}|u(x,t)||p(x,t)-\bar p_{n+1}(t)|\, dx \leq \qquad \qquad \qquad \qquad \qquad \qquad \qquad \qquad \qquad \qquad \qquad \qquad
\end{equation}
$$\qquad \begin{array}{cl}
\leq &\displaystyle{ \sum_{j=1}^3 \int_{B^{n+1}}|u(x,t)||p^{j,n}(x,t)-\bar p^{j,n}_{n+1}(t)|\, dx}\\\\
\leq & \displaystyle{ \|u(\cdot,t)\|_{3;B^{n+1}}\sum_{j=1}^3\|p^{j,n}(\cdot,t)-\bar p^{j,n}_{n+1}(t)\|_{\frac 32;B^{n+1}}  } \\\\
\stackrel{
{\tiny
\begin{array}{c}(\ref{poincareball}),\\ (\ref{lqavg}),\\ \textrm{H\"older} \end{array}}}{\lesssim} & \displaystyle{
\|u(t)\|_{3;B^{n+1}}\left(\|p^{1,n}(t)\|_{\frac 32;B^{n+1}} + |B^{n+1}|\sum_{j=2}^3  \|\nabla p^{j,n}(t)\|_{\infty;B^{n+1}} \right) } \\\\
\stackrel{{\tiny
\begin{array}{c}(\ref{firstpest}),\\ (\ref{secondpest}), \\ (\ref{thirdpest}),\\ \textrm{H\"older}\end{array}}
}{\lesssim} &  \displaystyle{
\|u(t)\|_{3;B^{n+1}}\left(\|J(t)\|_{\frac 32;B^n}
 +r_{n+1}^3\left\{\left(\sum_{k=1}^{n-1} 2^k\avint{B^k}|J(t)|\, dy\right) +\|J(t)\|_{\frac 32;B^1}+\|p(t)\|_{\frac 32;B^1} \right\}\right)\, .}
\end{array}
\qquad \qquad \qquad \qquad \qquad \qquad \qquad \qquad $$
Note further that, setting
\begin{equation}\label{ljkdef}
\mathbb{L}_{J,k}:=\left\|\avint{B^k}|J(t)|\, dy\right\|_{L^\infty_t(I^{k})}\, ,
\end{equation}
we have
$$\left\|\sum_{k=1}^{n-1} 2^k\avint{B^k}|J(t)|\, dy\right\|_{L_t^{\frac 32}(I^{n+1})} \leq
|I^{n+1}|^{\frac 23}\left(\max_{1\leq k\leq n-1} \mathbb{L}_{J,k} \right) \sum_{k=1}^{n-1} 2^k
\leq
r_{n+1}^{\frac 13}\max_{1\leq k\leq n-1} \mathbb{L}_{J,k}\, ,
$$
since $|I^{n+1}| = r^2_{n+1}$ and
$$\sum_{k=1}^{n-1} 2^k =\frac{ 2^n-2}{2-1} < 2^n = r_n^{-1}\, .$$
Integrating over $t\in I^{n+1}$ in (\ref{calcuppbara}), applying H\"older in the variable $t$ and recalling by (\ref{uthreelnest}) that $\|u\|_{3;Q^{n+1}}  \lesssim r_{n+1}^{\frac 53} L_{n+1}^{1/2}$, we obtain
\begin{equation}\label{calcuppbaraJ}
\intt{Q^{n+1}}|u||p-\overline{p_{n+1}}|\, dz \lesssim \qquad  \qquad  \qquad  \qquad  \qquad  \qquad  \qquad  \qquad  \qquad  \qquad
\end{equation}
$$  \lesssim r_{n+1}^{\frac 53} L_{n+1}^{1/2}\left\{\|J\|_{\frac 32;Q^n} +
r_{n+1}^{\frac{10}3} \max_{1\leq k\leq n-1} \mathbb{L}_{J,k} +r_{n+1}^3 \left(
\|J\|_{\frac 32;Q^1}+\|p\|_{\frac 32;Q^1}\right)\right\}\, .$$
It follows now from (\ref{Juu}) that
\begin{equation}\label{Juua}
\|J\|_{\frac 32;Q^k} \leq \|u\|^2_{3;Q^k}+\|\n d\|^2_{3;Q^k} \stackrel{(\ref{uthreelnest})}{\lesssim} \left(r_{k}^{\frac 53} L_{k}^{1/2} \right)^2 =
r_k^{\frac {10}3}L_k
\end{equation}
and
\begin{equation}\label{Juuc}
\mathbb{L}_{J,k} \stackrel{(\ref{ljkdef})}{\leq} \left\|\avint{B^k}\left(|u(\cdot)|^2 + |\n d(\cdot)|^2\right)\, dy\right\|_{\infty; I^{k}} \leq L_k\, .
\end{equation}
Now from (\ref{Juu}), (\ref{calcuppbaraJ}), (\ref{Juua}), (\ref{Juuc})  and the simple fact that $\tfrac 12 r_n = r_{n+1} \leq 1$
 we obtain
$$\begin{array}{rcl}
\displaystyle{r_{n+1}^{\frac 13} \avinttt{Q^{n+1}} |u||p-\bar p_{n+1}|\, dz  }   & \lesssim & \displaystyle{ L_{n+1}^{1/2}\bigg\{r_{n+1}^{\frac 13}L_n +
 r_{n+1}^{\frac 13}\max_{1\leq k\leq n-1} L_k +
\underbrace{r_1^{\frac {10}3}}_{\leq 1}L_1 +\|p\|_{\frac 32;Q^1}\bigg\}} \\\\
&\lesssim & \displaystyle{ L_{n+1}^{1/2}\left\{
 \max_{1\leq k\leq n} L_k  +\|p\|_{\frac 32;Q^1} \right\}\, .}
\end{array}
$$
Since
$$\avinttt{Q^{n+1}} \left(|u|^3 + |\n d|^3\right) \, dz \stackrel{(\ref{a})}{\lesssim} L_{n+1}^{\frac 32}\, ,
$$
adding the previous estimates and recalling (\ref{Lk}) and (\ref{Rk}) we have
$$R_{n+1} \lesssim L_{n+1}^{\frac 32} + L_{n+1}^{1/2}\left( \max_{1\leq k\leq n} L_k  +\|p\|_{\frac 32;Q^1}\right)$$
(where the constant is universal). This along with (\ref{lnolessln})
easily implies (\ref{c}) and proves  \linebreak Proposition \ref{la}. \hfill $\Box$
\ \\\\

\subsection{\bf Proof of Proposition \ref{lb}}
\noindent
For simplicity, take $\bar z = z_0 = (0,0)$, so that (recall (\ref{balls})) $Q^k = Q^k(0,0)$, etc., as the rest can be obtained by appropriate shifts.
\\\\
We want to take the test function $\phi$ in (\ref{locent}) such that  $\phi = \phi^n:= \chi \psi^n$, where (recall that here $Q^1 = Q^1(0,0) = B_{\frac 12}(0)\times (-\tfrac 14,0)$ so $\chi$ will be zero in a neighborhood of the ``parabolic boundary'' of $Q^1$)
\begin{equation}\label{chizero}
\chi \in \mathcal{C}^\infty_{0} \left(B_{\frac 12}(0)\times \left(-\tfrac 14,\infty\right)\right)\ , \quad \chi \equiv 1\ \textrm{in}\ Q^2\ , \quad 0 \leq\chi \leq 1
\end{equation}
and
\begin{equation}\label{defnpsin}
\psi^n(x,t):= \frac 1{(r_n^2 -t)^{3/2}} e^{-\frac{|x|^2}{4(r_n^2-t)}}  \quad \textrm{for} \quad t \leq 0\, .
\end{equation}
Note that the singularity of $\psi^n$ would naturally be at $(x,t)=(0,r_n^2) \notin Q^1$, so $\psi^n\in \mathcal{C}^\infty(\overline{Q_1})$ and we may extend $\psi^n$ smoothly to $t>0$ (where it's values will actually be irrelevant) for each $n$ so that, in particular, $\phi^n \in \mathcal{C}^\infty_{0} \left(B_{1}(0)\times \left(-1,\infty\right)\right)$ as required\footnote{In (\ref{locent}) as well, the values of $\phi$ for $t>\bar t$ are actually irrelevant.} in (\ref{locent}) (with $(\bar x,\bar t)=(0,0)$). Furthermore, we have
\begin{equation}\label{bh}
\nabla \psi^n(x,t) = -\frac x{2(r_n^2-t)}\psi^n(x,t) \quad \textrm{and} \quad
\psi_t^n + \Delta \psi^n \equiv 0 \quad \textrm{in}\ \ Q^1\ .
\end{equation}
Note first that for $(x,t)\in Q^n$ ($n\geq 2$), we have
$$0\leq |x| \leq r_n \quad \textrm{and} \quad r_n^2 \leq [r_n^2 -t]\leq 2r_n^2$$
so that
$$r_n^3 = \left(r_n^2\right)^{\frac 32} e^{\frac 0{8r_n^2}}\leq (r_n^2 -t)^{3/2} e^{\frac{|x|^2}{4(r_n^2-t)}} \leq \left(2r_n^2\right)^{\frac 32}
e^{\frac{r_n^2}{4r_n^2}} = 2^{\frac 32} e^{\frac 14} r_n^3\, . $$
Hence
\begin{equation}\label{psia}
\frac{1}{2^{\frac 32} e^{\frac 14} } \cdot \frac 1{r_n^3}
\leq
\psi^n(x,t) \leq \frac 1{r_n^3} \quad \forall \ (x,t) \in Q^n
\end{equation}
and therefore (as $r_n^2 - t> 0$)
\begin{equation}\label{psic}
|\nabla_x \psi^n(x,t)| =  \frac {|x|}{2(r_n^2 - t)} |\psi^n(x,t)|
\lesssim
\frac {r_n}{r_n^2}\cdot \frac {1}{r_n^3}
= \frac 1{r_n^4} \quad \forall \ (x,t) \in Q^n\, .
\end{equation}
Next, note similarly that for $2\leq k\leq n$ and  $(x,t)\in Q^{k-1}\setminus Q^k$, we have
$$r_k \leq |x| \leq r_{k-1} = 2r_k$$
and
$$r_k^2 \leq r_n^2 +r_k^2 \leq [r_n^2 -t]\leq r_{n}^2 +r_{k-1}^2 \leq 2r_{k-1}^2 = 8r_k^2\, ,$$
so that
$$e^{\frac 1{32}}r_k^3 = \left(r_k^2\right)^{\frac 32} e^{\frac {r_k^2}{32r_k^2}}
\leq (r_n^2 -t)^{3/2} e^{\frac{|x|^2}{4(r_n^2-t)}}
\leq \left(8r_k^2\right)^{\frac 32}e^{\frac{(2r_k)^2}{4r_k^2}} = 2^{\frac 92} e r_k^3\, .$$
Therefore
\begin{equation}\label{psib}
\frac {1}{2^{\frac 92}e}\cdot \frac 1{r_k^3} \leq \psi^n(x,t) \leq \frac {1}{e^{\frac 1{32}}}\cdot \frac {1}{r_k^3} \quad \forall \ (x,t) \in Q^{k-1} \backslash Q^k \ \ (2\leq k\leq n)
\end{equation}
and hence, as in (\ref{psic}),
\begin{equation}\label{psid}
|\nabla_x \psi^n(x,t)|
\lesssim
\frac {r_k}{r_k^2}\cdot \frac {1}{r_k^3}
= \frac 1{r_k^4}
 \quad \forall \ (x,t) \in Q^{k-1}\backslash Q^k \ \ (2\leq k\leq n)\ .
\end{equation}
We can therefore estimate (for $n\geq 2$ where $\phi^n = \psi^n$ in $Q^n$):
$$\boxed{
\begin{array}{l}
\begin{array}{l}
\displaystyle{\frac{1}{2^{\frac 32} e^{\frac 14} } \cdot \frac 1{r_n^3} \left[
\esssup_{I^n} \int_{B^n}\left(|u|^2 + |\n d|^2\right) +  \intt{Q^n} \left(|\n u|^2 + |\n^2 d|^2\right)
\right]}
\end{array}\\\\
\qquad \qquad \begin{array}{l}\displaystyle{\stackrel{(\ref{psia})}{\leq}
\esssup_{I^n} \int_{B^n}\left(|u|^2 + |\n d|^2\right) \phi^n + \intt{Q^n}\left(|\n u|^2 + |\n^2 d|^2\right) \phi^n}
\\\\
\displaystyle{\stackrel{(\ref{locent})}{\leq}
\bar C\,  \bigg\{\intt{Q^1} \left[ \left(|u|^2 + |\n d|^2\right)|\phi_t^n + \D \phi^n|  + (|u|^3 + |\n d|^3)|\nabla \phi^n| + \bar \rho |d|^2|\n d|^2 \phi^n \right] }
\\\\
\qquad \qquad \qquad \qquad \qquad \qquad \displaystyle{+ \  \int_{ I^1} \bigg|\int_{B^{1} } pu\cdot \nabla \phi^n\bigg|\ \  \bigg\} \, .}
\end{array}
\end{array}
}$$

\noindent
Note that
$$\phi^n_t + \D \phi^n \stackrel{(\ref{bh})}{=} \psi^n(\chi_t + \D \chi) + 2\nabla \chi \cdot \nabla \psi^n \stackrel{(\ref{chizero})}{\equiv} 0 \ \textrm{in}\ Q^2$$
and hence, taking $k=2$ in (\ref{psib}) and (\ref{psid}), we see that
\begin{equation}\label{phiheatest}
\left|\phi^n_t + \D \phi^n\right| \lesssim \frac 1{r_2^3} + \frac 1{r_2^4} \lesssim 1 \qquad \textrm{in }\ Q^1\ ,
\end{equation}
so that
$$\boxed{\intt{Q^1} \left(|u|^2 + |\n d|^2\right)|\phi^n_t + \Delta \phi^n|  \stackrel{(\ref{phiheatest})}{\lesssim} \intt{Q^1} \left(|u|^2 + |\n d|^2\right)
\stackrel{(\ref{etdefn})}{\lesssim} E_{3,q}^{2/3}}$$
by H\"older's inequality.  Note similarly that
$$|\nabla \phi^n |= |\chi \nabla \psi^n  +  \psi^n \nabla \chi| \stackrel{(\ref{chizero})}{\lesssim} | \nabla \psi^n| + |\psi^n| \ \textrm{in}\ Q^1$$
so that (since $r_n^4 < r_n^3$) (\ref{psia}), (\ref{psic}) and (\ref{psib}), (\ref{psid}), respectively, give
\begin{equation}\label{gradphiest}
|\nabla \phi^n|  \lesssim  \frac 1{r_n^4} \quad \textrm{in} \ \ Q^n\ , \quad |\nabla \phi^n|  \lesssim \frac 1{r_k^4} \quad \textrm{in} \ \ Q^{k-1}\backslash Q^k
\end{equation}
for any $n \geq 2$ and $2\leq k\leq n$.
Therefore
$$
\sum_{k=2}^n \intt{Q^{k-1}\backslash Q^k} \left( |u|^3 + |\n d|^3\right) \underbrace{|\nabla \phi^n|}_{\lesssim r_k^{-4}}  \stackrel{(\ref{gradphiest})}{\lesssim} \left[ \max_{1\leq k \leq n-1} (r_k)^{1-\alpha}\avintt{Q^k}\left(|u|^3 + |\n d|^3\right) \right]
\sum_{k=2}^n (r_k)^{\alpha}
$$
and similarly
$$ \intt{Q^{n}} \left( |u|^3 + |\n d|^3\right) \underbrace{|\nabla \phi^n|}_{\lesssim r_n^{-4}}\stackrel{(\ref{gradphiest})}{\lesssim} \left[(r_n)^{1-\alpha} \avintt{Q^n}\left(|u|^3 + |\n d|^3\right)\right](r_n)^{\alpha}
$$
for any $\alpha \in (0,1]$, and we note that
\begin{equation}\label{rkmsum}
\sum_{k=1}^\infty (r_k)^{\alpha}
=\sum_{k=1}^\infty \left(2^{-\alpha}\right)^{k} =  \frac {1}{2^{\alpha}-1} <\infty \qquad \forall \ \ \alpha >0 \, .
\end{equation}
Hence in view of the disjoint union
\begin{equation}\label{disjun}
Q^1=\left(\bigcup_{k=2}^n Q^{k-1}\backslash Q^k\right)\cup Q^{n}
\end{equation}
we have (taking $\alpha=1$ in (\ref{rkmsum}))
$$\boxed{\intt{Q^1} \left( |u|^3 + |\n d|^3\right) |\nabla \phi^n | \lesssim  \max_{1\leq k \leq n} \avintt{Q^k} \left(|u|^3 + |\n d|^3\right)\, .}$$
Similarly, setting
$$\alpha_q:=\frac{2(q-5)}{q-2}$$
(note $\alpha_q \in (0,\tfrac 12]$ for $q\in (5,6]$), we have
$$
\bar \rho \intt{Q^1}  |d|^2 |\n d|^2\phi^n
\leq \tfrac 2q \underbrace{\intt{Q^1}  |d|^q|\n d|^{3(1-\frac q6)}}_{\leq E_{3,q}} + (1-\tfrac 2q)\intt{Q^1}  |\n d|^3 (\phi^n)^{\frac 13(5-\alpha_q)}
$$
uniformly, of course, over $\bar \rho \in (0,1]$.  Since
$$ \intt{Q^{n}} |\n d|^3 (\underbrace{ \phi^n}_{\lesssim r_n^{-3}})^{\frac 13(5-\alpha_q)} \stackrel{(\ref{psia})}{\lesssim}
(r_n)^{\alpha_q -5}\intt{Q^{n}} |\n d|^3 \lesssim
(r_n)^{\alpha_q} \avintt{Q^n}  |\n d|^3
$$
for $n\geq 2$ and similarly
$$
\intt{Q^{k}\backslash Q^{k+1}}
 |\n d|^3 (\underbrace{ \phi^n}_{\lesssim r_k^{-3}})^{\frac 13(5-\alpha_q)} \stackrel{(\ref{psib})}{\lesssim}
 ( r_k)^{\alpha_q -5}\intt{Q^{k}} |\n d|^3 \lesssim
 (r_k)^{\alpha_q }\avintt{Q^k}  |\n d|^3
$$
for $1\leq k\leq n-1$, we see that (\ref{rkmsum}) with $\alpha = \alpha_q$ and (\ref{disjun}) again give
$$\intt{Q^1}  |\n d|^3 (\phi^n)^{\frac 13(5-\alpha_q)}
\leq (2^{\alpha_q}-1)^{-1}  \max_{1\leq k \leq n} \avintt{Q^k}  |\n d|^3\, .$$
We therefore see that
$$
\boxed{ \bar \rho \intt{Q^1}  |d|^2 |\n d|^2\phi^n
\lesssim  \tfrac 25 E_{3,q}  + \tfrac 23 (2^{\alpha_q}-1)^{-1}\max_{1\leq k \leq n} \avintt{Q^k}  |\n d|^3 \quad \textrm{with} \ \alpha_q:=\frac{2(q-5)}{q-2}\, ,}$$
uniformly for any $\bar \rho \in (0,1]$ and $q\in (5,6]$.
\\\\
Putting all of the above together and recalling (\ref{Lk}), we see that for $n\geq 2$ we have
\begin{equation}\label{lnfirstestabc}
\!\!\!\!\!\!\!\!\!\!\! \frac{L_n}{\bar C}=\frac{1}{\bar C} \left[\esssup_{I^n} \avint{B^n} \left(|u|^2 + |\n d|^2\right) + \int_{I^n} \avint{B^n} \left(|\n u|^2 + |\n^2 d|^2\right)\right]  \qquad \qquad \qquad \qquad
\end{equation}
$$\qquad \qquad  \qquad \lesssim E_{3,q} + E_{3,q}^{2/3}  +
 (2^{\alpha_q}-1)^{-1}\max_{1\leq k \leq n} \avintt{Q^k} \left(|u|^3 + |\n d|^3\right)+  \int_{ I^1} \bigg|\int_{B^{1} } pu\cdot \nabla \phi^n\bigg|\, .
$$
Furthermore we claim that for $1\leq k_0 \leq n-1$ we have
\begin{equation}\label{pressure}
\!\!\!\!\!\!\!\!\!\!\!\!\!\!\!\!\!\!\!\!\boxed{\  \int_{I^1} \bigg|\int_{B^{1} } pu\cdot \nabla \phi^n\bigg|  \lesssim
 \max_{k_0\leq k \leq n}\left(r_k^{1/3}\avintt{Q^k} |p-\bar p_k ||u|\right) +  k_02^{4k_0}\intt{Q^1} |p||u|
\ .}
\end{equation}
Assuming this for the moment and continuing,  for $n\geq 2$, (\ref{lnfirstestabc}), (\ref{pressure}) and Young's convexity inequality along with the fact that, for any $k_1\geq 1$, we can estimate
$$\max_{1\leq k \leq k_1} \avintt{Q^k} \left(|u|^3 + |\n d|^3\right) \lesssim k_1 2^{5k_1}\intt{Q^1}\left(|u|^3 + |\n d|^3\right)$$
imply (recalling (\ref{Rk})) that
$$\!\!\!\!\!\!\!\!\!\!\!\!\!\!\! \frac{L_n}{\bar C}
\lesssim  E_{3,q} +  E_{3,q}^{2/3}  +
 (2^{\alpha_q}-1)^{-1}\max_{k_0\leq k \leq n} R_k + k_02^{5k_0}\underbrace{\intt{Q^1} |u|^3 + |\n d|^3+ |p|^{3/2}}_{\leq E_{3,q}}
$$
for any $k_0\in \{1, \dots, n-1\}$, and hence Proposition \ref{lb} is proved.
\\\\
To prove (\ref{pressure}), we consider additional functions $\chi_k$ (so that $\chi_k \phi^n = \chi_k \chi \psi^n$) satisfying
(recall that here $Q^k = Q^k(0,0) = B_{r_k}(0)\times (-r_k^2,0)$, so $\chi_k$ will be zero in a neighborhood of the ``parabolic boundary'' of $Q^k$)
\begin{equation}\label{chiprop}
\begin{array}{c}
\chi_k \in \mathcal{C}_{0}^\infty(\widetilde Q_{r_k}) \quad \textrm{with} \quad \widetilde Q_r:=B_{r}(0)\times (-r^2,r^2) \ \ \textrm{for} \ r>0\, ,\\\\
\chi_k \equiv 1 \quad \textrm{in} \ \ \widetilde Q_{\frac 78 r_k}\, ,\quad 0\leq \chi_k \leq 1 \quad \textrm{and} \quad |\nabla \chi_k | \lesssim \frac 1{r_k}
\end{array}\end{equation}
($\left .\chi_k\right|_{\{t>0\}}$ will again actually be irrelevant) so that in particular (as $\widetilde Q_{r_{k+2}} \subset \widetilde Q_{\frac 78 r_{k+1}}$ where \linebreak $\chi_k \equiv \chi_{k+1} \equiv 1$)
\begin{equation}\label{suppchi}
\mathrm{supp}\left(\chi_k - \chi_{k+1}\right) \subset \widetilde Q_{r_k} \backslash \widetilde Q_{r_{k+2}}\, .
\end{equation}
Then since $Q^1  = Q_{1/2}(0,0) \subset  Q_{\frac 78 }(0,0) = Q_{\frac 78 r_0}(0,0)$, we have $\chi_0 \equiv 1$ on $Q^1$ and hence for any $n\geq 2$, writing
$$\chi_0 = \chi_n + \sum_{k=0}^{n-1}(\chi_{k} - \chi_{k+1})\, ,$$
for any fixed $k_0 \in \N\cap [1,n-1]$ and at each fixed $\tau \in I^1$ we have
\begin{equation}\label{bigcalc}
\!\begin{array}{rcl}
\displaystyle{\int_{B^{1} } pu\cdot \nabla \phi^n} & \stackrel{(\ref{chiprop})}{=} &\displaystyle{\int_{B^{1} } pu\cdot \nabla [\chi_0 \phi^n] }\\\\
& = &
\displaystyle{\int_{B^{1} } pu\cdot \nabla [\chi_n \phi^n]
+\sum_{k=0}^{n-1}\int_{B^{1} } pu\cdot \nabla [(\chi_{k} - \chi_{k+1}) \phi^n]
}\\\\
& \stackrel{(\ref{chiprop}),(\ref{suppchi})}{=} &
\displaystyle{\int_{B^{n} } pu\cdot \nabla [\chi_n \phi^n]
+\sum_{k=0}^{n-1}\int_{[B^k\setminus B^{k+2}] } pu\cdot \nabla [(\chi_{k} - \chi_{k+1}) \phi^n]
}\\\\
& \stackrel{(\ref{udivfree})}{=} &
 \displaystyle{\int_{B^{n} } (p-\bar p_n)u\cdot \nabla [\chi_n \phi^n]+\sum_{k=0}^{k_0-1}\int_{[B^k\setminus B^{k+2}] } pu\cdot \nabla [(\chi_{k} - \chi_{k+1}) \phi^n]} \\\\
&&\qquad \displaystyle{+\sum_{k=k_0}^{n-1}\int_{[B^k\setminus B^{k+2}] } (p-\bar p_k)u\cdot \nabla [(\chi_{k} - \chi_{k+1}) \phi^n]\, ,}
\end{array}
\end{equation}
where
$$\bar p_k = \bar p_k (\tau)=\avint{B^k}p(x,\tau)\ dx\ .$$
Note first that (\ref{psib}), (\ref{psid}) and (\ref{suppchi}) imply (since $r_{j+1} = 2r_{j}$ for any $j$)
 that
$$|\nabla [(\chi_{k} - \chi_{k+1}) \phi^n]| \leq |\chi_{k} - \chi_{k+1}|| \nabla \phi^n| +|\phi^n||\nabla (\chi_{k} - \chi_{k+1}) | \lesssim r_k^{-4}
$$
$$
\textrm{on} \quad Q^k\setminus Q^{k+2} =
(Q^k\setminus Q^{k+1}) \cup (Q^{k+1}\setminus Q^{k+2})
$$
for any $k$, and similarly
$$|\nabla [\chi_n \phi^n]| \leq |\chi_n|| \nabla \phi^n| +|\phi^n||\nabla \chi_n | \lesssim r_n^{-4}
\quad \textrm{on} \ \ Q^n\, .$$
Therefore we can estimate (recalling again (\ref{chiprop}) and (\ref{suppchi}) when integrating $|(\ref{bigcalc})|$ over $\tau \in I^1$)
$$\int_{\tau\in I^1} \bigg|\int_{B^{1} \times \{\tau\}} pu\cdot \nabla \phi^n\bigg| \lesssim
 k_02^{4k_0}\intt{Q^1} |p||u| +\sum_{k=k_0}^{n}r_k\avintt{ Q^k}
|p-\bar p_k||u|
$$
which, along with  (\ref{rkmsum}) with $q=\frac 32$ implies (\ref{pressure}) for any $ k_0\in [1,n-1]$ as desired. \hfill $\Box$
\ \\\\

\subsection{\bf Proof of Proposition \ref{lc}}\label{claimsproofs}
\noindent
In this section we prove the technical decay estimate (Proposition \ref{lc}) used to prove Lemma \ref{thmb}.   In all of what follows, recall the definitions in (\ref{athroughgdef}) and (\ref{athroughgdefm}) of $A_{z_0}$, $B_{z_0}$, $C_{z_0}$, $D_{z_0}$, $E_{z_0}$, $F_{z_0}$, $ G_{q,z_0}$ and $M_{q,z_0}$.  We will require the following three claims which essentially appear in \cite{linliu} and which generalize certain lemmas in \cite{caf}; however we include full proofs in order to clarify certain details, and to highlight the role of $G_{q,z_0}$ (not utilized in \cite{linliu}) in Claim \ref{clmalocen} which is therefore\footnote{Note that $G_{z_0}(r) \lesssim \|d\|_\infty$ uniformly in $r$ (and $z_0$), though in our setting we may have $d\notin L^\infty$.} a  slightly  refined version of  what appears in \cite{linliu}.
\begin{clm}[General estimates (cf. Lemmas 5.1 and 5.2 in \cite{caf})]\label{clmgenfunabce}
There exist constants $c_1,c_2>0$ such that for any  $u$ and $d$ which have the regularities in (\ref{enspaces}) for  $\Omega_T:=\Omega \times (0,T)$ as in Theorem \ref{mainthm}, the estimates
\begin{equation}\label{m}
C_{z_0}(\gamma \rho)\leq  c_1 \left[\gamma^3 A_{z_0}^{\frac 3 2} + \gamma^{-3} A_{z_0}^{\frac 3 4}B_{z_0}^{\frac 3 4}\right](\rho)
\end{equation}
and
\begin{equation}\label{ag}
E_{z_0}(\gamma \rho) \leq  c_2 \left[C_{z_0}^{\frac 13 } A^{\frac 1 2}_{z_0}B^{\frac 1 2}_{z_0}\right](\gamma \rho)
\end{equation}
hold for  any  $z_0\in \R^{3+1}$ and $\rho>0$ such that $Q_\rho^*(z_0)\subseteq \Omega_T$ and any $\gamma\in (0,1]$.
\end{clm}
\begin{clm}[Estimates requiring the pressure equation (cf. Lemmas 5.3 and 5.4 in \cite{caf})]\label{clma}
There exist constants $ c_3, c_4 >0$ such that for any  $u$, $d$ and $p$ which have the regularities in (\ref{enspaces}) and (\ref{pspace}) for  $\Omega_T:=\Omega \times (0,T)$ as in Theorem \ref{mainthm} and which satisfy the pressure equation (\ref{preseq}), the estimates
\begin{equation}\label{n}
D_{z_0}(\gamma \rho)\leq c_3 \left[\g (D_{z_0}+  A_{z_0}^{\frac 3 4}B_{z_0}^{\frac 3 4} + C_{z_0}^{\frac 12}) + \g^{-5}A_{z_0}^\frac{3}{4}B_{z_0}^{\frac 32}\right](\rho)
\end{equation}
and
\begin{equation}\label{p}
F_{z_0}(\gamma \rho)\leq c_4\left[\gamma^\frac1{12}(A_{z_0} +D_{z_0}^{\frac 43}+ C_{z_0}^{\frac 23} ) + \gamma^{-10} A_{z_0}(B_{z_0}^{\frac 12}+B_{z_0}^2)\right](\rho)\, .
\end{equation}
hold for   any  $z_0\in \R^{3+1}$ and $\rho>0$ such that $Q_\rho^*(z_0)\subseteq \Omega_T$ and any $\gamma\in (0,\tfrac 12]$.
\end{clm}
\ \\
The crucial aspect of the estimates (\ref{m}), (\ref{ag}), (\ref{n})  and (\ref{p}) (which control $M_{q,z_0}(\g \rho)$)  in proving Lemma \ref{thmb} (through Proposition \ref{lc}) is that whenever a negative power of $\g$ appears, there is always a factor of $B_{z_0}$ as well, which will be small when proving Lemma \ref{thmb}.  Positive powers of $\gamma$ will similarly be small; in each term evaluated at $\rho$ (see also (\ref{o}) below), we must have either $\gamma^\alpha$ or $B_{z_0}^\alpha$ for some $\alpha >0$.
\\\\
To complete the proof of Proposition \ref{lc}, we require the following:

\begin{clm}[Estimate requiring the local energy inequality (cf. Lemma 5.5 in \cite{caf})]\label{clmalocen}
There exists a constant $ c_5 >0$ such that for any  $u$, $d$ and $p$ which have the regularities in (\ref{enspaces}) and (\ref{pspace}) for  $\Omega_T:=\Omega \times (0,T)$ as in Theorem \ref{mainthm} and such that $u$ satisfies the weak divergence-free property (\ref{divfree})  and   the local energy inequality (\ref{locenta}) holds for some constant ${\bar C}\in (0,\infty)$, the estimate
\begin{equation}\label{q}
\!\!\!\!\!\!A_{z_0}(\tfrac \rho 2)\leq c_5 \cdot {\bar C}  \left[C^\frac23 +  E+F_{z_0} + (1+[\, \cdot \, ]^2)G_q^{\frac{4}{6-q}}+  (G_q^{\frac{2}{6-q}}+C^{\frac 13}) B^{\frac 12}\right](\rho)
\end{equation}
holds for any $q\in [2,6)$ and any  $z_0\in \R^{3+1}$ and $\rho>0$ such that $Q_\rho^*(z_0)\subseteq \Omega_T$.
\end{clm}
\ \\
Postponing the proof of the claims, let us use them to prove the proposition.
\\\\
In all of what follows, we note the simple facts that, for any $\rho >0$ and $\alpha \in (0,1]$,
\begin{equation}\label{ca}
\begin{array}{c}
{\displaystyle \mathcal{K} \in \{A_{z_0},B_{z_0}\} \quad \Longrightarrow \quad \mathcal{K}(\alpha \rho) \leq \alpha^{-1}\mathcal{K}(\rho)}\ , \\\\
{\displaystyle \mathcal{K} \in \{C_{z_0},D_{z_0},E_{z_0},F_{z_0}\} \quad \Longrightarrow \quad \mathcal{K}(\alpha \rho)\leq \alpha^{-2}\mathcal{K}(\rho)}\\\\
\textrm{and} \qquad {\displaystyle G_{q,z_0}(\alpha \rho)\leq \alpha^{-2-\frac q2}G_{q,z_0}(\rho)}\ .
\end{array}
\end{equation}
\ \\
{\bf Proof of Proposition \ref{lc}.} \quad Fixing $z_0$ and $\rho_0$ as in Proposition \ref{lc}, under the assumptions in the proposition we see that estimates (\ref{m}), (\ref{ag}), (\ref{n}), (\ref{p}) and (\ref{q}) hold for all $\rho \in (0,\rho_0]$, $\gamma \in (0,\frac 12]$ and $q\in [2,6)$ by Claims \ref{clmgenfunabce},  \ref{clma} and  \ref{clmalocen}.
\\\\
Note first that (\ref{m}), (\ref{ag}) and (\ref{ca}) imply that
$$E_{z_0}(\gamma \rho) \lesssim \left[A_{z_0}B_{z_0}^{\frac 12} + \gamma^{-2}A_{z_0}^{\frac 34}B_{z_0}^{\frac 34} \right](\rho)$$
and hence, for example, there exists some $c_6>0$ such that
\begin{equation}\label{o}
E_{z_0}(\gamma \rho)\leq  c_6 \left[\gamma^2 A_{z_0} + \gamma^{-2} \left(A_{z_0}^{\frac 12}B_{z_0}^{\frac 12} + A_{z_0}B_{z_0}\right)\right](\rho)\, ,
\end{equation}
for $\rho \in (0,\rho_0]$ and $\gamma \in (0,\tfrac 12]$ (in fact, for $\gamma \in (0,1]$) and that  it follows from (\ref{q}), the assumption (\ref{grhofinite}) and the assumption that $\rho_0 \leq 1$ that there exists $c_7>0$ such that
$$
\!\!\!\!\!\! ({\bar C})^{-1}A_{z_0}(\tfrac \rho 2)\leq c_7  \left[C_{z_0}^{\frac 23}+ E_{z_0} +F_{z_0} +G_{q,z_0}^{\frac{4}{6-q}}+  (G_{q,z_0}^{\frac{2}{6-q}}+C{z_0}^{\frac 13})B_{z_0}^{\frac 12}\right](\rho)\, ,
$$
and hence, recalling (\ref{athroughgdefm}), we have that, for some $c_8>0$,
\begin{equation}\label{r}
({\bar C})^{-\frac 3 2}A^{\frac 3 2}_{z_0}\left(\rho/ 2\right)\leq c_8  \left[M_{q,z_0}(\rho) +M_{q,z_0}^{\frac 1 2}(\rho)B_{z_0}^{\frac 34}(\rho)\right]
\end{equation}
for $\rho \in (0,\rho_0]$.  We note as well that, as in (\ref{gsiginterpesta}), if $\sigma \in [q, 6)$ and if (\ref{grhofinite}) holds for some ${\bar g} \geq 1$, then
\begin{equation}\label{ggsigab}
G^{\frac{6}{6-q}}_{q,z_0}(\g \rho) \stackrel{(\ref{gsiginterpest})}{\leq}
{\bar g}^{\frac{6}{6-\sigma}}\cdot C_{z_0}^{\alpha_{\sigma,q}}(\g \rho)
\stackrel{(\ref{m})}{\leq}
{\bar g}^{\frac{6}{6-\sigma}}\cdot \left[\gamma^3 A_{z_0}^{\frac 3 2} + \gamma^{-3} A_{z_0}^{\frac 3 4}B_{z_0}^{\frac 3 4}\right]^{\alpha_{\sigma,q}}(\rho)
\end{equation}
for $\rho \in (0,\rho_0]$.
Now, writing $\gamma \rho = 2\gamma \cdot \tfrac \rho 2$ for $2\gamma \leq \frac 12$  it follows from (\ref{m}), (\ref{n}), (\ref{p}), (\ref{o}), (\ref{ggsigab}) and (\ref{athroughgdefm}) followed by an application of (\ref{ca}) (with $\alpha = \frac 12$) to all terms except for $A_{z_0}$ along with the facts that $\gamma, B_{z_0}(\rho) \leq 1$ (so that you can always estimate positive powers by $1$) as well as the fact that $\alpha_{\sigma,q}\in (0,1)$ that
$$\!\!\!\!\!\!\!\!\!\!\!\!\!\!\!\!\!\!\!\!\!\!\!\!\!\!\!\!\!\!\!\!\!M_{q,z_0}(\g \rho) \leq
[C_{z_0} + G^{\frac 6{6-q}}_{q,z_0} + D^2_{z_0} + E_{z_0}^{\frac 3 2} + F_{z_0}^{\frac 3 2}](\g \rho) \qquad \qquad \qquad \qquad \qquad \qquad \qquad \qquad \qquad \qquad \qquad
$$
$$\lesssim
\left[\gamma^3 A_{z_0}^{\frac 3 2}(\tfrac \rho 2) + \gamma^{-3} A_{z_0}^{\frac 3 4}(\tfrac \rho 2)B_{z_0}^{\frac 3 4}(\rho)\right]
+
{\bar g}^{\frac{6}{6-\sigma}}\cdot\left[\gamma^3 A_{z_0}^{\frac 3 2}(\tfrac \rho 2) + \gamma^{-3} A_{z_0}^{\frac 3 4}(\tfrac \rho 2)B_{z_0}^{\frac 3 4}(\rho)\right]^{\alpha_{\sigma,q}}
\qquad \qquad \qquad $$
$$
+ \left[\g M_{q,z_0}^{\frac 12}(\rho)  + \g^{-5}A_{z_0}^\frac{3}{4}(\tfrac \rho 2)\left(B_{z_0}^{\frac 3 4}(\rho)+B_{z_0}^{\frac 32}(\rho)\right)\right]^2
$$
$$
\qquad  + \left[\gamma^2 A_{z_0}(\tfrac \rho 2) + \gamma^{-2} \left(A_{z_0}^{\frac 12}(\tfrac \rho 2)B_{z_0}^{\frac 12}(\rho) + A_{z_0}(\tfrac \rho 2)B_{z_0}(\rho)\right)\right]^{\frac 3 2}
$$
$$
\qquad \qquad \qquad + \left[\gamma^\frac1{12}\left(A_{z_0}(\tfrac \rho 2) +M_{q,z_0}^{\frac 23}(\rho) \right) + \gamma^{-10} A_{z_0}(\tfrac \rho 2)\left(B_{z_0}^{\frac 12}(\rho)+B_{z_0}^2(\rho)\right)\right]^{\frac 3 2}
$$
$$\!\!\!\!\!\!\!\!\!\!\!\!\!\!\!\lesssim
(1+{\bar g}^{\frac{6}{6-\sigma}})\left[\gamma^{\frac {\alpha_{\sigma,q}}8} \left(M_{q,z_0}(\rho)+[A_{z_0}^{\frac {3} 2}(\tfrac \rho 2)]^{\alpha_{\sigma,q}}+[A_{z_0}^{\frac 3 2}(\tfrac \rho 2)]\right) + \gamma^{-15} \left([A_{z_0}^{\frac 3 2}(\tfrac \rho 2)]^{\frac {\alpha_{\sigma,q}}2}+[A_{z_0}^{\frac 3 2}(\tfrac \rho 2)]^{\frac 12}+[A_{z_0}^\frac{3}{2}(\tfrac \rho 2)]\right)B_{z_0}^{\frac {3\alpha_{\sigma,q}} 4}(\rho)\right]
$$
so long as $\gamma \in (0,\tfrac 14]$.  Noting that $1\leq {\bar g}^{\frac{6}{6-\sigma}}$, the estimate (\ref{l}) for such $\gamma$ and for $\rho \in (0,\rho_0]$ now follows from the estimate above along with (\ref{r})  as, in particular, (\ref{r}) implies (as $\gamma, B_{z_0}(\rho) \leq 1$ and $\alpha_{\sigma,q}\in (0,1)$) that
$$
({\bar C})^{-\frac 3 2}A^{\frac 3 2}_{z_0}\left(\tfrac \rho 2\right)\lesssim M_{q,z_0}(\rho) +\g^{-15-\frac{\alpha_{\sigma,q}}8}M_{q,z_0}^{\frac 1 2}(\rho)B_{z_0}^{\frac {3\alpha_{\sigma,q}}4}(\rho)
$$
which we apply to the terms above with the positive power of $\gamma$, and that
$$
({\bar C})^{-\frac 3 2} A^{\frac 3 2}_{z_0}\left(\tfrac \rho 2\right)\lesssim M_{q,z_0}(\rho) +M_{q,z_0}^{\frac 1 2}(\rho)\, ,
$$
which we apply to the terms above with the negative power of $\gamma$.  This completes the proof of \linebreak Proposition \ref{lc}. \hfill $\Box$
\\\\\\
Let us now prove the claims:
\\\\
{\bf Proof of Claim \ref{clmgenfunabce}:}  \quad For simplicity, we will suppress the dependence on $z_0=(x_0,t_0)$ in what follows.
\\\\
Let us first prove (\ref{m}).  Note that for any $r\leq \rho$, at any fixed $t\in I_r^*$, taking $v\in \{u,\n d\}$ we have
$$
\int_{B_r}|v|^2\, dx \leq
\int_{B_\rho}\left||v|^2 - \overline{|v|^2}^\rho \right|\, dx \ + \ |B_r|\, \overline{|v|^2}^\rho
\lesssim \rho \int_{B_\rho}\left|\n |v|^2\right| \, dx \ + \ \left(\frac r\rho \right)^3 \int_{B_\rho}|v|^2\, dx
$$
due to Poincar\'e's inequality (\ref{poincareball}).  Since $\left|\n |v|^2\right|\leq |v||\n v|$ almost everywhere,  H\"older's inequality then implies that
\begin{equation}\label{ubrtworrho}
\|v\|_{2;B_r}^2 \lesssim \rho \|v\|_{2;B_\rho}\|\n v\|_{2;B_\rho} + \left(\frac r\rho \right)^3\|v\|_{2;B_\rho}^2\, .
\end{equation}
Therefore
$$\|v\|_{3;B_r}^3 \stackrel{(\ref{interpsobball})}{\lesssim}
\frac 1{r^{\frac 32}}\left(\|v\|_{2;B_r}^2 \right)^{\frac 32} + \|v\|_{2;B_r}^{\frac 32}\|\n v\|_{2;B_r}^{\frac 32}\qquad \qquad \qquad
$$
$$
\qquad \qquad \quad \ \stackrel{(\ref{ubrtworrho})}{\lesssim}
\left(1+\left(\frac \rho r\right)^{\frac 32} \right)\|v\|_{2;B_\rho}^{\frac 32}\|\n v\|_{2;B_\rho}^{\frac 32}
+\frac 1{r^{\frac 32}}\left(\frac r \rho \right)^{\frac 92}\|v\|_{2;B_\rho}^{3}\, .
$$

\noindent
Summing over $v\in \{u,\n d\}$, we see that
$$\|u\|_{3;B_r}^3 + \|\n d\|_{3;B_r}^3 {\lesssim}
\left(1+\left(\frac \rho r\right)^{\frac 32} \right)\left(\|u\|_{2;B_\rho}^2 + \|\n d\|_{2;B_\rho}^2\right)^{\frac 34}\left(\|\n u\|_{2;B_\rho}^2 + \|\n^2 d\|_{2;B_\rho}^2\right)^{\frac 34}
$$
$$\qquad \qquad \qquad \qquad +\frac {r^3} {\rho^{\frac 92}}\left(\|u\|_{2;B_\rho}^2+\|\n d\|_{2;B_\rho}^2\right)^{\frac 32}\, .
$$
Now integrating over $t\in I_r^*$ (where $|I_r^*|=r^2$), H\"older's inequality implies that
$$r^2C(r)
\lesssim
|I^*_r|^{\frac 14}\left(1+\left(\frac \rho r\right)^{\frac 32} \right)\left\|\|u\|_{2;B_\rho}^2 + \|\n d\|_{2;B_\rho}^2\right\|_{\infty;I_r^*}^{\frac 34}\left(\|\n u\|_{2;Q_\rho^*}^2 + \|\n^2 d\|_{2;Q_\rho^*}^2\right)^{\frac 34}
$$
$$\qquad \qquad \qquad \qquad +|I^*_r| \frac {r^3} {\rho^{\frac 92}}\left\|\|u\|_{2;B_\rho}^2+\|\n d\|_{2;B_\rho}^2\right\|_{\infty;I_r^*}^{\frac 32}
$$
$$\lesssim
r^{\frac 12}\left(1+\left(\frac \rho r\right)^{\frac 32} \right) (\rho A(\rho))^{\frac 34}(\rho B(\rho))^{\frac 34}
+\frac {r^5} {\rho^{\frac 92}}(\rho A(\rho))^{\frac 32}\, ,\qquad \qquad \qquad
$$
which, upon  dividing both sides by $r^2$, setting $\gamma:=\frac r\rho$ and noting that $1\leq \gamma^{-\frac 32}$, precisely gives (\ref{m}).
\\\\
Next, to prove (\ref{ag}), we use the Poincar\'e-Sobolev inequality
$$\|g-\overline{g}^r\|_{q^*;B_r} \leq c_q \|\n g\|_{q;B_r}$$
(the constant is independent of $r$ due to the relationship between $q$ and $q^*$) corresponding to the embedding ${W^{1,q}\hookrightarrow L^{q^*}}$  for $q<3$ (in $\R^3$) and $q^*=\frac{3q}{3-q}$.  Taking $q=1$, at any $t\in I_r^*$ and for $v\in \{u,\n d\}$  the H\"older and Poincar\'e-Sobolev inequalities give us
$$\int_{B_r}|u|\left||v|^2 - \overline{|v|^2}^r \right|\, dx \leq  \|u\|_{3;B_r} \|\ |v|^2 - \overline{|v|^2}^r\ \|_{\frac 3 2;B_r}
\qquad \qquad \qquad \qquad $$
$$
\qquad \qquad \qquad \qquad \qquad \qquad \lesssim \|u\|_{3;B_r} \|\n (|v|^2)\|_{1;B_r}
\lesssim \|u\|_{3;B_r} \|v\|_{2;B_r} \|\n v\|_{2;B_r}\, . $$
Summing this first over $v\in \{u,\n d\}$ at a fixed $t$ and then integrating over $t\in I_r^*$, we see that
$$
\begin{array}{rcl}
{\displaystyle r^2 E(r)}
& \lesssim & {\displaystyle \int_{I_r^*} \|u\|_{3;B_r}\left( \|u\|_{2;B_r}^2+\|\n d\|_{2;B_r}^2\right)^{\frac 12} \left(\|\n d\|_{2;B_r}^2+\|\n^2 d\|_{2;B_r}^2\right)^{\frac 12} \, dt}\\\\
&\lesssim & \|u\|_{3;Q_r^*} \left\|\left( \|u\|_{2;B_r}^2+\|\n d\|_{2;B_r}^2\right)^{\frac 12}\right\|_{6;I_r^*}
\left(\|\n u\|_{2;Q_r^*}^2 + \|\n^2 d\|_{2;Q_r^*}^2\right)^{\frac 12} \\\\
&\lesssim & |I^*_r|^{\frac 16}\left(\|u\|_{3;Q_r^*}^3\right)^{\frac 13} \left\|\|u\|_{2;B_r}^2+\|\n d\|_{2;B_r}^2\right\|_{\infty;I_r^*}^{\frac 12}
\left(\|\n u\|_{2;Q_r^*}^2 + \|\n^2 d\|_{2;Q_r^*}^2\right)^{\frac 12} \\\\
&\lesssim & r^{\frac 13} (r^2 C(r))^{\frac 13} (rA(r))^{\frac 12}
\left(rB(r)\right)^{\frac 12}\ \  =\ \  r^2 [C^{\frac 13} A^{\frac 12}
B^{\frac 12}](r)
\end{array}
$$
which proves (\ref{ag}) and completes the proof of Claim \ref{clmgenfunabce}. \hfill $\Box$
\\\\\\
{\bf Proof of Claim \ref{clma}:}
\\\\
As in (\ref{clasrepformcv}) of Claim \ref{smoothlocalpressure}, for any  $t\in I_r^*(z_0)$ ($r\leq \rho$) we use Remark \ref{smoothlocalpressurerk} to decompose $\Pi:=p(\cdot, t)$ for almost every $x\in B_{\frac{3\rho}4}(x_0)$ using a smooth cut-off function $\psi$ equal to one in $\Omega_1:=B_{\frac{3\rho}4}(x_0)$ and supported in $\Omega_2:=B_\rho(x_0)$, so that
\begin{equation}\label{ac}
 |\n \psi| \lesssim \rho^{-1} \quad  \textrm{and} \quad  |\D \psi| \lesssim \rho^{-2}\, ,
\end{equation}
as
$$
p(x,t) = \underbrace{-\int \n G^x \cdot v(t)\psi\, dy}_{=:p_1(x,t)} + \underbrace{\int G^x_{\psi,1} \cdot v(t)\, dy}_{=:p_2(x,t)} + \underbrace{\int G^x_{\psi,2} p(\cdot, t)\, dy}_{=:p_3(x,t)}
$$
with
$$G^x_{\psi,1}:=-G^x \n \psi\, , \qquad
G^x_{\psi,2}:=2  \n G^x \cdot \n \psi + G^x  \D \psi
$$
and
$$v(t):=[\n^T \cdot (u \otimes u + \n d \odot \n d)](\cdot,t)\, .$$
Our goal is  to estimate $p(x,t)$ for $x\in B_{\frac \rho 2}(x_0)$.
\\\\
Both $p_2$ and $p_3$ contain derivatives of $\psi$ in each term so that the integrand can only be non-zero when $|y-x_0|>\frac {3\rho}4$, and hence for $x\in B_{\frac \rho 2}(x_0)$ one has
\begin{equation}\label{acgx}
|x-y| \geq \frac \rho 4 \quad \Longrightarrow \quad |G^x(y)| \lesssim \rho^{-1} \quad \textrm{and} \quad |\n G^x(y)| \lesssim \rho^{-2}\, .
\end{equation}
In view of (\ref{ac}) and (\ref{acgx})  and the fact that $\psi$ is supported in $B_\rho(x_0)$, we have (omitting the dependence on $t$, and noting that the constants in the inequalities are independent of $t$ as they come only from $G^x$ and $\psi$)

\pagebreak
\begin{equation}\label{ba}
\!\!\!\!\!\sup_{x\in B_{\frac \rho 2}(x_0)}|p_2(x)| \lesssim \rho^{-2} \int_{B_\rho(x_0)} (|u||\n u| + |\n d||\n^2 d|)\, dy \qquad \qquad \qquad \qquad \qquad
\end{equation}
$$
\qquad \qquad \ \  \qquad \qquad \lesssim \rho^{-2} \left(\int_{B_\rho(x_0)} (|u|^2 + |\n d|^2)\, dy \right)^{\frac12} \left(\int_{B_\rho(x_0)} (|\n u|^2 + |\n^2 d|^2)\, dy \right)^{\frac12}$$
and similarly
\begin{equation}\label{ah}
\sup_{x\in B_{\frac \rho 2}(x_0)}|p_3(x)| \lesssim \rho^{-3}\int_{B_\rho(x_0)} |p|\, dy\, .
\end{equation}
For $p_1$, Young's inequality for convolutions (where we set $R:=2\rho$ as in Remark \ref{smoothlocalpressurerk}) with \linebreak $2/3 +1=3/4 + 11/12$ gives
$$
\|p_1\|_{\frac32;B_\rho(x_0)} \lesssim \left\|\frac{1}{|\cdot|^2} \right\|_{\frac{4}{3};B_{2\rho}(0)} \left\|(|u| + |\n d|)( |\n u|+|\n^2 d|) \right\|_{\frac{12}{11};B_\rho(x_0)}
$$
$$
\lesssim \rho^{\frac14} \left\|(|u| + |\n d|)( |\n u|+|\n^2 d|) \right\|_{\frac{12}{11};B_\rho(x_0)}
$$
and then H\"older's inequality with $11/12 = 1/4 + 1/6 + 1/2$ gives
$$
\!\!\!\!\!\!\!\!\!\!\!\!\!\!\!\|p_1\|_{\frac32;B_\rho(x_0)}^{\frac32}
\lesssim
\left(
\rho^{\frac14}
\left\|(|u| + |\n d|)^{\frac12} \right\|_{4;B_\rho(x_0)}
\left\|(|u| + |\n d|)^{\frac12} \right\|_{6;B_\rho(x_0)}
\left\|\ |\n u|+|\n^2 d|\ \right\|_{2;B_\rho(x_0)}
\right)^{\frac32}
$$
\begin{equation}\label{be}
\lesssim
\rho^{\frac38}
\left(\rho A(\rho)  \right)^{\frac38}
\left\|\ |u| + |\n d|\  \right\|^{\frac34}_{3;B_\rho(x_0)}
\left\|\ |\n u|+|\n^2 d|\ \right\|^{\frac32}_{2;B_\rho(x_0)}\, . \qquad \qquad \ \ \quad
\end{equation}
\ \\
For the following, we fix now any $r \in (0, \frac \rho 2]$, and omit the dependence on $x_0$, $t_0$ and $z_0$ in $B_r(x_0)$, $B_\rho(x_0)$, $I^*(t_0)$, $A_{z_0}$, $B_{z_0}$, $C_{z_0}$ and $D_{z_0}$ (we will retain $z_0$ in the notation for  $F_{z_0}$ to distinguish it from $F=\n f$).
\ \\\\
To first prove (\ref{n}), we note that (\ref{ba}) implies (since $r\leq \frac \rho 2$) that

$$\int_{B_r}|p_2|^{\frac32}\, dx \lesssim r^3
\rho^{-3} \left(\int_{B_\rho} (|u|^2 + |\n d|^2)\, dy \right)^{\frac34} \left(\int_{B_\rho} (|\n u|^2 + |\n^2 d|^2)\, dy \right)^{\frac34}
$$
$$
\leq r^3
\rho^{-3} \left(\rho A(\rho) \right)^{\frac34} \left(\int_{B_\rho} (|\n u|^2 + |\n^2 d|^2)\, dy \right)^{\frac34} \ \ \quad
$$
so that, integrating over $t\in I_r^*$  and using H\"older's inequality, we have
\begin{equation}\label{ad}
r^{-2}\intt{Q_r^*}|p_2|^{\frac32}\, dz \lesssim r^{-2}r^3 \rho^{{-\frac94}}  A^{\frac34}(\rho)\cdot |I^*_\rho|^{\frac 14}\left(\rho B(\rho)\right)^{\frac34}
= \frac r \rho \cdot [(AB)^{\frac34}](\rho)\ ,
\end{equation}
and that (\ref{ah}) similarly implies that
\begin{equation}\label{ab}
r^{-2} \intt{Q_r^*} |p_3|^{\frac 3 2}\, dz \lesssim r \rho^{-\frac 9 2}  \int_{I^*_r} \left(\int_{B_\rho}|p|\, dy \right)^{\frac 3 2}
\lesssim \frac r \rho \cdot D(\rho)\, .
\end{equation}
Finally, integrating (\ref{be}) over $t\in I_r^*$, H\"older  with $1=1/4+  3/4$  gives
$$
\begin{array}{rcl}
r^{-2} \|p_1\|_{\frac32;Q_r^*}^{\frac32}
&\lesssim &
r^{-2} \rho^{\frac34}
 A^{\frac38}(\rho)
\left\|\ |u| + |\n d|\  \right\|^{\frac34}_{3;Q^*_\rho}
\left\|\ |\n u|+|\n^2 d|\ \right\|^{\frac32}_{2;Q^*_\rho}
\\\\
&\lesssim &
r^{-2} \rho^{\frac34}
 A^{\frac38}(\rho)
\left(\rho^2 C(\rho)  \right)^{\frac14}
\left( \rho B(\rho) \right)^{\frac34}
=  \left( C^{\frac14}(\rho)\right) \cdot \left(\left(\frac r \rho \right)^{-2} A^{\frac38}(\rho)B^{\frac34}(\rho)\right)\ .
\end{array}
$$
Multiplying and dividing by $( r/ \rho )^{\frac \alpha 2}$ for any $\alpha \in \R$, Cauchy's inequality gives
\begin{equation}\label{ae}
r^{-2} \|p_1\|_{\frac32;Q_r^*}^{\frac32} \lesssim
\left(\frac r \rho  \right)^{\alpha }  C^{\frac12}(\rho)  + \left(\frac r \rho  \right)^{- \alpha  -4}  A^{\frac34}(\rho)B^{\frac 32}(\rho)\ .
\end{equation}
Since we want a positive power of $\gamma = r/\rho$ in the first term and a negative one on the second (because it contains $B$ which will be small), we want to take $\alpha>0$.  Choosing $\alpha = 1$ purely to make the following expression simpler, since $p = p_3 + p_2 + p_1$, we see from (\ref{ad}), (\ref{ab}) and (\ref{ae}) that
$$D(r) \lesssim \frac r \rho \cdot [D +   (AB)^{\frac34} +  C^{\frac12}](\rho)  + \left(\frac r \rho  \right)^{-5}  \left[A^{\frac34}B^{\frac 32}\right](\rho)
$$
which implies (\ref{n}) for $\gamma :=\frac r \rho \leq \frac 12$.
\\\\
To prove (\ref{p}), we note that $F_{z_0}(r) \leq  F_1(r) + F_2(r) +F_3(r)$, where we set
$$F_j(r) := \frac1{r^2} \intt{Q_r} |p_j||u|\ dz\ .$$
To estimate $F_1$ we use H\"older and (\ref{be})
to see that (in fact, for $r\leq \rho$)
$$\int_{B_r}|p_1||u|\ dx \leq \|u\|_{3;B_\rho} \|p_1\|_{\frac32;B_\rho} \qquad  \qquad  \qquad  \qquad  \qquad  \qquad  \qquad  \qquad
$$
$$
 \qquad  \qquad  \qquad \lesssim
\|u\|_{3;B_\rho}\cdot \rho^{\frac14}
\left(\rho A(\rho)  \right)^{\frac14}
\left\|\ |u| + |\n d|\  \right\|^{\frac12}_{3;B_\rho}
\left\|\ |\n u|+|\n^2 d|\ \right\|_{2;B_\rho}
$$
$$
\ \leq
\rho^{\frac12}
 A^\frac14(\rho)
\left\|\ |u| + |\n d|\  \right\|^{\frac32}_{3;B_\rho}
\left\|\ |\n u|+|\n^2 d|\ \right\|_{2;B_\rho}
$$
and hence Cauchy-Schwarz in time gives
$$\qquad \qquad \qquad  F_1(r) \lesssim r^{-2}\rho^{\frac12}
 A^\frac14(\rho)
\left\|\ |u| + |\n d|\  \right\|^{\frac32}_{3;Q_\rho^*}
\left\|\ |\n u|+|\n^2 d|\ \right\|_{2;Q_\rho^*}
$$
$$
  \lesssim
r^{-2}\rho^{\frac12}
 A^\frac14(\rho)
(\rho^2C(\rho))^{\frac12} (\rho B(\rho))^\frac12
$$
$$
 \quad  \qquad \qquad  =
\left(\left(\frac r\rho \right)^\alpha
C^{\frac12}(\rho)\right)\cdot \left(\left(\frac r\rho \right)^{-2-\alpha}   [A^\frac14 B^\frac12](\rho)\right)
$$
$$ \quad  \qquad  \qquad  \qquad \lesssim
\left(\left(\frac r\rho \right)^\alpha
C^{\frac12}(\rho)\right)^\frac43 +\left(\left(\frac r\rho \right)^{-2-\alpha}   [A^\frac14 B^\frac12](\rho)\right)^4
$$
for any $\alpha \in \R$.  Taking, say, $\alpha=\frac12$, we have
\begin{equation}\label{bg}
F_1(r) \lesssim \left(\frac r\rho \right)^\frac23
C^{\frac23}(\rho) +\left(\frac r\rho \right)^{-10}   [A B^2](\rho)\, .
\end{equation}
Now for $F_2$ note that, using (\ref{ba}), we have (since $r\leq \frac \rho 2$)
$$
\begin{array}{rcl}
\displaystyle{\int_{B_r}|p_2||u|\, dx}& \lesssim & \displaystyle{  \rho^{-2} \int_{B_\rho} (|u||\n u| + |\n d||\n^2 d|)\, dy \int_{B_r}|u|\, dx
}\\\\
 & \lesssim &   \rho^{-2} \|\ |u| + |\n d|\ \|_{2;B_\rho} \|\ |\n u|+|\n^2 d|\ \|_{2;B_\rho} (r^3)^\frac12\|u\|_{2;B_r}
\\\\
&\lesssim &
   \rho^{-2}r^\frac32 (\rho A(\rho)) \|\ |\n u|+|\n^2 d|\ \|_{2;B_\rho}
\end{array}
$$
so that integrating over $t\in I^*_r$ and using H\"older in time we have
\begin{equation}\label{bc}
F_2(r) \lesssim \frac1{r^2} \frac{r^\frac32}{\rho^2}(\rho A(\rho))(\rho B(\rho))^\frac12 (r^2)^\frac12 = \left(\frac r\rho \right)^\frac12[AB^\frac12](\rho)\ .
\end{equation}
For $F_3$, using (\ref{ah}) and H\"older, we see that
$$\frac1{r^2}\int_{B_r} |p_3||u|\ dx \leq \frac{1}{r^2\rho^3}\left(\int_{B_\rho}|p|\, dy \right)\left(\int_{B_r}|u|\, dx \right) \qquad  \qquad  \qquad  \qquad  \qquad \qquad  \qquad  \qquad
$$
$$
 \qquad  \qquad \leq \frac{1}{r^2\rho^3}\left(\int_{B_\rho}|p|^{\frac32}\, dx \right)^{\frac23}(\rho^3)^{\frac13}\left(\int_{B_r}(|u|^{\frac12})^4\, dx \right)^{\frac14}
\left(\int_{B_r}(|u|^{\frac12})^6\, dx \right)^{\frac16}(r^3)^{\frac7{12}}
$$
which gives us (setting $\gamma:=\frac r\rho$)
$$F_3(r) \lesssim \frac{1}{r^\frac14\rho^2}(rA(r))^\frac14 \left(\intt{Q_\rho^*}|p|^{\frac32}\ dx \right)^{\frac23}\left(\intt{Q_r^*}|u|^3\ dx \right)^{\frac16}(r^2)^\frac16 \qquad  \qquad  \qquad
$$
$$
\leq
\frac{1}{r^\frac14\rho^2}(rA(r))^\frac14 \left(\rho^2D(\rho) \right)^{\frac23}\left(r^2C(r) \right)^{\frac16}(r^2)^\frac16 \qquad  \qquad  \qquad  \qquad  \qquad
$$
$$
\ \leq \left(\frac{r}{\rho}\right)^\frac23 (\g^{-1}A)^\frac14(\rho) D^\frac23(\rho) (\g^{-2}C)^\frac16(\rho)
= \left(\frac{r}{\rho}\right)^\frac1{12} A^\frac14(\rho) D^\frac23(\rho) C^\frac16(\rho)
$$
by (\ref{ca}).  Hence Young's inequality implies
\begin{equation}\label{am}
F_3(r) \lesssim
\left(\frac{r}{\rho}\right)^\frac1{12} \left(A(\rho)+ D^\frac43(\rho) +C^\frac23(\rho)\right)\, .
\end{equation}
Adding  (\ref{bg}),  (\ref{bc}) and (\ref{am}) and passing to the smallest powers of   $\g=\frac r\rho \, (< 1)$ we see that
$$
F_{z_0}(r) \lesssim
\left(\frac{r}{\rho}\right)^\frac1{12} \left(A+ D^\frac43 +C^\frac23\right)(\rho)
 +\left(\frac r\rho \right)^{-10}   [A(B^\frac12+ B^2)](\rho)
$$
which implies (\ref{p}), and completes the proof of Claim \ref{clma}. \hfill $\Box$
\ \\\\\\
{\bf Proof of Claim \ref{clmalocen}:} \quad We will again omit the dependence on $z_0$ (except in $F_{z_0}$).
\\\\
To estimate $A(\tfrac \rho 2)$, we use the local energy inequality  (\ref{locenta})  with a non-negative cut-off function $\phi \in C^\infty_0(Q^*_\rho)$ which is equal to $1$ in $Q^*_{\frac \rho 2}$, with
$$|\nabla \phi| \lesssim \rho^{-1} \qquad \textrm{and} \qquad |\phi_t|, |\n^2 \phi| \lesssim \rho^{-2}\ .$$
We'll need to estimate terms which control those that appear on the right-hand side of the local energy inequality (\ref{locenta}), which we'll call $I$ - $V$ (all of which depend on $\rho$) as follows:
$$I:= \intt{Q^*_\rho} (|u|^2 + |\nabla d|^2)|\phi_t + \D \phi|\ dz
\lesssim
\rho^{-2}\|\ |u|^2 + |\n d|^2 \ \|_{\frac32;Q^*_\rho}(\rho^5)^\frac13
\qquad \qquad \qquad $$
\begin{equation}\label{cd}
\qquad \qquad \qquad \qquad \qquad \quad \lesssim
\rho^{-2}(\rho^2C(\rho))^\frac23(\rho^5)^\frac13 = \rho C^\frac23(\rho)\ .
\end{equation}
Using the assumption (\ref{divfree}) that  $\n \cdot u = 0$ weakly and indicating by $\overline{g}^\rho$ the average of a function $g$ in $B_\rho$, we have
$$
\begin{array}{rrl}
II &:= &\displaystyle{\int_{I^*_\rho} \left|\int_{B_{\rho}}(|u|^2 + |\n d|^2)u\cdot \nabla \phi\, dx \right| \, dt
}\\\\
 &= &\displaystyle{ \int_{I^*_\rho} \left|\int_{B_{\rho}}\left[(|u|^2 - \overline{|u|^2}^\rho) + (|\n d|^2 - \overline{|\n d|^2}^\rho)\right]u\cdot \nabla \phi\, dx \right| \, dt}
\end{array}$$
hence
\begin{equation}\label{ce}
II\lesssim \rho^{-1} (\rho^2 E(\rho)) = \rho E(\rho)\ .
\end{equation}
Clearly we have
\begin{equation}\label{cf}
III:= \intt{Q^*_\rho}|pu\cdot \n \phi|\ dz \lesssim \rho^{-1}(\rho^2F_{z_0}(\rho)) = \rho F_{z_0}(\rho)\ .
\end{equation}
Using the weak divergence-free condition  $\n \cdot u=0$ in (\ref{divfree}) to write (see (\ref{vecprodrule}))
$$(u\cdot \n)d = \n^T \cdot (d\otimes u)$$
(at almost every $x$) and  integrating by parts we have
$$IV:
= \int_{I^*_\rho} \left|\int_{B_{\rho}} u\otimes \n \phi : \n d \odot \n d \, dx \right| \, dt
= \int_{I^*_\rho} \left|\int_{B_{\rho}} [(u\cdot \n)d]\cdot[(\n \phi \cdot \n)d] \, dx \right| \, dt \qquad  \qquad  \qquad
$$
$$
= \int_{I^*_\rho} \left|\int_{B_{\rho}} [\n^T \cdot (d\otimes u)]\cdot[(\n \phi \cdot \n)d] \, dx \right| \, dt
= \int_{I^*_\rho} \left|-\int_{B_{\rho}}  d\otimes u : \n^T [(\n \phi \cdot \n)d] \, dx \right| \, dt\, ,
$$
and clearly
$$|\n^T [(\n \phi \cdot \n)d]| \lesssim |\n^2 \phi||\n d| + |\n \phi | |\n^2 d|\, .$$

\noindent
Therefore,  for $q\in [2,6]$ we have\footnote{Note that it is only the appearance of $\n^2 d$ in the estimate of term $IV$ which forces us to include $u$ in the definition of $G_{q,z_0}$.  Indeed, switching the roles of $u$ (which appears in $C_{z_0}$ along with $\n d$) and $\n d$ (which appears in $G_{q,z_0}$ even with $u$ omitted), one could otherwise control  term  $IV$ in precisely the same way.  If $u$ is omitted in $G_{q,z_0}$, one could still obtain the same estimate of $IV$ if one takes $q=6$, but this would dramatically weaken the statement of Theorem \ref{mainthm}.  The remainder of the proof of Theorem \ref{mainthm} does not require (but is not harmed by) the inclusion of $u$ in $G_{q,z_0}$.}
$$IV  \lesssim
  \intt{Q^*_\rho} |d||u|\left(\rho^{-2}|\n d| + \rho^{-1}|\n^2 d|\right)\ dz
 \qquad  \qquad  \qquad
   \qquad
 $$
$$
 \qquad \qquad \quad \
\begin{array}{cl}
\leq &
 \|\, |d||u|\, \|_{2;Q^*_\rho}\left( \rho^{-2}\|\n d\|_{2;Q^*_\rho}
+
\rho^{-1} \|\n^2 d\|_{2;Q^*_\rho}\right)
\\\\
\lesssim &
 \|\, |d||u|\, \|_{2;Q^*_\rho}\left( \rho^{-2}\cdot \rho^{\frac 56}\|\n d\|_{3;Q^*_\rho}
+
\rho^{-1} \|\n^2 d\|_{2;Q^*_\rho}\right)
\\\\
\leq &
 \left(\rho^3G_2(\rho)\right)^{\frac 12}\left( \rho^{-2}\cdot \rho^{\frac 56}(\rho^2 C(\rho))^{\frac 13}
+
\rho^{-1} (\rho B(\rho))^{\frac 12}\right)
\\\\
= &
 \rho\,  \left(G_2(\rho) \right)^{\frac 12}\left(  C^{\frac 13}(\rho)
+
 B^{\frac 12}(\rho)\right)
\\\\
\stackrel{(\ref{gsiginterpest})}{\leq}&
 \rho\,  \left(G_q^{\frac 2q}(\rho)C^{1-\frac 2q}(\rho) \right)^{\frac 12}\left(  C^{\frac 13}(\rho)
+
 B^{\frac 12}(\rho)\right)
\, ,
\end{array}
$$
so that
\begin{equation}\label{cg}
IV \lesssim
 \rho\left[G_q^{\frac 1q}\left(C^{\frac 56-\frac 1q}+ C^{\frac 12-\frac 1q} B^{\frac 12}\right)\right](\rho)
\, .
\end{equation}
Similarly,  for $q\in [2,6]$ we  have
\begin{equation}\label{ch}
V:=   \intt{Q^*_\rho} |d|^2|\n d|^2\phi  \ dz
\lesssim \rho^3 G_2(\rho)
\stackrel{(\ref{gsiginterpest})}{\leq} \rho^3 G_q^{\frac 2q}(\rho)C^{1-\frac 2q}(\rho)\, .
\end{equation}
Finally, using (\ref{cd}) - (\ref{ch}), the local energy inequality (\ref{locenta}) (with constant ${\bar C}$) gives
$$
\begin{array}{rcl}
({\bar C})^{-1}\tfrac \rho 2A(\tfrac \rho 2) &\lesssim & I + II  + III+ IV  + V
\\\\
& \lesssim &\rho \left[ C^\frac23 +  E+F_{z_0} + G_q^{\frac 1q}\left(C^{\frac 56-\frac 1q}+ C^{\frac 12-\frac 1q} B^{\frac 12}\right) + [\, \cdot\, ]^2 G_q^{\frac 2q}C^{1-\frac 2q} \right](\rho)
\\\\
& \lesssim &\rho \left[ C^\frac23 +  E+F_{z_0} + (1+[\, \cdot\, ]^2)G_q^{\frac{4}{6-q}}+  (G_q^{\frac{2}{6-q}}+C^{\frac 13}) B^{\frac 12} \right](\rho)
\end{array}
$$
as long as $2\leq q<6$, as in that case we have
$$
G_q^{\frac 1q}C^{\frac 56-\frac 1q}=(G_q^{\frac{4}{6-q}})^{\frac{6-q}{4q}}(C^{\frac 23})^{\frac{5q-6}{4q}}
\leq \left(\frac{6-q}{4q}\right) G_q^{\frac{4}{6-q}}+\left(\frac{5q-6}{4q}\right) C^{\frac 23}
\leq \tfrac{3}{4} G_q^{\frac{4}{6-q}}+\tfrac{5}{4} C^{\frac 23}
\, ,
$$
$$
G_q^{\frac 1q}C^{\frac 12-\frac 1q}=(G_q^{\frac{2}{6-q}})^{\frac{6-q}{2q}}(C^{\frac 13})^{\frac{3q-6}{2q}}
\leq \left(\frac{6-q}{2q}\right)G_q^{\frac{2}{6-q}}+\left(\frac{3q-6}{2q}\right)C^{\frac 13}
\leq \tfrac{3}{2}G_q^{\frac{2}{6-q}}+\tfrac{3}{2}C^{\frac 13}
$$
and
$$
G_q^{\frac 2q}C^{1-\frac 2q}= (G_q^{\frac{4}{6-q}})^{\frac{6-q}{2q}}(C^{\frac 23})^{\frac{3q-6}{2q}}
\leq \left(\frac{6-q}{2q}\right)G_q^{\frac{4}{6-q}} + \left(\frac{3q-6}{2q}\right)C^{\frac 23}
\leq \tfrac{3}{2}G_q^{\frac{4}{6-q}} + \tfrac{3}{2}C^{\frac 23}\, .
$$
This implies (\ref{q}) and proves Claim \ref{clmalocen}. \hfill $\Box$


\begin{thebibliography}{CKN82}

\bibitem[CKN82]{caf}
L.~Caffarelli, R.~Kohn, and L.~Nirenberg.
\newblock Partial regularity of suitable weak solutions of the
  {N}avier-{S}tokes equations.
\newblock {\em Comm. Pure Appl. Math.}, 35(6):771--831, 1982.

\bibitem[DHW19]{duhuwang}
H. Du, X. Hu, and C. Wang.
\newblock Suitable weak solutions for the co-rotational {B}eris-{E}dwards system in dimension three.
\newblock {\em arXiv:1905.08440} (preprint), 2019.

\bibitem[GT01]{gilbarg}
D.~Gilbarg and N.~S. Trudinger.
\newblock {\em Elliptic Partial Differential Equations of Second Order}.
\newblock Springer-Verlag, Berlin, 2001.

\bibitem[LS99]{ladyser}
O.~A. Ladyzhenskaya and G.~A. Seregin.
\newblock On partial regularity of suitable weak solutions to the
  three-dimensional {N}avier-{S}tokes equations.
\newblock {\em J. Math. Fluid Mech.}, 1(4):356--387, 1999.

\bibitem[Ler34]{leray}
J.~Leray.
\newblock Sur le mouvement d'un liquide visqueux emplissant l'espace.
\newblock {\em Acta Math.}, 63:193--248, 1934.

\bibitem[Lin98]{lin}
Fanghua Lin.
\newblock A new proof of the {C}affarelli-{K}ohn-{N}irenberg theorem.
\newblock {\em Comm. Pure Appl. Math.}, 51(3):241--257, 1998.


\bibitem[LL95]{linliu95}
Fang-Hua Lin and Chun Liu.
\newblock Nonparabolic dissipative systems modeling the flow of liquid
  crystals.
\newblock {\em Comm. Pure Appl. Math.}, 48(5):501--537, 1995.


\bibitem[LL96]{linliu}
Fang-Hua Lin and Chun Liu.
\newblock Partial regularity of the dynamic system modeling the flow of liquid
  crystals.
\newblock {\em Discrete Contin. Dynam. Systems}, 2(1):1--22, 1996.

\bibitem[LW14]{linwang14}
Fanghua Lin and Changyou Wang.
\newblock Recent developments of analysis for hydrodynamic flow of nematic
  liquid crystals.
\newblock {\em Philos. Trans. R. Soc. Lond. Ser. A Math. Phys. Eng. Sci.},
  372(2029):20130361, 18, 2014.

\bibitem[Sch77]{scheffer77}
Vladimir Scheffer.
\newblock Hausdorff measure and the {N}avier-{S}tokes equations.
\newblock {\em Comm. Math. Phys.}, 55(2):97--112, 1977.

\bibitem[Sch80]{scheffer80}
Vladimir Scheffer.
\newblock The {N}avier-{S}tokes equations on a bounded domain.
\newblock {\em Comm. Math. Phys.}, 73(1):1--42, 1980.


\bibitem[Sch85]{scheffer3}
Vladimir Scheffer.
\newblock A solution to the {N}avier-{S}tokes inequality with an internal
  singularity.
\newblock {\em Comm. Math. Phys.}, 101(1):47--85, 1985.

\bibitem[Vas07]{vasseur}
Alexis~F. Vasseur.
\newblock A new proof of partial regularity of solutions to {N}avier-{S}tokes
  equations.
\newblock {\em NoDEA Nonlinear Differential Equations Appl.}, 14(5-6):753--785,
  2007.


\bibitem[WZ77]{wheeden}
Richard~L. Wheeden and Antoni Zygmund.
\newblock {\em Measure and integral}.
\newblock Marcel Dekker, Inc., New York-Basel, 1977.
\newblock An introduction to real analysis, Pure and Applied Mathematics, Vol.
  43.


\end{thebibliography}


\end{document}